\documentclass[12pt]{article}
\usepackage{amssymb,amsmath,amsfonts,graphicx}
\newtheorem{theorem}{Theorem}[section]
\newtheorem{definition}[theorem]{Definition}
\newtheorem{corollary}[theorem]{Corollary}
\newtheorem{proposition}[theorem]{Proposition}
\newtheorem{remark}[theorem]{Remark}
\newtheorem{lemma}[theorem]{Lemma}

\newtheorem{problem}[theorem]{Problem}
\newcommand {\Kc}      {{\mathcal K}}
\newcommand {\Mc}      {{\mathcal M}}

\newcommand {\Wc}      {{\mathcal W}}

\newcommand {\Ac}      {{\mathcal A}}

\newcommand {\Cc}      {{\mathcal C}}

\newcommand {\Qc}      {{\mathcal Q}}

\newcommand {\Oc}      {{\mathcal O}}
\newcommand {\Sc}      {{\mathcal S}}
\newcommand {\WP}      {W^1_p(\RN)}
\newcommand {\WPO}     {W^1_p(\Omega)}
\newcommand {\LOPO}    {L^1_p(\Omega)}
\newcommand {\LOP}     {L^1_p(\RN)}

\newcommand {\TW}      {{\mathbb W}(\Omega)}

\newcommand {\tQ}      {\widetilde{Q}}
\newcommand {\tF}      {\widetilde{F}}

\newcommand {\tf}      {\tilde{f}}

\newcommand {\tx}      {\tilde{x}}
\newcommand {\ty}      {\tilde{y}}
\newcommand {\R}       {{\bf R}}
\newcommand {\N}       {{\bf N}}
\newcommand {\RN}      {\R^n}
\newcommand {\ve}      {\varepsilon}
\newcommand {\tom}     {\tilde{\omega}}
\newcommand {\emp}     {\emptyset}
\newcommand {\cw}      {\curlywedge}
\newcommand {\clO}     {\overline{\Omega}}
\newcommand {\bc}      {\bar{c}}
\newcommand {\intl}    {\int\limits}
\newcommand {\DO}      {\partial\Omega}
\newcommand {\DOEP}    {\Oc_\ve(\DO)}

\newcommand {\dt}      {\rho_{\alpha,\Omega}}
\newcommand {\dtc}     {\rho_{\alpha,\clO}}
\newcommand {\dtcb}    {\rho_{\beta,\clO}}
\newcommand {\db}      {\rho_{\beta,\Omega}}
\newcommand {\dao}     {d_{\alpha,\Omega}}
\newcommand {\dbo}     {d_{\beta,\Omega}}
\newcommand {\dco}     {d_{\alpha,\clO}\,}
\newcommand {\dcb}     {d_{\beta,\clO}\,}

\newcommand {\CSQ}     {\Cc[\Omega,\dt]}
\newcommand {\OA}      {\Omega^{*,\,\alpha}}
\newcommand {\OB}      {\Omega^{*,\,\beta}}
\newcommand {\DOA}     {(\DO)_\alpha}
\newcommand {\DOB}     {(\DO)_\beta}
\newcommand {\EN}      {\|\cdot\|}
\newcommand {\IPO}     {I_\alpha(\Omega)}
\newcommand {\pn}      {(p-n)/(p-1)}
\newcommand {\qn}      {(q-n)/(q-1)}

\newcommand {\baf}     {\bar{f}_\alpha}
\newcommand {\tra}     {\trs_{\DOA}\!}
\newcommand {\D}[1]    {\dist(#1,\DO)}
\newcommand {\eo}      {\ell(\omega)}
\newcommand {\oqa}     {\omega_{Q,\alpha}}
\newcommand {\oqp}     {\omega_{Q',\alpha}}
\newcommand {\WQ}      {\widetilde{\Qc}}
\newcommand {\sa}      {\stackrel{\alpha}{\sim}}
\newcommand {\sbs}     {\stackrel{\beta}{\sim}}
\newcommand {\dlbeta}  {\delta_{\beta,\clO}\,}

\newcommand {\supp}    {\operatorname{supp}}
\newcommand {\Ex}      {\operatorname{Ext}}
\newcommand {\diam}    {\operatorname{diam}}
\newcommand {\dist}    {\operatorname{dist}}

\newcommand {\Conv}    {\operatorname{Conv}}
\newcommand {\Lip}     {\operatorname{Lip}}
\newcommand {\len}     {\operatorname{len}}
\newcommand {\trs}     {\operatorname{tr}}
\newcommand {\lng}     {\operatorname{length}}
\newcommand {\esssup}  {\operatorname{ess\,\sup}}
\newcommand {\bx}      {\hspace{10mm}$\Box$}
\newcommand {\BX}      {\hspace{10mm}\Box}
\newcommand {\nn}      {\nonumber}
\newcommand {\rf}[1]    {(\ref{#1})}   
\newcommand {\reff}[1]  {\ref{#1}}     
\newcommand{\lbl}[1]      {\label{#1}}       
\newcommand{\be}          {\begin{eqnarray}}
\newcommand{\bel}[1]      {\begin{eqnarray} \label{#1}}
\newcommand{\ee}           {\end{eqnarray}}

\begin{document}
\parindent 1em
\parskip 0mm
\medskip
\centerline{\large{\bf On the boundary values of Sobolev $W^1_p$-functions}}
\vspace*{10mm} \centerline{By  {\it Pavel Shvartsman}}
\vspace*{3mm}
\centerline {\it Department of Mathematics, Technion -
Israel Institute of Technology,}
\centerline{\it 32000 Haifa, Israel}\vspace*{3 mm}
\centerline{\it e-mail: pshv@tx.technion.ac.il}
\vspace*{12 mm}
\renewcommand{\thefootnote}{ }
\footnotetext[1]{{\it\hspace{-6mm}Math Subject
Classification} 46E35}
\begin{abstract} 
For each $p>n$ we use local oscillations to give intrinsic characterizations of the trace of the Sobolev space $\WPO$ to the boundary of an \emph{arbitrary} domain $\Omega\subset\RN$.
\end{abstract}
\section*{{\normalsize 1. Introduction}}
\setcounter{section}{1}
\setcounter{theorem}{0}
\setcounter{equation}{0}
\par {\bf 1.1. The trace problem for the Sobolev space
$\WPO$.}
\par Let $\Omega$ be a domain in $\RN$. We recall that, for each  $p\in [1,\infty]$, the Sobolev space $\WPO$ consists of all (equivalence classes of) real valued functions $f\in L_{p}(\Omega)$ whose first order distributional partial derivatives on $\Omega$ belong to $L_{p}(\Omega)$. $W_{p}^{1}(\Omega)$ is normed by
$$ \|f\|_{\WPO}:=\|f\|_{L_{p}(\Omega)}+ \|\nabla
f\|_{L_{p}(\Omega )}. $$
\par By $\LOPO$ we denote the corresponding homogeneous
Sobolev space, defined by the finiteness of the seminorm
$$
\|f\|_{\LOPO}:=\|\nabla f\|_{L_{p}(\Omega )}.
$$
\par In this paper we study the following trace
problem.
\begin{problem}\label{PR1} {\em Given $p\in[1,\infty]$  and an arbitrary function $f:\DO\to\R$, find a necessary and sufficient condition for $f$ to be the trace to $\DO$
of a function $F\in\WPO$. \vspace*{2mm}}
\end{problem}
\par This problem is of great interest, mainly
due to its important applications to boundary-value problems in partial differential equations where it is essential to be able to characterize the functions defined on $\DO$, which appear as traces to $\DO$ of Sobolev functions.
\par There is an extensive literature devoted to the theory of boundary traces in Sobolev spaces. Among the multitude of results we mention the monographs of Grisvard \cite{Gr}, Lions and Magenes \cite{LM}, Maz'ya and Poborchi \cite{M,MP}, and the papers by Gagliardo \cite{Ga}, Nikol'skii \cite{N}, Besov \cite{Be}, Aronszajn and  Szeptycki \cite{AS}, Yakovlev \cite{Y1,Y2}, Jonsson and Wallin \cite{J1,J2,JW}, Maz'ya and Poborchi \cite{MP1,MP2,MP3,P1}, Maz'ya, Poborchi and Netrusov \cite{MPN}, Vasil'chik \cite{V}; we refer the reader to these works and references therein for numerous results and techniques concerning this topic. In these monographs and papers Problem \reff{PR1} is investigated and solved for various families of smooth, Lipschitz and non-Lipschitz domains in $\RN$ with different types of singularities on the boundary.
\par In this paper we characterize the traces to the
boundary of Sobolev $\WPO$-functions whenever $p>n$ and
$\Omega$ is an {\it arbitrary} domain in $\RN$.
\par The first challenge that we face in the study of Problem \reff{PR1} for an arbitrary domain $\Omega\subset\RN$ is the need to find a ``natural" definition of the trace of a Sobolev $\WPO$-function to the boundary of the domain which is compatible with the structure of Sobolev functions on the domain $\Omega$. More specifically, we want to choose this definition in such a way that every $f\in\WPO$ possesses a well-defined ``trace'' to $\DO$ which in a certain sense characterizes the behavior  of the function $f$ near the boundary.
\par We recall that, when $p>n$, it follows from the Sobolev embedding theorem that every function $f\in\LOPO$
coincides almost everywhere with a function satisfying the {\it local} H\"{o}lder condition of order $\alpha:=1-\frac{n}{p}$. I.e., after possibly modifying $f$ on a set of Lebesgue measure zero, we have, for every cube
$Q\subset \Omega$, that
\bel{S-Lip-OM} |f(x)-f(y)|\le C(n,p)\|f\|_{\LOPO}
\|x-y\|^{1-\frac{n}{p}}~~~\text{ for every}~~x,y\in Q. \ee
\par This fact enables us {\it to identify each element of $\LOPO$ with its unique continuous representative}. Thus we are able to restrict our attention to the case of {\it
continuous} Sobolev functions defined on $\Omega$.
\par Since we deal only with continuous Sobolev functions defined on $\Omega$, given $f\in \LOPO$ it would at first seem quite natural to try to define the ``boundary values'' of $f$ on $\DO$ to be the {\it continuous} extension of $f$ from $\Omega $ to $\DO $. In other words, we could try to extend the domain of definition of $f$ to the closure of $\Omega $, the set $\clO$, by letting
\bel{DEF-CONT}
{\bar{f}}(y):=\lim_{x\to y,x\in\Omega}f(x),
~~~~\text{for each }y\in\clO.
\ee
\par This indeed is the natural definition to use for certain classes of domains in $\RN$ such as Lipschitz domains or $(\ve,\delta)$-domains, see Jones \cite{Jn}.
But in general it does not work. For an obvious example showing this, consider the planar domain which is a ``slit square'' $\Omega =(-1,1)^{2}\backslash J$, where $J$ is the line segment $[(-1/2,0),(1/2,0)]$. The reader can easily
construct a $C^{\infty}$-function $f\in \WPO$ which equals
zero on the upper ``semi-square"
$\{x=(x_1,x_2)\in[-1/4,1/4]^2:~x_2>0\}$ and takes the value $1$ on the lower ``semi-square"
$\{x=(x_1,x_2)\in[-1/4,1/4]^2:~x_2<0\}$. Clearly,
\rf{DEF-CONT} cannot provide a well defined function $\bar{f}$ on the segment $[(-1/4,0),(1/4,0)]$.
The obvious reason for the existence of such kinds of
counterexamples is the fact that the continuity of a
$\WPO$-function does not imply its {\it uniform} continuity on $\Omega$.
\par In order to define a notion of ``trace to the boundary" which will work for \emph{all} domains $\Omega$ we have to adopt a somewhat different approach. Its point of departure is an important property of Sobolev functions which will be recalled in more detail below (see definition \rf{DEF-TDA} and inequality \rf{UC-DA}), namely that
{\it every $f\in\WPO$, $p>n,$ is uniformly continuous with respect to a certain intrinsic metric $\dt$ defined on $\Omega$.}
\par This property motivates us to define the {\it completion} of $\Omega$ with respect to this intrinsic metric. We can then define the ``trace to the boundary" of each function $f\in\WPO$ by first extending $f$ by continuity (with respect to $\dt$) to a continuous function $\tf$ defined on this completion, and then taking the
restriction of $\tf$ to the appropriately defined {\it boundary} of this completion.
\par As is of course to be expected, in all cases where the definition \rf{DEF-CONT} is applicable, these new notions of boundary and trace coincide with the ``classical" ones.
\par We begin our formal development of this approach in the next subsection.
\bigskip
\par {\bf 1.2. Subhyperbolic metrics in $\Omega$ and their Cauchy completions.} Following the terminology of Buckley and Stanoyevitch \cite{BSt3}, given $\alpha \in [0,1]$ and a rectifiable curve $\gamma \subset \Omega $, we define the {\it subhyperbolic length} of $\gamma $ by the line integral
\bel{L-SH} \len_{\alpha ,\Omega}(\gamma) :=\intl_{\gamma }
\D{z}^{\alpha-1}\,ds(z). \ee
(Here $\D{z}$ denotes the usual Euclidean distance from the point $z$ to the boundary of $\Omega$, and $ds$ denotes usual arc length.)
\par Then we let $d_{\alpha ,\Omega }$ denote the corresponding {\it subhyperbolic metric} on $\Omega$ given, for each $x,y\in \Omega $, by
\bel{DEF-D} d_{\alpha ,\Omega }(x,y):=\inf_{\gamma }
\len_{\alpha ,\Omega }(\gamma ) \ee
where the infimum is taken over all rectifiable curves
$\gamma \subset \Omega $ joining $x$ to $y$.
\par The metric $d_{\alpha ,\Omega }$ was introduced and
studied by Gehring and Martio in \cite{GM}. See also
\cite{AHHL,L,BKos} for various further results using this
metric.  Note also that $\len_{0,\Omega }$ and $d_{0,\Omega }$ are the well-known {\it quasihyperbolic length} and {\it quasihyperbolic distance}, and $d_{1,\Omega }$ is the {\it inner (or geodesic) metric} on $\Omega $.
\par The subhyperbolic metric $d_{\alpha ,\Omega }$ with
$\alpha =(p-n)/(p-1)$ arises naturally in the study of
Sobolev $W_{p}^{1}(\Omega )$-functions for $p>n$. In
particular, Buckley and Stanoyevitch \cite{BSt2} proved
that the local H\"{o}lder condition \rf{S-Lip-OM} is
equivalent to the following H\"{o}lder-type condition: for
every $x,y\in \Omega$
\bel{H-LOC} |f(x)-f(y)|\le C(n,p)\| f\|_{\LOPO}\{d_{\alpha
,\Omega }(x,y)^{1-\frac{1}{p}}+\| x-y\| ^{1-\frac{n}{p}}\}
\ee
provided $f\in\LOPO$ and $\alpha =(p-n)/(p-1)$.
\par In view of this result it is convenient for us to introduce a new metric $\dt$ on $\Omega$ for each $\alpha \in (0,1]$, by simply putting
\bel{DEF-TDA} \dt(x,y):=\dao(x,y) +\|x-y\| ^{\alpha}, \ee
for each $x,y\in\Omega$. Then \rf{H-LOC} can be rewritten in the following form:
\bel{UC-DA} |f(x)-f(y)|\le C(n,p)\|f\|_{\LOPO}\,\dt(x,y)^{1-\frac{1}{p}}, ~~~x,y\in\Omega, \ee
where $\alpha=(p-n)/(p-1)$. Thus every function $f\in\LOPO$ is {\it uniformly continuous with respect to the metric $\dt$}.
\par This observation immediately implies the following
important fact for each $p>n$ and $\alpha=(p-n)/(p-1)$:
\par \smallskip{\it Every Sobolev function $f\in\LOPO$ admits a unique continuous extension from the metric space
$(\Omega,\dt)$ to its Cauchy completion.}
\medskip
\par Let us now recall several standard facts concerning Cauchy completions and fix the notation that we will use here for the particular case of the Cauchy completion of $(\Omega,\dt)$.
\par Throughout this paper we will use the notation  $(x_i)$ for a sequence $\{x_i\in\Omega:i=1,2,...~\}$ of points in $\Omega$. Let $\CSQ$ be the family of all Cauchy sequences in $\Omega$ with respect to the metric $\dt$:
\bel{D-CSQ}
\CSQ:=\{(x_i):~x_i\in\Omega,~ \lim_{i,j\to\infty}\dt(x_i,x_j)=0\}.
\ee
Observe that, by definition \rf{DEF-TDA}, the set $\CSQ$ consists of sequences $(x_i)\subset\Omega$ which converge in $\clO$ (in the Euclidean norm) and are fundamental with respect to the metric $\dao$.
\par By $``\sa"$ we denote the standard equivalence relation on $\CSQ$,
\bel{ESIM}
(x_i)\sa (y_i)~\Leftrightarrow  \lim_{i\to\infty}\dt(x_i,y_i)=0.
\ee
For each sequence $(x_i)\in\CSQ$, we use the notation
\bel{KE}
[(x_i)]_\alpha~~\text{{\it is the equivalence class of }}~~(x_i)~~\text{{\it with respect to}}~~``\sa".
\ee
\par Let
\bel{CMP-OM}
\OA:=\{[(x_i)]_\alpha:~(x_i)\in\CSQ\}
\ee
be the set of all equivalence classes with respect to $\sa$. We let $\dtc$ denote the standard metric on $\OA$ defined by the formula
\bel{DR-S}
\dtc([(x_i)]_\alpha,[(y_i)]_\alpha)
:=\lim_{i\to\infty}\dt(x_i,y_i),
~~~~[(x_i)]_\alpha,[(y_i)]_\alpha\in\OA.
\ee
\par As usual we identify every point $x\in\Omega$ with the equivalence class $\hat{x}=[(x,x,...)]_\alpha$ of the constant sequence. This identification enables us to consider the domain $\Omega$ as a subset of $\OA$. Observe that
$$
\dtc\mid_{\Omega\times\Omega}=\dt\,,
$$
i.e., the mapping $\Omega\ni x\mapsto \hat{x}\in\OA$ is an isometry.
\begin{remark}\lbl{R-CPM} {\em Since
$$
\|x-y\|^\alpha\le\dt(x,y)=\dao(x,y) +\|x-y\| ^{\alpha},~~~x,y\in\Omega,
$$
every Cauchy sequence $(x_i)\in\CSQ$ is also a Cauchy sequence with respect to the Euclidean distance. Consequently it converges to a point in $\clO$. Moreover, all sequences from any given equivalence class $\omega=[(x_i)]_\alpha\in\OA$ converge (in $\EN$) to the same point. We denote the common (Euclidean) limit point of all these sequences by $\eo$; thus
$$
y_i\stackrel{\EN}{\longrightarrow} \eo~~~
\text{as}~~ i\to\infty ~~~\text{for every sequence} ~~~(y_i)\in\omega.
$$
\par Now \rf{DEF-TDA} and \rf{DR-S} imply the following formula: for every $\omega_1=[(x_i)]_\alpha,$ $\omega_2=[(y_i)]_\alpha\in \OA$ we have
\bel{CDR-L}
\dtc(\omega_1,\omega_2)
=\lim_{i\to\infty}\dao(x_i,y_i)+
\|\ell(\omega_1)-\ell(\omega_2)\|^{\alpha}.
\ee
}
\end{remark}
\begin{remark}\lbl{OPEN} {\em As we shall see below,
$
\dt(u,v)\sim\|u-v\|^\alpha
$
provided $u,v$ belong to a sufficiently small neighborhood of a point $x\in\Omega$. This shows that the metric $\dt$ and the Euclidean metric determine the same local topology on $\Omega$. In particular, this implies that $\Omega$ is an {\it open subset} of $\OA$ (in the $\dtc$\,\,-topology).}
\end{remark}
\smallskip
\par We are now ready to define a kind of ``boundary" of $\Omega$, which is the appropriate replacement of the usual boundary for our purposes here, and will be one of the main objects to be studied in this paper.
\begin{definition}\lbl{B-DA} {\em Let $\alpha\in(0,1]$ and let $\Omega$ be a domain in $\RN$. We let $\DOA$ denote the boundary of $\Omega$ (as an open subset of $\OA$) in the topology of the metric space $(\OA,\dtc)$. We call $\DOA$ the {\it $\alpha$-boundary} of the domain $\Omega$.}
\end{definition}
\par We observe that, by Remark \reff{OPEN},
\bel{B-A}
\DOA=\OA\setminus\Omega.
\ee
Thus $\DOA$ consists of the {\it new elements} which appear as a result of taking the completion of $\Omega$ with respect to the metric $\dt$.
\par By Remark \reff{R-CPM},
$$
 \eo\in\DO~~~\text{for each}~~~ \omega\in\DOA.
$$
This means that every element of the $\alpha$-boundary can be identified with a point $x\in\DO$ and an equivalence class $[(x_i)]_\alpha$ of Cauchy sequences (with respect to the metric $\dt$) which converge to $x$ in the Euclidean norm.
\par However, as of course is to be expected from the preceding discussion, in general the set $\DOA$ will not be in one to one correspondence with $\DO$ because there
may be points $x\in\DO$ which ``split" into a family of elements $\omega\in\DOA$ all of which satisfy $\eo=x$. Such families may even be infinite. Every $\omega$ in such a family can be thought of as a certain ``approach" to the point $x$ by elements of $\Omega$ whose $\dt$-distance to $x$ tends to $0$. For example, in Fig. 1 below the point $z_1$ splits into $6$ different elements of $\DOA$, while the points $z_2,z_3$ and $z_4$ split, respectively, into $4,3$ and $2$ such elements.
\medskip
\begin{figure}[h]
\center{\includegraphics[scale=0.7]{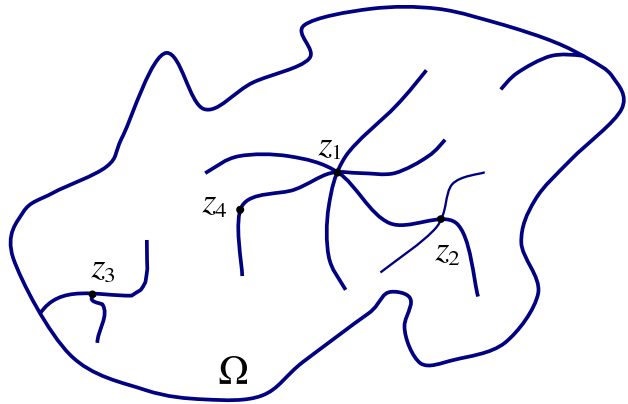}}
\caption{A domain $\Omega$ with agglutinated parts of the boundary.}
\end{figure}
\par We will use the terminology \emph{agglutinated point}
for points (like $z_i,\,i=1,2,3,$ in Fig. 1) which ``split" into multiple elements of the $\alpha$-boundary. Formally, a point $z\in\DO$ will be called an $\alpha$-agglutinated point of $\DO$ if there exist at least two different equivalence classes $\omega_1,\omega_2\in\OA,\,\omega_1\ne \omega_2,$ such that $z=\ell(\omega_1)=\ell(\omega_2)$.
We refer to the set of all $\alpha$-agglutinated points as the $\alpha$-agglutinated part of the boundary.
\par The reader may care to think of the $\alpha$-boundary of $\Omega$ as a kind of ``bundle" of the regular boundary $\DO$ under cutting of $\DO$ along its ``agglutinated" parts, see Fig. 2.
\medskip
\begin{figure}[h]
\center{\includegraphics[scale=1.1]{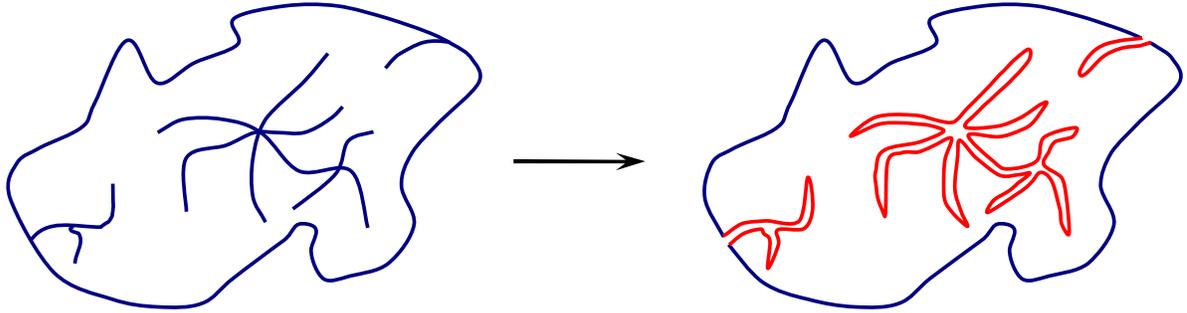}}
\caption{Cutting the domain $\Omega$ along the agglutinated parts of the boundary}
\end{figure}
\begin{remark}\lbl{R-IP} {\em In general, when we associate elements of $\DOA$ to elements of $\DO$ in the way described above, we can also lose a part of $\DO$. I.e., there may exist points in $\DO$ which do not arise as $\eo$ for any $\omega\in\DOA$. We refer to the set of all such points as the {\it $\alpha$-inaccessible part of $\DO$}. More formally, we set
$$
\IPO= \{x\in\DO:~\nexists~~\omega\in\OA
~\text{such that}~x=\ell(\omega)\}.
$$
\par Thus $\IPO$ is the set of all points $x\in\DO$ such that every sequence $(x_i)$ in $\Omega$ which converges to $x$ in $\EN$-norm, is not a Cauchy sequence with respect to the metric $\dt$. Roughly speaking $\IPO$ consists of all points $x\in\DO$ for which  $\dtc(x,\Omega)=+\infty$.
\par We call the set $\DO\setminus\IPO$ {\it the $\alpha$-accessible part of $\DO$}.
\smallskip
\par It may be helpful to give an explicit example of a domain $\Omega$ which has a non-empty inaccessible part. Figure 3 shows such a domain, which contains an infinite sequence of rectangular portions with slits. We shall choose $\alpha=1$ so that we are dealing with the geodesic metric. In this picture the line segment $[A,B]$ is a part of the boundary of $\Omega$. Clearly, for every $x\in[A,B],y\in\Omega,$ and every sequence $(x_i)$ in $\Omega$ such that $\|x_i-x\|\to 0$, the intrinsic distance $d_{1,\Omega}(x_i,y)\to\infty$ as $i\to\infty$. (Of course, every such a sequence $(x_i)$ in $\Omega$ is not a Cauchy sequence with respect to the geodesic distance in $\Omega$.) Thus $[A,B]$ is the $1$-inaccessible part of $\DO$.
\smallskip
\begin{figure}[h]
\center{\includegraphics[scale=1.1]{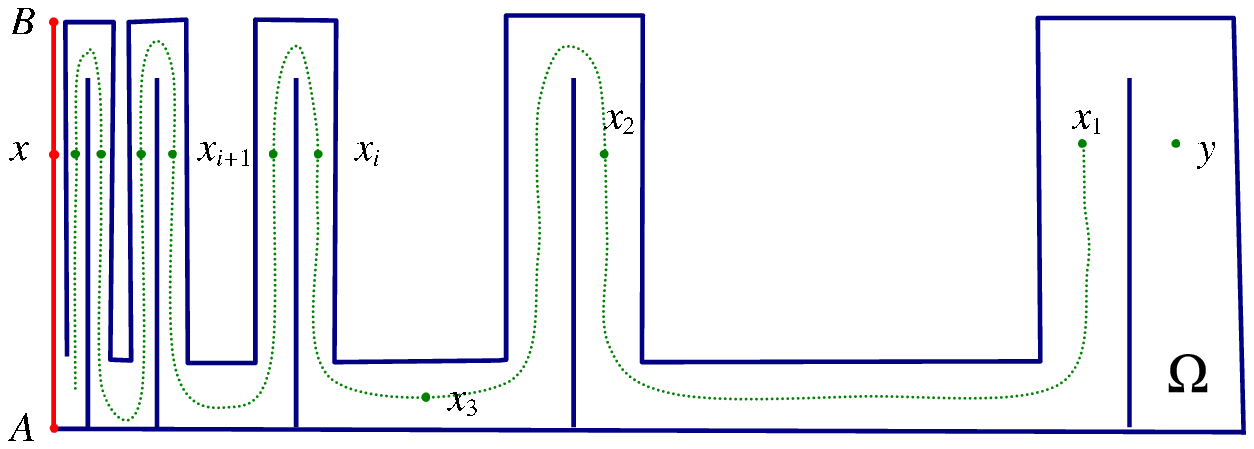}}
\caption{A domain $\Omega$ with non-empty set of $1$-inaccessible points of $\DO$.}
\end{figure}
\par The need to deal (or not deal) appropriately with parts of the boundary of some domain $\Omega$ which are
not accessible with respect to an intrinsic metric on
$\Omega$ arises naturally in the study of boundary values of Sobolev functions.
\par For instance, consider the space $L^1_\infty(\Omega)$ which coincides with the space $\Lip(\Omega,d_{1,\Omega})$ of functions on $\Omega$ satisfying a Lipschitz condition with respect to the geodesic metric $d_{1,\Omega}$. This metric can be extended ``by continuity" to the boundary of $\Omega$. We denote this extension of $d_{1,\Omega}$ by $d_{1,\clO}$. In this case the inaccessible part of $\DO$, the set  $I_1(\Omega)$, consists of points $x\in\DO$ such that the geodesic distance from $x$ to $\Omega$, i.e., $d_{1,\clO}(x,\Omega)$, equals $+\infty$. Since $d_{1,\clO}(I_1(\Omega),\Omega)=+\infty$, there are no points in $\Omega$ which are {\it close to $I_1(\Omega)$}. Consequently, the notion of the trace of a Sobolev $L^1_\infty(\Omega)$-function to $I_1(\Omega)$ is meaningless. Conversely, for the same reason, for every function $f:\DO\to\R$ its values on $I_1(\Omega)$ do not have any influence on whether $f$ can be extended to a Sobolev $L^1_\infty(\Omega)$-function on all of the domain $\Omega$. (Note that the trace of the space $L^1_\infty(\Omega)$ to the accessible part of the boundary coincides with the space $\Lip(\DO\setminus I_1(\Omega),d_{1,\clO})$.)
\par These observations motivate our approach to the notion of the boundary values of Sobolev functions on $\Omega$ as the restriction to the accessible part of $\DO$ rather than the restriction to all of the boundary of $\Omega$.}
\end{remark}
\medskip
\par {\bf 1.3. The trace to the Sobolev boundary of a domain.} We are now in a position to define the trace of $\LOPO$ to the $\alpha$-boundary of $\Omega$ whenever $\alpha=\pn$.
\par As already mentioned above, by inequality \rf{UC-DA}, every function $f\in\LOPO$ is {\it uniformly continuous} with respect to $\dt$. Since $\Omega$ is a dense subset of $\OA$ (in $\dt$-metric), there exists {\it a (unique) continuous extension $\baf$ of $f$ from $\Omega$ to $\OA$.}
\par From here onwards we will find it convenient to refer to the set $\DOA$ introduced in Definition \ref{B-DA}, when $\alpha=\pn$, as {\it the Sobolev $W^1_p$-boundary of $\Omega$ in $\RN$.}
\par We let $\tra f$  denote the restriction of $\baf$ to $\DOA$, i.e.,
\bel{DTR}
\tra f:=\baf|_{\DOA}.
\ee
We refer to the function $\tra f$ as {\it the trace of $f$ to the Sobolev $W^1_p$-boundary of\, $\Omega$.}
\medskip
\par More specifically,  $\tra f$ is a function on $\DOA$ defined as follows: Let $\omega\in\OA$ be an equivalence class and let $(y_i)\in\omega$ be an arbitrary sequence. Then
\bel{TR-S}
\tra f(\omega):=\lim_{i\to\infty} f(y_i).
\ee
Since $f\in\LOPO$ is uniformly continuous with respect to $\dt$, the trace $\tra f$ is well defined and does not depend on the choice of the sequence $(y_i)\in\omega$ in \rf{TR-S}.
\par An equivalent definition of the trace $\tra f$ is given by the formula:
\bel{TR-LM}
\tra f(\omega):=\lim\{f(x):~\dtc(x,\omega)\to 0,~ x\in\Omega\}.
\ee
\par For every domain $\Omega $ and for each $p>n$, we can define the Banach space $\tra(\WPO)$ of all traces to the Sobolev boundary of $\Omega$:
$$
\tra(\WPO)=\{f:\exists~\text{a continuous}~F\in\WPO~\text{such that}~\tra F=f\}.
$$
We equip the space $\tra(\WPO)$ with the norm
$$
\|f\|_{\tra(\WPO)}=\inf\{\|F\|_{\WPO}:~F\in\WPO~~\text{and continuous},~ \tra F=f\}.
$$
We define the space  $\tra(\LOPO)$ in an analogous way.
\par We now turn to a formulation of the main results of the paper. At this stage, for simplicity, we will present these results for a domain $\Omega$ in $\RN$ whose boundary does not contain ``agglutinated" parts. In other words, we assume for the moment that $\DO$ does not ``split" under Cauchy completion of $\Omega$ with respect to $\dt$. This simplification will allow us to interpret the trace to the boundary as a function defined on (the accessible part of) $\DO$ rather than a function defined on a certain set of equivalence classes.
\medskip
\begin{definition}\lbl{A-C} {\em We say that a domain $\Omega\subset\RN$ satisfies the condition $(A_\alpha)$ if for every ``accessible" point $x\in\DO$, i.e., for every $x\in\DO\setminus\IPO$, (see Remark \reff{R-IP}), there exists a {\it unique} equivalence class $\omega\in\OA$ such that $\eo=x$.}
\end{definition}
\medskip
\par Thus for a domain $\Omega$ satisfying the condition $(A_\alpha)$ we may identify the set $\DOA$ with $\DO\setminus\IPO$. Also we may simplify several main definitions and notions introduced above. In particular, we can introduce a metric on the set $\clO\setminus \IPO$ of all accessible points by letting
\bel {D-AO}
\dco(x,y):=\liminf\{\len_{\alpha,\Omega}(\tx,\ty):~
\tx\stackrel{\EN}{\longrightarrow} x,\ty\stackrel{\EN}{\longrightarrow} y,~\tx,\ty\in\Omega\},
\ee
see \rf{L-SH}. Here $x,y\in\clO\setminus \IPO$. Clearly, $\dco=\dao$ on $\Omega$. Moreover, it can be easily seen that $\dco$ coincides with the metric of the Cauchy completion of the metric space $(\Omega,\dao)$, and that, for every $x,y\in\clO\setminus \IPO$,
$$
\dtc(x,y):=\dco(x,y) +\|x-y\| ^{\alpha},
$$
cf. \rf{DEF-TDA}.
\par In this setting the set $\IPO$ of inaccessible points, see Remark \reff{R-IP}, coincides with the set
$$
\IPO=\{x\in\DO:~\dco(x,\Omega)=\infty\}.
$$
In turn, the trace to $\DOA$, cf. \rf{DTR}, can be defined as
$$
\tra f(x)=\lim\{f(y):~\dco(y,x)\to 0,\, y\in\Omega\},
~~~x\in\DOA,
$$
provided $f\in\LOPO$, $p>n$, and $\alpha=\pn$.
\par Recall that every function $f\in\LOPO$ satisfies inequality \rf{H-LOC}. In Section 5 we show that this inequality implies the following property of the trace
$\tf=\tra f$: for every $x,y\in\DOA$
$$
|\tf(x)-\tf(y)|\le C(n,p)\|f\|_{\LOPO}\,\dco(x,y)^{1-\frac{1}{p}}.
$$
\par In particular, for every function $f\in\LOPO$, $p>n$, its trace $\tra f$ to the Sobolev $W^1_p$-boundary $\DOA$ is a {\it continuous (with respect to the metric $\dco$) function} on $\DOA$.
\bigskip
\par {\bf 1.4. Main results: a variational criterion of the trace and a characterization via sharp maximal functions.}
\par In \cite{S3} we studied the problem of characterizing the trace spaces $\LOP|_S$ and $\WP|_S$ to an arbitrary closed set $S\subset\RN$ whenever $p>n$. We gave various intrinsic characterizations of these trace spaces in terms of local oscillations and doubling measures supported on $S$. The approach introduced and used in \cite{S3} was based on an important property of the classical Whitney extension operator, namely that this operator provides an almost optimal extension of each function on $S$ which is the restriction to $S$ of a function in $\WP$.
\par Our approach here to Problem \reff{PR1} is an adaptation of the main ideas and methods of \cite{S3}. In particular, we show that the Whitney extension operator has a similar property to the one just mentioned: {\it it provides an almost optimal extension of every function defined on $\DO$ to a function from $\WPO$ whenever $p>n$.}
\par This enables us to characterize the boundary values of Sobolev $L^1_p$-functions and then of $W^1_p$-functions in ways similar to those presented in \cite{S3}. In particular, our first main result, Theorem \reff{VC-D}, is an analog of a trace criterion for the space $\LOP$ given in \cite{S3}.
\par  Before we recall that result, we will need to specify some more notation: Throughout this paper, the terminology ``cube'' will mean a closed cube in $\RN$ whose sides are parallel to the coordinate axes. We let $Q(x,r)$ denote the cube in $\RN$ centered at $x$ with side length $2r$. Given $\lambda >0$ and a cube $Q$ we let $\lambda Q$ denote the dilation of $Q$ with respect to its center by a factor of $\lambda $. (Thus $\lambda Q(x,r) =Q(x,\lambda r)$.) The Lebesgue measure of a measurable set $A\subset \RN$ will be denoted by $\left|A\right|$.
\par We proved in \cite{S3} that $f\in \LOP|_S$ for $n<p<\infty$, if and only if there exists a constant $\lambda >0$ such that for every finite family
$\{Q_{i}:i=1,...,m\}$ of pairwise disjoint cubes in $\RN$ and every choice of points $x_{i},y_{i}\in(\eta Q_{i}) \cap S$ the inequality
\bel{VDC}
\sum_{i=1}^{m}\frac{|f(x_{i})-f(y_{i})|^{p}}
{(\diam Q_{i})^{p-n}}\le \lambda
\ee
holds. Here $\eta$ is an absolute constant satisfying $\eta\le 11$. The special case of this result where $S=\RN$ and $\eta=1$  was treated earlier by Yu. Brudnyi \cite{Br1}.
\par We call the criterion expressed by \rf{VDC} {\it the variational or discrete characterization} of the trace space $\LOP|_S.$
\par Theorem \reff{VC-D}, which we will formulate in a moment, provides a ``variational" intrinsic characterization of the trace  space $\tra(\LOPO)$. This result can be thought of as a slight modification of the variational criterion of \rf{VDC}, where the requirement $x_{i},y_{i}\in (\eta Q_{i})\cap S$ for the points $x_{i}$ and $y_{i}$ is augmented by a certain  additional geometrical requirement, which we shall call  ``$Q_{i}$-visibility".
\begin{definition}\lbl{Q-V} {\em Given a point $x\in \clO$ and a cube $Q\subset \Omega $, we say that {\it $x$ is $Q$-visible in $\Omega $} if for each $y\in Q$ the semi-open line segment $(x,y]$ lies in $\Omega $.}
\end{definition}
\par In other words, a point $x\in \clO$ is $Q$-visible in $\Omega$ if
$$
\Conv(Q\cup\{x\})\setminus \{x\}\subset\Omega.
$$
Here $\Conv$ denotes the convex hull of a set. See Fig. 4.
\medskip
\begin{figure}[h]
\center{\includegraphics[scale=0.6]{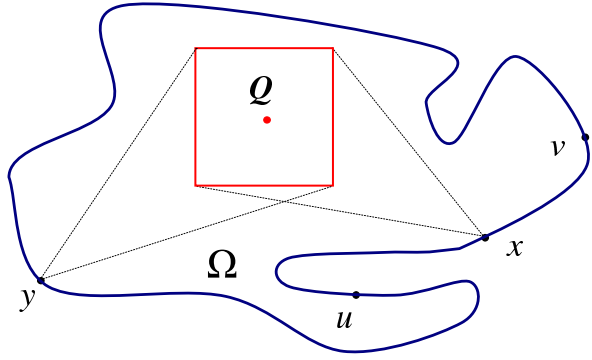}}
\caption{$Q$-visibility:\,The points $x$ and $y$ are $Q$-visible in $\Omega$ while $u$ and $v$ are not.}
\end{figure}
\begin{theorem}\lbl{VC-D} Let $\Omega$ be a domain in
$\RN$ satisfying condition $(A_\alpha)$. Let $p\in(n,\infty)$ and let $\alpha=\pn$. Let $\eta$ be a constant satisfying $\eta \ge 41$.
\par
A function $f:\DOA\to\R$ is the trace to $\DOA$ of a (continuous) function $F\in L_{p}^{1}(\Omega )$ if and only if $f$ is continuous (with respect to $\dco$) and there exists a constant $\lambda>0$ such that for every finite family $\{Q_{i}:i=1,...,m\}$ of pairwise disjoint cubes in $\Omega$ and every choice of $Q_{i}$-visible points
$$x_{i},y_{i}\in (\eta Q_{i})\cap \DOA, $$
the inequality
\bel{V-IN}
\sum_{i=1}^{m}\frac{|f(x_{i})-f(y_{i})|^{p}}
{(\diam Q_{i})^{p-n}}\le \lambda
\ee
holds. Moreover,
$$
\|f\|_{\tra(\LOPO)}\sim \inf\lambda ^{\frac{1}{p}}
$$
with constants of equivalence depending only on $n$, $p$ and $\eta$.
\end{theorem}
\bigskip
\begin{figure}[h]
\center{\includegraphics[scale=0.7]{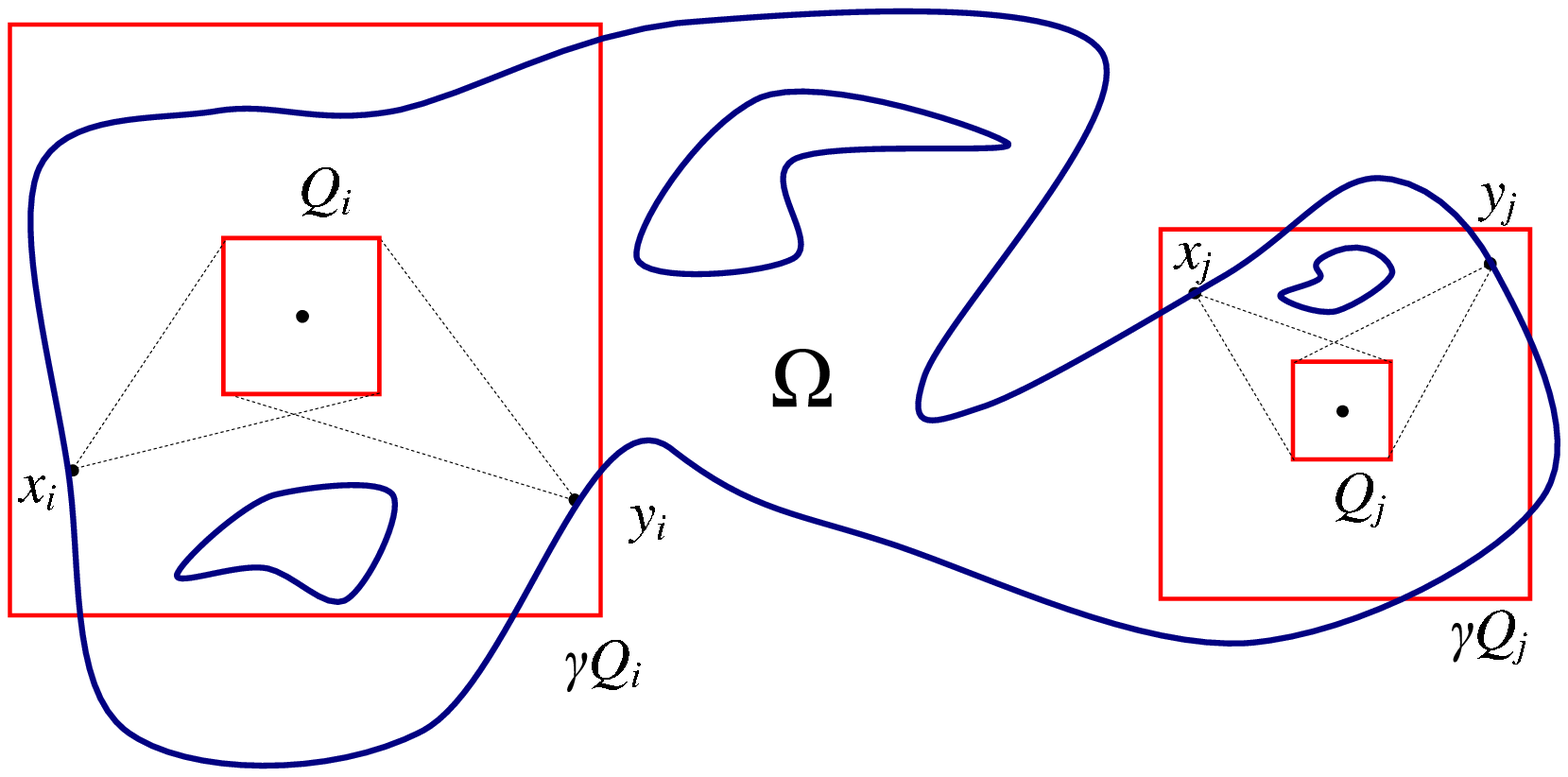}}
\caption{A variational criterion for the trace space $\tra(\LOPO)$.}
\end{figure}
\par We prove Theorem \reff{VC-D} in Section 5 as a corollary of Theorem \reff{VS-GN} which provides a variational trace criterion for an {\it arbitrary} domain $\Omega\subset\RN$. As we shall see, the more general criterion which appears in Theorem \reff{VS-GN} is a natural modification of the criterion in Theorem \reff{VC-D}, where the notion of $Q$-visibility is adapted to the case of a domain whose boundary admits ``agglutinated" points.
\par
We now turn to the second main result of the paper,
Theorem \reff{VC-WD}, which describes the traces of $\WPO$-functions to the boundary of $\Omega$. Here again, we first need some more terminology. Given $\ve>0$ we let
\bel{OC}
\DOEP:=\{x\in\Omega:~\D{x}<\ve\}
\ee
denote the $\ve$-neighborhood of $\DO$ in $\Omega$.
\par Let $\theta>1$ and let $\Qc=\{Q\}$ be a covering of $\Omega$ by non-overlapping cubes satisfying the following condition:
\bel{TH-W} \frac{1}{\theta}\diam(Q)\le\dist(Q,\DO)\le{\theta}\diam(Q). \ee
\par Let $Q$ be a cube in $\Omega$ and let $a_{Q}$ be a point in $\DO$ which is nearest to $Q$ (in the Euclidean metric). Let $T_{\Qc}:\Omega\to\DO$ be a mapping defined by the formula
\bel{T-C} T_{\Qc}|_Q:=a_Q,~~~~Q\in\Qc. \ee
Since $\Qc$ is a family of non-overlapping cubes, the
mapping $T_{\Qc}$ is well defined a.e. on $\Omega$.
\begin{theorem}\lbl{VC-WD} Let $\Omega$ be a domain in
$\RN$ satisfying condition $(A_\alpha)$ and let $p\in(n,\infty)$. Fix constants $\ve>0,\theta>1$, $\eta \ge 22\theta^2$, and an arbitrary covering $\Qc$ of $\Omega$ consisting of
non-overlapping cubes $Q\subset\Omega$ each satisfying inequality \rf{TH-W}.
\par A function $f:\DOA\to\R$ is an element of\,\,
$\tra(\WPO)$ if and only if $f$ is continuous (with respect to $\dco$), $f\circ T_{\Qc}\in L_{p}(\DOEP)$ and there exists a constant $\lambda>0$ such that the inequality \rf{V-IN} holds for every finite family $\{Q_{i}:i=1,...,m\}$ of pairwise disjoint cubes contained in
$\DOEP$, and every choice of $Q_i$-visible points
$x_{i},y_{i}\in (\eta\,Q_{i})\cap\DOA$.
\par Moreover,
$$ \|f\|_{\tra(\WPO)}\sim \|f\circ
T_{\Qc}\|_{L_{p}(\DOEP)} +\inf \lambda ^{\frac{1}{p}}
$$
with constants of equivalence depending only on $n$, $p$, $\ve$, $\theta$ and $\eta$.
\end{theorem}
\par A general version of this result, Theorem \reff{VCR-AO}, which characterizes the trace space $\tra(\WPO)$ for an {\it arbitrary} domain $\Omega\subset\RN$, is proven in Section 5.
\medskip
\begin{remark}\lbl{FQR} {\em Let $\Wc_\Omega$ be a
fixed Whitney decomposition of $\Omega$, i.e., a
covering of $\Omega$ by non-over\-lap\-ping
cubes such that $\diam Q\sim\dist(Q,\DO),$ $Q\in\Wc_\Omega$. Then the main statements of Theorems \reff{VC-D} and \reff{VC-WD}
remain true if we only consider cubes $\{Q_i\}$ which belong to $\Wc_\Omega$. In fact, the only modification is that
in this case the corresponding constants in the formulations of these theorems will also depend on    parameters of the Whitney decomposition $\Wc_\Omega$.}
\end{remark}
\par Our next result provides a different kind of
description of trace spaces. It is expressed in terms of a
certain kind of maximal function with respect to the metric $\dao$. This maximal function is a variant of $f_{p}^{\sharp}$, the familiar {\it fractional sharp maximal function in} $L_{p}$:
$$
f_{p}^{\sharp }(x):=\sup_{r>0}\frac{1}{r}
\left( \frac{1}{|Q(x,r)|}%
\intl\limits_{Q(x,r)}|f(y)-f_{Q(x,r)}|^{p}\,dy\right)
^{\frac{1}{p}}.
$$
Here $f_{Q}:=|Q|^{-1}\intl_{Q}f\,dx$ denotes the average of
$f$ on the cube $Q$. For $p=\infty $ the corresponding definition is
$$
f_{\infty }^{\sharp }(x):=\esssup
\{|f(y)-f(z)|/r:~r>0,\,y,z\in Q(x,r)\}.
$$
In \cite{C1} Calder\'{o}n proved that, for $1<p\le \infty$, a function $f$ is in $\WP$ if and only if $f$ and
$f_{1}^{\sharp }$ are both in $L_{p}(\RN)$. See also \cite{CS}. We observe that inequality \rf{S-Lip-OM} allows us to replace $f_{1}^{\sharp }$ in this statement with the (bigger) fractional sharp maximal function $f_{\infty }^{\sharp }$ so that for $p>n$ and for every $f\in\WP$, we have
$$
\|f\|_{\WP}\sim \|f\|_{L_{p}(\RN)}+ \|f_{\infty}^{\sharp
}\|_{L_{p}(\RN)}.
$$
\par In \cite{S3} we introduced a variant
of the fractional sharp maximal function $f_{\infty }^{\sharp }$  defined with respect to an {\it arbitrary} closed set $S\subset\RN$. For a function $f:S\to\R$ and $x\in\RN$ it is given by
$$
f_{\infty,S}^{\sharp }(x):=
\esssup\{|f(y)-f(z)|/r:~r>0,\,y,z\in Q(x,r)\cap S\}.
$$
We proved that $\|f\|_{\LOP|_S}\sim \|f_{\infty,S}^{\sharp }\|_{L_p(\RN)}$ for every $f$ defined on $S$,
provided that $p>n$.
\par Theorem \reff{MF-OM} formulated below presents an analog of this result for the space $\tra(\LOPO)$ whenever $\alpha=\pn$ and $p>n$.
\par Let us introduce yet another variant of the fractional sharp maximal function, this time defined for functions on a domain $\Omega\subset\RN$.  Fix $q, n<q<p,$ and put $\beta:=\qn$. For simplicity, we will assume that $\Omega$ satisfies the conditions  $(A_\alpha)$ and $(A_\beta)$, see Definition \reff{A-C}. Thus we can identify the $\alpha$-boundary and the $\beta$-boundary of $\Omega$ with (possibly different) subsets of $\DO$. In Section 2 we show that $\dt\le C\db$ for some constant $C=C(\beta,n)$, see Corollary \reff{AB-CH}, so that $\DOB\subset\DOA\subset\DO$.
\par We recall that the metric $\dcb$ on $\clO$ is defined by formula \rf{D-AO}. We introduce a quasi-metric $\dlbeta$ on $\clO$ by setting
$$
\dlbeta(x,y):=d^{\frac{1}{\beta}}_{\beta,\clO}(x,y),
~~~x,y\in\clO.
$$
\par Given $x\in\clO$ and $r>0$, we let $B(x,r:\dlbeta)$ denote the closed ball in the quasi-metric space $(\clO,\dlbeta)$ with center $x$ and radius $r$:
$$
B(x,r:\dlbeta):=\{y\in\clO:~\dlbeta(x,y)\le r\}.
$$
\par Our new fractional sharp maximal function $f_{\infty,\beta,\Omega}^{\sharp}$ is defined for each $f:\DOA\to\R$ by
$$
f_{\infty,\beta,\Omega}^{\sharp }(x):=
\esssup\left\{\frac{|f(y)-f(z)|}{r}:r>0,\,y,z\in  B(x,r:\dlbeta)\cap\DOB\right\},~~~x\in\Omega.
$$
\begin{theorem}\lbl{MF-OM} Let $n<q<p$ and let
$\beta=\qn$, $\alpha=\pn$. Let $\Omega$ be a domain in $\RN$ satisfying the conditions $(A_\alpha)$ and $(A_\beta)$.
\par (i). A function $f\in\tra(\LOPO)$ if and only if
$f$ is continuous on $\DOA$ (with respect to $\dco$) and
also $f_{\infty,\beta,\Omega}^{\sharp}\in L_p(\Omega)$.
Moreover,
\bel{E-L}
\|f\|_{\tra(\LOPO)}\sim
\|f_{\infty,\beta,\Omega}^{\sharp }\|_{L_p(\Omega)}.
\ee
\par (ii). Fix $\ve>0,\theta>1$ and a covering $\Qc$ of $\Omega$ consisting of non-overlapping cubes $Q\subset\Omega$ satisfying inequality \rf{TH-W}.
Let $T_{\Qc}$ be the mapping defined by \rf{T-C}.
\par A function $f:\DOA\to\R$ is an element of \,\,
$\tra(\WPO)$ if and only if  $f$ is continuous on $\DOA$ (with respect to $\dco$), and $f\circ T_{\Qc}$ and
$f_{\infty,\beta,\Omega}^{\sharp }$ are both in
$L_{p}(\DOEP)$. Furthermore,
\bel{E-W}
\|f\|_{\tra(\WPO)}\sim \|f\circ T_{\Qc}\|_{L_p(\DOEP)}+ \|f_{\infty,\beta,\Omega}^{\sharp }\|_{L_p(\DOEP)}.
\ee
\par The constants of equivalence in \rf{E-L} depend only
on $n,p$ and $q$, and in \rf{E-W} they only depend on
$n,p,q,\ve$ and $\theta$.
\end{theorem}
\smallskip
\par This theorem is a particular case of a corresponding result for an {\it arbitrary} domain, Theorem \reff{MF-C}, which we prove in Section 6.
\par As already mentioned above, the classical Whitney
extension operator provides an almost optimal extension of
a function $f$ defined on $\DOA$ to a function $F$ in
$\WPO$, whenever $p>n$ and $\alpha=\pn$. Since this extension operator is {\it linear}, we obtain the following
\begin{theorem} \lbl{LINEXT} Let $\alpha=\pn$. For every domain $\Omega\subset\RN$ and every $p>n$ there exists a
continuous linear operator $E:\tra(\WPO)\to \WPO$ such that $ \tra(Ef)=f$ for each $f\in\tra(\WPO)$. Its operator norm is bounded by a constant depending only on $n$ and $p$.
\end{theorem}
\par A similar result holds for the space $\LOPO$.
\section*{{\normalsize 2. Subhyperbolic metrics on a domain and chains of cubes}}
\setcounter{section}{2}
\setcounter{theorem}{0}
\setcounter{equation}{0}
\par Throughout the paper $C,C_1,C_2,...$ will be
generic positive constants which depend only on parameters
determining sets (say, $n,\alpha,\theta,$ etc.) or
function spaces ($p,q,$ etc). These constants can change
even in a single string of estimates. The dependence of a
constant on certain parameters is expressed, for example,
by the notation $C=C(n,p)$. We write $X\sim Y$ if there is
a constant $C\ge 1$ such that $X/C\le Y\le CX$.
\par Recall that by $\left|A\right|$ we denote the Lebesgue measure of a measurable set $A\subset\RN$.
\par It will be convenient for us to measure distances in
$\RN$ in the uniform norm
$$
\|x\|:=\max\{|x_i|:~i=1,...,n\}, \ \ \
x=(x_1,...,x_i)\in\RN.
$$
Thus every cube
$$ Q=Q(x,r):=\{y\in\RN:\|y-x\|\le r\} $$
is a ``ball" in $\|\cdot\|$-norm  of ``radius" $r$ centered at $x$. Given subsets $A,B\subset \RN$, we put
$$ \diam A:=\sup\{\|a-a'\|:~a,a'\in A\} $$
and
$$
\dist(A,B):=\inf\{\|a-b\|:~a\in A, b\in B\}.
$$
For $x\in \RN$ we also set $\dist(x,A):=\dist(\{x\},A)$. By $\overline{A}$ we denote the closure of $A$, and by $A^{\circ}$ its interior.
\par Let $\Qc=\{Q\}$ be a family of cubes in $\RN$. By
$M(\Qc)$ we denote its {\it covering multiplicity}, i.e.,
the minimal positive integer $M$ such that every point
$x\in\RN$ is covered by at most $M$ cubes from $\Qc$. Finally, given a cube $Q\subset\Omega$, by $a_{Q}$ we denote a point in $\DO$ which is nearest to $Q$ in the Euclidean metric.
\par In this section we present a series of results related to the proof of the necessity part of the main theorems. We begin with geometrical characterizations of the intrinsic metrics $\dao$ and $\dt$ introduced in subsection 1.2, see \rf{L-SH}, \rf{DEF-D} and \rf{DEF-TDA}.
\par In the next two lemmas we estimate the subhyperbolic length of a line segment in a domain.
\begin{lemma}\lbl{Q-OM} Let $\Omega$  be a domain in $\RN$ and let $x,y\in\Omega$. Assume that
\bel{E-LXY}
\|x-y\|\le\max\{\D{x},\D{y}\}.
\ee
Then for every $\alpha\in(0,1]$ we have
\bel{DXY}
\intl_{[x,y]}
\D{z}^{\alpha-1}\,ds(z)
\le\tfrac{1}{\alpha}\|x-y\|^\alpha.
\ee
\end{lemma}
\par {\it Proof.} Suppose that $\|x-y\|\le\D{x}$ so that
$Q(x,\|x-y\|)\subset \Omega.$
\par Let $z\in[x,y]$. Then
$\|x-y\|=\|x-z\|+\|z-y\|$
so that for every $u\in\Omega$ such that $\|u-z\|\le\|y-z\|$ we have
$$
\|x-u\|\le\|x-z\|+\|z-u\|\le \|x-z\|+\|z-y\|=\|x-y\|.
$$
We obtain $u\in Q(x,\|x-y\|)\subset \Omega$ so that $Q(z,\|y-z\|)\subset\Omega$. This proves that $\D{z}\ge \|y-z\|$. Hence
$$
\intl_{[x,y]}
\D{z}^{\alpha-1}\,ds(z)
\le\intl_{[x,y]}\|y-z\|^{\alpha-1}\,ds(z)
=\tfrac{1}{\alpha}\|x-y\|^\alpha.\BX
$$
\begin{lemma}\lbl{DS-QV} Let $\Omega$  be a domain in $\RN$ and let  $Q=Q(a,r)$ be a cube in $\Omega$. Let $x\in\clO$ be a $Q$-visible point, see Definition \reff{Q-V}, and let  $\beta\in(0,1]$. Then:
\par (i). For every $b_1,b_2\in(x,a]$  we have
\bel{FP}
\intl_{[b_1,b_2]}
\D{z}^{\beta-1}\,ds(z)
\le C(\beta)\left(\frac{\|a-x\|}{\diam Q}\right)^{1-\beta}\|b_1-b_2\|^\beta;
\ee
\par (ii). Every sequence $(x_i)\subset(x,a]$ which tends to $x$ (in $\|\cdot\|$-norm) is fundamental with respect to the metric $\db$. Moreover, every two sequences $(x_i),(y_i)\in(x,a]$ tending to $x$ in the Euclidean norm are equivalent with respect to $\db$, i.e.,
$$
\lim_{i\to\infty}\db(x_i,y_i)=0.
$$
\end{lemma}
\par {\it Proof.} Prove (i). Since $x$ is $Q$-visible,  $\Conv(x,Q)\setminus\{x\}\subset\Omega$ so that $tx+(1-t)Q\subset\Omega$
for every $t\in[0,1)$. Therefore for every  $z\in[a,x)$ we have
$$
Q(z,r_z)\subset\Omega~~~\text{where}~~~
r_z:=\frac{\|z-x\|}{\|a-x\|}\,r
$$
so that
$$
\D{z}\ge r_z=\frac{\|z-x\|}{\|a-x\|}\,r,~~~z\in[a,x).
$$
Suppose that $\|b_1-x\|\le \|b_2-x\|$. We have
\be
I&=&\intl_{[b_1,b_2]}
\D{z}^{\beta-1}\,ds(z)\le
\intl_{[b_1,b_2]}
\left(\frac{\|z-x\|r}{\|a-x\|}\right)^{\beta-1}\,ds(z)\nn\\
&=&\left(\frac{\|a-x\|}{r} \right)^{1-\beta}\intl_{[b_1,b_2]}
\|z-x\|^{\beta-1}\,ds(z)
=\left(\frac{\|a-x\|}{r}\right)^{1-\beta} \intl_{\|b_1-x\|}^{\|b_2-x\|}s^{\beta-1}\,ds.\nn
\ee
Hence
$$
I\le
\frac{1}{\beta}\left(\frac{\|a-x\|}{r}\right)^{1-\beta}
(\|b_2-x\|^{\beta}-\|b_1-x\|^{\beta})\le
\frac{2^{1-\beta}}{\beta}
\left(\frac{\|a-x\|}{\diam Q} \right)^{1-\beta}\|b_1-b_2\|^{\beta}
$$
proving \rf{FP}.
\par Prove (ii). By \rf{FP}, for every $u,v\subset(x,a]$ we have
$$
\dbo(u,v)\le\intl_{[u,v]}
\D{z}^{\beta-1}\,ds(z)
\le C(\beta)\left(\frac{\|a-x\|}{\diam Q}\right)^{1-\beta}\|u-v\|^\beta
$$
so that
\bel{SUV}
\db(u,v)=\dbo(u,v)+\|u-v\|^\beta\le
A\,\|u-v\|^\beta
\ee
where
$$
A=1+C(\beta)\left(\frac{\|a-x\|}{\diam Q}\right)^{1-\beta}.
$$
\par Now let $(x_i)\subset(x,a]$ and let $x_i\stackrel{\EN}{\longrightarrow} x$. Then, by \rf{SUV}, 
$$
\db(x_i,x_j)\le A\,\|x_i-x_j\|^\beta\to 0~~~\text{as}~~i,j\to\infty,
$$
proving that $(x_i)$ is fundamental with respect to $\db$.
\par Let $(x_i),(y_i)\in(x,a]$ and $(x_i),(y_i)$ tend to $x$ (in $\EN$-norm). Applying \rf{SUV} with $u=x_i$ and $v=y_i$ we obtain
$$
\db(x_i,y_i)\le A\,\|x_i-y_i\|^\beta\to 0~~~\text{as}~~i\to\infty.
$$
\par The lemma is proved.\bx
\begin{lemma}\lbl{L-QV} Let  $Q=Q(x,r),Q'=Q(x',r')$ be cubes in a domain $\Omega$ such that $Q\cap Q'\ne\emp$. Assume that either $\diam Q'\le\dist(Q,\DO)$ or $Q'\subset Q.$
\par Let $a_Q$ be a point on $\DO$ nearest to $Q$ (in the Euclidean norm). Then $a_Q$ is $Q'$-visible in $\Omega$ and for each $y\in(a_Q,x]$ there exists a point  $y'\in(a_Q,x']$ such that
\bel{R1}
\|y-y'\|\le C\|y-a_Q\|,
\ee
and
$$
\dao(y,y')\le C\|y-a_Q\|^\alpha,~~~\alpha\in(0,1],
$$
see Fig. 6. Here $C$ is a constant depending only on $x,x',r,r',a_Q$, and $\alpha$.
\end{lemma}
\begin{figure}[h]
\center{\includegraphics[scale=0.74]{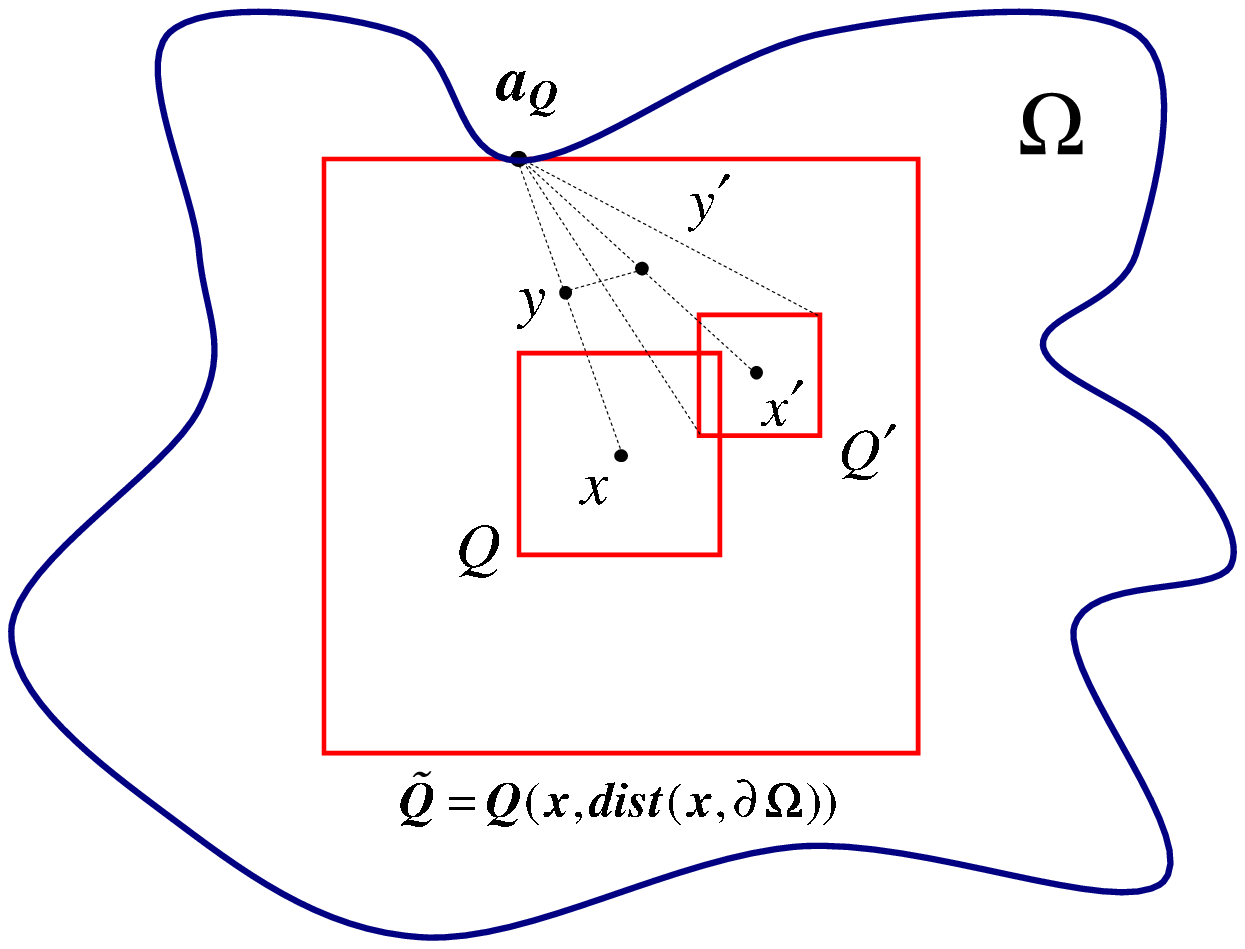}}
\caption{Lemma \reff{L-QV}: the case $Q\cap Q'\ne\emp$ and $\diam Q'\le\dist(Q,\DO)$.}
\end{figure}
\par {\it Proof.} Put $\lambda=\|a_Q-y\|/\|a_Q-x\|$ and $y'=a_Q+\lambda(x'-a_Q)$. Since $y\in(a_Q,x]$, we have $0<\lambda\le 1$ so that $y'\in(a_Q,x']$. In addition,
$y=a_Q+\lambda(x-a_Q).$ Hence
$$
\|y-y'\|=\lambda\|x-x'\|=\|a_Q-y\|\,\|x-x'\|/\|a_Q-x\|,
$$
proving \rf{R1}.
\par Let $z\in[y,y']$ and let $\tQ:=Q(x,\dist(x,\DO))$. Then $a_Q\in\partial\tQ$ and $\tQ^{\circ}\subset \Omega$.
\par If $\diam Q'\le\dist(Q,\DO)$, then  $Q'\subset\tQ^{\circ}.$ Clearly, the same is true whenever $Q'\subset Q$. The open cube $\tQ^{\circ}$ is a convex set and the point $a_Q$ lies on its boundary so that $a_Q$ is $Q'$-visible in $\Omega$, see Definition \reff{Q-V}.
\par Moreover, by dilation with respect to $a_Q$, we obtain
$$
Q_1:=Q(y,\lambda r),Q_2:=Q(y',\lambda r')\subset \tQ^{\circ}\subset\Omega.
$$
Since $\tQ^{\circ}$ is a convex set, the convex hull $\Conv(Q_1,Q_2)\subset \tQ^{\circ}\subset\Omega.$ In particular, if $\theta\in[0,1]$ and $z=\theta y+(1-\theta)y'$, then
$$
\theta Q_1+(1-\theta)Q_2\subset\Omega.
$$
Let $r''=\min\{r,r'\}.$ Then
$$
Q(z,\lambda r'')\subset\theta Q_1+(1-\theta)Q_2\subset\Omega,
$$
so that $\D{z}\ge\lambda r''.$ Hence
$$
\dao(y,y')\le\intl_{[y,y']}
\D{z}^{\alpha-1}\,ds(z)
\le (\lambda r'')^{\alpha-1}\intl_{[y,y']}
1\,ds(z)
$$
so that
\be
\dao(y,y')
&\le& (\lambda r'')^{\alpha-1}\|y-y'\|=(\lambda r'')^{\alpha-1}(\lambda\|x-x'\|)\nn\\&=&
\lambda^{\alpha} (r'')^{\alpha-1}\|x-x'\|=(\|a_Q-y\|/\|a_Q-x\|)^{\alpha} (r'')^{\alpha-1}\|x-x'\|.\nn
\ee
\par The lemma is proved.\bx
\begin{definition}\lbl{CH} {\em Let $\Omega$ be a domain in $\RN$ and let $x,y\in\Omega$. A finite family of cubes $\{Q_i\subset\Omega:~i=1,...,m\}$ is said to be a {\it chain of cubes} joining $x$ to $y$ in $\Omega$ if $x\in Q_1,y\in Q_m$ and $Q_i\cap Q_{i+1}\ne\emp$ for every $i=1,...,m-1.$}
\end{definition}
\begin{lemma}\lbl{CH-A} Let $x,y\in\Omega$ and let $\alpha\in(0,1]$. Then for every chain of cubes $\{Q_i:~i=1,...,m\}$ joining $x$ to $y$ in $\Omega$ the following inequality
$$
\dt(x,y)\le(1+2/\alpha)\sum_{i=1}^m(\diam Q_i)^\alpha
$$
holds.
\end{lemma}
\par {\it Proof.} Let $Q_i=Q(x_i,r_i),i=1,...,m,$ and let $z_i\in Q_i\cap Q_{i+1}$, $i=1,...,m-1.$ Put $z_0:=x$ and
$z_m:=y$. Then
$$
\|x-y\|^\alpha=\|z_0-z_m\|^\alpha\le
\left(\sum_{i=0}^{m-1}\|z_i-z_{i+1}\|\right)^\alpha\le
\sum_{i=0}^{m-1}\|z_i-z_{i+1}\|^\alpha\le
\sum_{i=1}^m(\diam Q_i)^\alpha.
$$
\par Let $\gamma$ be a broken line with nodes $\{z_0,x_1,z_1,x_2,z_2,...,x_{m-1},z_{m-1},x_m,z_m\}$. Then
\be
\dao(x,y)&\le&\intl_{\gamma}
\D{z}^{\alpha-1}\,ds(z)\nn\\
&=&
\sum_{i=0}^{m-1}\left\{\,\,\intl_{[z_i,x_{i+1}]}
\D{z}^{\alpha-1}\,ds(z)+
\intl_{[x_{i+1},z_{i+1}]}
\D{z}^{\alpha-1}\,ds(z)\,\right\}.\nn
\ee
Since $z_i\in Q_{i+1}=Q(x_{i+1},r_{i+1})\subset\Omega$,
$$
\|z_i-x_{i+1}\|\le r_{i+1}\le\D{x_{i+1}}
$$
so that, by Lemma \reff{Q-OM},
$$
\intl_{[z_i,x_{i+1}]}
\D{z}^{\alpha-1}\,ds(z)\le \tfrac{1}{\alpha}\|z_i-x_{i+1}\|^\alpha\le \tfrac{1}{\alpha}\left(\diam Q_{i+1}\right)^\alpha.
$$
In a similar way we prove that
$$
\intl_{[x_{i+1},z_{i+1}]}
\D{z}^{\alpha-1}\,ds(z)\le \tfrac{1}{\alpha}\left(\diam Q_{i+1}\right)^\alpha.
$$
Hence
$$
\dao(x,y)\le
(2/\alpha)\sum_{i=0}^{m-1}\left(\diam Q_{i+1}\right)^\alpha.
$$
Finally,
$$
\dt(x,y)=\|x-y\|^\alpha+\dao(x,y)\le
(1+2/\alpha)\sum_{i=1}^{m}\left(\diam Q_{i}\right)^\alpha
$$
proving the lemma.\bx
\begin{lemma}\lbl{CH-DO}(\cite{S4}) Let $x,y\in\Omega$ and let $\gamma\subset\Omega$ be a continuous curve joining $x$ to $y$. There exists a chain of cubes $\Qc=\{Q_i=Q(z_i,r_i):~i=1,...,m\}$ joining $x$ to $y$ in $\Omega$ such that:
\par (i). $z_i\in\gamma$ and $r_i=\tfrac18\D{z_i}$ for every $i=1,...,m$;
\par (ii). $2Q_i\subset\Omega,$ $i=1,...,m$;
\par (iii). The covering multiplicity $M(\Qc)$ of the family of cubes $\Qc$ is bounded by a constant $C=C(n)$.
\end{lemma}
\begin{lemma}\lbl{CH-XY} For every $x,y\in\Omega$ and every $\alpha\in(0,1]$ there exists a chain of cubes
$$
Ch_{\alpha,\Omega}(x,y)=\{Q_i\subset\Omega:~i=1,...,m,\}
$$
joining $x$ to $y$ in $\Omega$ with covering multiplicity $M(Ch_{\alpha,\Omega}(x,y))\le C(n)$ such that the following inequality
$$
\sum_{i=1}^{m}\left(\diam Q_{i}\right)^\alpha\le C(\alpha,n)\dt(x,y)
$$
holds.
\end{lemma}
\par {\it Proof.} By \rf{L-SH} and \rf{DEF-D}, there exists a rectifiable curve $\gamma\subset\Omega$ joining $x$ to $y$ such that
\bel{C-GAM}
\intl_{\gamma}
\D{z}^{\alpha-1}\,ds(z)\le 2\dao(x,y).
\ee
\par Let $\Qc=\{Q_i=Q(z_i,r_i):~i=1,...,m,\}$ be a chain of cubes joining $x$ to $y$ in $\Omega$ and satisfying conditions (i)-(iii) of Lemma \reff{CH-DO}.
\par Let us consider two cases. First suppose that $x,y\in Q_{k}=Q(z_{k},r_{k})\subset\Omega$ for some $k\in\{1,...,m\}$. Put $a:=(x+y)/2$ and $\tQ:=Q(a,\|x-y\|/2)$. Clearly, $x,y\in\tQ$. Since $\|x-y\|\le 2r_{k}$ and $a\in Q_{k}$, we have $\tQ\subset Q(a,r_{k})\subset 2Q_{k}$. But, by property (ii) of Lemma \reff{CH-DO}, $2Q_{k}\subset\Omega$ so that $\tQ\subset\Omega$.
\par Moreover,
$$
(\diam \tQ)^\alpha=\|x-y\|^\alpha\le \dt(x,y)=\|x-y\|^\alpha+\dao(x,y),
$$
so that in this case one can put $Ch_{\alpha,\Omega}(x,y):=\{\tQ\}.$
\par Now suppose that that for every $Q=Q(z,r)\in \Qc$ either $x\notin Q$ or $y\notin Q$. This implies the existence of a point $u\in\partial Q\cap\gamma$ such that $\gamma_{zu}\subset Q$. Here $\gamma_{zu}$ denotes the arc of $\gamma$ joining $z$ to $u$.
\par On the other hand, by property (i) of Lemma of \reff{CH-DO}, $8r=\D{z}$ so that for every $v\in Q$ we have
$$
\D{v}\le \D{z}+\|v-z\|\le 8r+r=9r.
$$
Hence,
\be
\intl_{\gamma\cap Q}
\D{v}^{\alpha-1}\,ds(v)&\ge& \intl_{\gamma_{zu}}
\D{v}^{\alpha-1}\,ds(v)\nn\\&\ge& (9r)^{\alpha-1}\intl_{\gamma_{zu}}
1\,ds(v)=(9r)^{\alpha-1}\lng(\gamma_{zu}).\nn
\ee
Since $u\in\partial Q$, we obtain
$$
\lng(\gamma_{zu})\ge\|z-u\|=r
$$
so that
$$
\intl_{\gamma\cap Q}
\D{v}^{\alpha-1}\,ds(v)
\ge (9r)^{\alpha-1}r=9^{\alpha-1}r^\alpha=9^{\alpha-1}(\diam Q/2)^\alpha.
$$
\par Thus for every cube $Q_i\in\Qc$, $i=1,...,m,$ we have
$$
(\diam Q_i)^\alpha\le C(\alpha)\intl_{\gamma\cap Q_i}
\D{v}^{\alpha-1}\,ds(v)
$$
so that
$$
\sum_{i=1}^{m}\left(\diam Q_{i}\right)^\alpha\le C(\alpha)\sum_{i=1}^{m}\,\intl_{\gamma\cap Q_i}
\D{v}^{\alpha-1}\,ds(v)\nn\\
\le C(\alpha)M(\Qc)\intl_{\gamma}
\D{v}^{\alpha-1}\,ds(v).
$$
Recall that $M(\Qc)$ denotes  the covering multiplicity of the family of cubes $\Qc$.
\par By property (iii) of Lemma \reff{CH-DO}, $M(\Qc)\le C(n)$ so that
$$
\sum_{i=1}^{m}\left(\diam Q_{i}\right)^\alpha\le C(\alpha,n)\intl_{\gamma}
\D{v}^{\alpha-1}\,ds(v).
$$
Combining this inequality with \rf{C-GAM}, we obtain
$$
\sum_{i=1}^{m}\left(\diam Q_{i}\right)^\alpha\le C(\alpha,n)\dao(x,y)
$$
proving the lemma.\bx
\par Lemma \reff{CH-A} and Lemma \reff{CH-XY} imply the following
\begin{proposition}\lbl{GE-CHR} For every $x,y\in\Omega$ and every $\alpha\in(0,1]$
$$
\dt(x,y)\sim\inf_{Ch}\,\,\sum_{i=1}^{m}\left(\diam Q_{i}\right)^\alpha
$$
where the infimum is taken over all chains of cubes $Ch=\{Q_i:i=1,...,m\}$ joining $x$ to $y$ in $\Omega$.
Moreover, the same equivalence holds whenever  $Ch$ runs over all chains of cubes joining $x$ to $y$ in $\Omega$ with the covering multiplicity $M(Ch)\le C$ where $C=C(n)$ is a constant depending only on $n$.
\par In both cases the constants of equivalence depend only on $n$ and $\alpha$.
\end{proposition}
\par This proposition implies the following
\begin{corollary}\lbl{AB-CH} For every $0<\beta<\alpha\le 1$ and every $x,y\in\Omega$ the following inequality
$$
\dt(x,y)\le C(\beta,n)\,\db(x,y)
$$
holds.
\end{corollary}
\begin{lemma}\lbl{CMP} For every $x,y\in\Omega$ and every $\alpha\in(0,1]$
$$
\dt(x,y)\sim\left \{
\begin{array}{ll}
\dao(x,y),& \|x-y\|\ge\min\{\D{x},\D{y}\},\\
\\
\|x-y\|^\alpha,& \|x-y\|<\min\{\D{x},\D{y}\},
\end{array}
\right.
$$
with constants of equivalence depend only on $\alpha$.
\end{lemma}
\par {\it Proof.} Assume that $\D{x}\ge\D{y}$. Let us consider the case
\bel{FCI}
\|x-y\|\ge\min\{\D{x},\D{y}\}=\D{y}.
\ee
Prove that in this case
\bel{IXYD}
\|x-y\|^\alpha\le 2^{1-\alpha}\dao(x,y).
\ee
\par In fact, for every $z\in[x,y]$ we have
$$
\D{z}\le\D{y}+\|z-y\|\le\D{y}+\|x-y\|.
$$
Hence, for every rectifiable curve $\gamma\subset\Omega$ joining $x$ to $y$ we obtain
\be
\intl_{\gamma}\D{z}^{\alpha-1}\,ds(z)
&\ge& \intl_{\gamma}(\D{y}+\|x-y\|)^{\alpha-1}\,ds(z)\nn\\
&=&(\D{y}+\|x-y\|)^{\alpha-1}\lng(\gamma)\nn\\&\ge&
(2\|x-y\|)^{\alpha-1}\|x-y\|=2^{\alpha-1}\|x-y\|^\alpha.\nn
\ee
Taking  in this inequality the infimum over all such $\gamma$ we obtain inequality \rf{IXYD}.
\par Thus
$$
\dao(x,y)\le\dt(x,y)=\dao(x,y)+\|x-y\|^\alpha\le (1+2^{\alpha-1})\dao(x,y)
$$
provided inequality \rf{FCI} is satisfied.
Let us consider the case
$$
\|x-y\|\le\min\{\D{x},\D{y}\}.
$$
In this case inequality \rf{E-LXY} is satisfied so that, by Lemma \reff{Q-OM},
$$
\dao(x,y)\le\intl_{[x,y]}\D{z}^{\alpha-1}\,ds(z) \le\tfrac{1}{\alpha}\|x-y\|^\alpha.
$$
Hence,
$$
\|x-y\|^\alpha\le \dt(x,y)=\dao(x,y)+\|x-y\|^\alpha\le (1+\tfrac{1}{\alpha})\|x-y\|^\alpha
$$
proving the lemma.\bx
\section*{{\normalsize 3. Sobolev-Poincar\'{e} type
inequalities on a domain.}}
\setcounter{section}{3}
\setcounter{theorem}{0}
\setcounter{equation}{0}
\par The following proposition presents the classical Sobolev imbedding inequality for the case $p>n$, see, e.g. \cite{M}, p. 61, or \cite{MP}, p. 55. This inequality is also known in the literature as Sobolev-Poincar\'e inequality (for $p>n$).
\begin{proposition}\lbl{EDIFF} Let $F\in \LOPO$ be
a continuous function defined on a domain
$\Omega\subset\RN$ and let $n<q\le p<\infty$. Then for
every cube $Q\subset\Omega$ and every $x,y\in Q$ the
following inequality
\bel{SPE}
|F(x)-F(y)|\le C(n,q)\, (\diam Q) \left(\frac{1}{|Q|} \intl_Q
\|\nabla F(z)\|^q\,dz\right)^{\frac{1}{q}}
\ee
holds.
\end{proposition}
\par  Clearly, inequality \rf{SPE} (for $p=q$) implies the local H\"older inequality \rf{S-Lip-OM}.
\par In this section we present several {\it global} versions of inequality \rf{SPE} related to the case of {\it arbitrary} points $x,y\in\Omega$. These variants of the  Sobolev-Poincar\'e inequality on a domain are a slight generalization of the global H\"older-type inequality \rf{H-LOC} proved by Buckley and Stanoyevitch \cite{BSt3}.
\par Fix $q\in(n,p]$ and put
$$
\beta=\frac{q-n}{q-1}\,.
$$
\begin{lemma}\lbl{E-CHF} Let $F\in \LOPO$ be
a continuous function and let $x,y\in\Omega$. Let  $Ch=\{Q_1,...,Q_m\}$ be a chain of cubes joining $x$ to $y$ in $\Omega$ with the covering multiplicity $M=M(Ch)<\infty$.
 Then
$$
|F(x)-F(y)|\le C(n,q)\,M^{\frac1q}\left(\,\sum_{i=1}^m\,(\diam Q_i)^\beta\right)^{1-\frac1q}\left(\,\,\intl_U
\|\nabla F(z)\|^q\,dz\right)^{\frac{1}{q}}.
$$
where $U:=\cup_{i=1}^m Q_i$.
\end{lemma}
\par {\it Proof.} Let $z_i\in Q_i\cap Q_{i+1}$, $i=1,...,m-1.$ Put $z_0:=x, z_m=y$. Then for every $i=0,...,m-1,$ we have  $z_i,z_{i+1}\in Q_{i+1}\subset\Omega$ so that by the Sobolev-Poincar\'e inequality \rf{SPE}
$$
|F(z_i)-F(z_{i+1})|\le C (\diam Q_{i+1}) \left(\frac{1}{|Q_{i+1}|} \intl_{Q_{i+1}}
\|\nabla F(z)\|^q\,dz\right)^{\frac{1}{q}}.
$$
Hence
\bel{IFZ}
|F(z_i)-F(z_{i+1})|\le C(\diam Q_{i+1})^{1-\frac{n}{q}} \left(\,\,\intl_{Q_{i+1}}
\|\nabla F(z)\|^q\,dz\right)^{\frac{1}{q}}.
\ee
\par Now, by the H\"older inequality,
\be
|F(x)-F(y)|&\le&\sum_{i=0}^{m-1}\frac{|F(z_i)-F(z_{i+1})|} {(\diam Q_{i+1})^{1-\frac{n}{q}}}\, (\diam Q_{i+1})^{1-\frac{n}{q}}\nn\\
&\le&\left(\sum_{i=0}^{m-1}\left(\frac{|F(z_i)-F(z_{i+1})|} {(\diam Q_{i+1})^{1-\frac{n}{q}}}\right)^q\right) ^{\frac{1}{q}}
\left(\sum_{i=0}^{m-1}\left(\left(
\diam Q_{i+1}\right)^{1-\frac{n}{q}}\right)^{\frac{q}{q-1}}
\right)^{1-\frac{1}{q}}\nn
\\
&=&
\left(\sum_{i=0}^{m-1}\frac{|F(z_i)-F(z_{i+1})|^q} {(\diam Q_{i+1})^{q-n}}\right)^{\frac{1}{q}}
\left(\sum_{i=0}^{m-1}\left(\diam Q_{i+1}\right)^\frac{q-n}{q-1}\right)^{1-\frac{1}{q}}.\nn
\ee
Recall that $\beta=(q-n)/(q-1)$. Combining the latter inequality  with \rf{IFZ}, we obtain
\be
|F(x)-F(y)|&\le& C\,\left(\sum_{i=0}^{m-1}\left(\diam Q_{i+1}\right)^\beta\right)^{1-\frac{1}{q}}
\left(\sum_{i=0}^{m-1}\intl_{Q_{i+1}}
\|\nabla F(z)\|^q\,dz\right)^{\frac{1}{q}}\nn\\
&\le& C\,\left(\sum_{i=1}^m(\diam Q_i)^\beta\right)^{1-\frac1q}\left(M\intl_U
\|\nabla F(z)\|^q\,dz\right)^{\frac{1}{q}}.\nn
\ee
\par The lemma is proved.
\begin{proposition}\lbl{SP-OM} Let $F\in \LOPO$ be
a continuous function and let $q\in(n,p]$,  $\beta=(q-n)/(q-1)$. There exist constants $\lambda=\lambda(n,q)$ and $C=C(n,q)$ such that for every $x,y\in \Omega$
$$
|F(x)-F(y)|\le C\db(x,y) ^{1-\frac{1}{q}} \left(\intl_B
\|\nabla F(z)\|^q\,dz\right)^{\frac{1}{q}}
$$
where
$$
B=\{z\in\Omega:~\db(x,z)\le \lambda\db(x,y)\}.
$$

\end{proposition}
\par {\it Proof.} By Lemma \reff{CH-XY} there exists a chain of cubes $Ch_{\alpha,\Omega}(x,y)=\{Q_1,...,Q_m\}$ joining $x$ to $y$ in $\Omega$ with covering multiplicity $M(Ch_{\beta,\Omega}(x,y))\le C(n)$ such that
\bel{S-R}
\sum_{i=1}^{m}\left(\diam Q_{i}\right)^\beta\le C(\beta,n)\db(x,y).
\ee
Then by Lemma \reff{E-CHF}
$$
|F(x)-F(y)|\le C\,M^{\frac1q}\left(\,\sum_{i=1}^m(\diam Q_i)^\beta\right)^{1-\frac1q}\left(\,\,\intl_U
\|\nabla F(z)\|^q\,dz\right)^{\frac{1}{q}}
$$
where $U=\cup_{i=1}^m Q_i$. Hence
$$
|F(x)-F(y)|\le C\,\db(x,y)^{1-\frac1q}\left(\,\intl_U
\|\nabla F(z)\|^q\,dz\right)^{\frac{1}{q}}.
$$
\par Let  $z\in U$ so that $z\in Q_k$ for some $k\in\{1,...,m\}$. Then $\{Q_1,...,Q_k\}$ is a chain of cubes joining $x$ to $z$ in $\Omega$. By Lemma \reff{CH-A}
$$
\db(x,z)\le(1+2/\beta)\sum_{i=1}^k(\diam Q_i)^\beta.
$$
Hence, by \rf{S-R},
$$
\db(x,z)\le(1+2/\beta)\sum_{i=1}^m(\diam Q_i)^\beta\le\lambda\db(x,y)
$$
with $\lambda=\lambda(n,q)$. Finally we obtain
$$
U=\bigcup_{i=1}^m Q_i\subset B=\{z\in\Omega:~\rho_{\beta,\Omega}(x,z)\le \lambda
\rho_{\beta,\Omega}(x,y)\}
$$
proving the proposition.\bx
\begin{proposition}\lbl{SP-Q} Let $F\in \LOPO$ be
a continuous function and let $x,y\in\Omega$. Let $q\in(n,p]$ and let $\beta=(q-n)/(q-1)$. There exist constants $\lambda_1=\lambda_1(n,q)\ge 1$ and $C=C(n,q)$ such that for every $R\ge \lambda_{1}\db(x,y)^{\frac1\beta}$ the following inequality
$$
|F(x)-F(y)|\le C\,R\left(\frac{1}{|Q(x,R)|} \intl_{Q(x,R)\cap \Omega}
\|\nabla F(z)\|^q\,dz\right)^{\frac{1}{q}}
$$
holds.
\end{proposition}
\par {\it Proof.} Let $\lambda=\lambda(n,q)$ be the constant from Proposition \reff{SP-OM} and let
$$
B=\{z\in\Omega:~\db(x,z)\le \lambda\db(x,y)\}.
$$
Then for every $z\in B$ we have
$$
\|x-z\|^\beta\le\db(x,z)=
\|x-z\|^\beta+d_{\beta,\Omega}(x,z)\le\lambda\db(x,y)
$$
so that $B\subset Q(x,r)$ with $r=\lambda^{\frac1\beta} \db(x,y)^{\frac1\beta}.$
\par Hence, by Proposition \reff{SP-OM},
\be
|F(x)-F(y)|&\le& C\db(x,y) ^{1-\frac{1}{q}}\left(\,\intl_B
\|\nabla F(z)\|^q\,dz\right)^{\frac{1}{q}}\nn\\
&\le&
C(r^\beta)^{1-\frac{1}{q}}\left(\,\,\intl_{Q(x,r)\cap \Omega}
\|\nabla F(z)\|^q\,dz\right)^{\frac{1}{q}}=
C\,r^{1-\frac{n}{q}} \left(\,\,\intl_{Q(x,r)\cap \Omega}
\|\nabla F(z)\|^q\,dz\right)^{\frac{1}{q}}.\nn
\ee
\par Now, for every $R\ge r=\lambda_1 \db(x,y) ^{\frac1\beta}$ we have
\be
|F(x)-F(y)|&\le&
C\,R^{1-\frac{n}{q}} \left(\,\,\intl_{Q(x,R)\cap \Omega}
\|\nabla F(z)\|^q\,dz\right)^{\frac{1}{q}}\nn\\&\le&
C\,R\left(\frac{1}{|Q(x,R)|} \intl_{Q(x,R)\cap \Omega}
\|\nabla F(z)\|^q\,dz\right)^{\frac{1}{q}}\nn
\ee
proving the proposition.\bx
\begin{lemma}\lbl{E-QVX} Let  $Q=Q(a,r)$ be a cube in $\Omega$ and let $t\ge 1$. Assume that $x\in (tQ)\cap\Omega$ and $x$ is a $Q$-visible point of $\Omega$, see Definition \reff{Q-V}. Let $F\in \LOPO$ be
a continuous function. Then for every $q\in(n,p]$ we have
$$
\left(\frac{|F(x)-F(a)|}{\diam Q}\right)^q\le C\,\frac{1}{|Q|}\,\intl_{(\gamma Q)\cap \Omega}
\|\nabla F(z)\|^q\,dz
$$
where $\gamma=\gamma(n,q,t)$ and $C=C(n,q,t)$.
\end{lemma}
\par {\it Proof.} Let  $\beta=(q-n)/(q-1)$. By Lemma \reff{DS-QV} (with $b_1=x\in\Omega, b_2=a$),
$$
d_{\beta,\Omega}(x,a) \le
\intl_{[x,a]}
\D{z}^{\beta-1}\,ds(z)
\le C(\beta)\left(\frac{\|a-x\|}{\diam Q}\right)^{1-\beta}\|x-a\|^\beta.
$$
Since $x\in tQ$, we have $\|a-x\|\le tr=\tfrac12 t\diam Q$ so that
$$
d_{\beta,\Omega}(x,a)\le C(\beta)
t^{1-\beta}\|x-a\|^\beta.
$$
Hence
$$
\db(x,a)=\|x-a\|^\beta+d_{\beta,\Omega}(x,a) \le
C(\beta) t^{1-\beta}\|x-a\|^\beta
\le C(\beta) t^{1-\beta}(tr)^{\beta}=
C(\beta)tr^{\beta}.
$$
Let $\lambda_1=\lambda_1(n,q)$ be the constant from Proposition \reff{SP-Q}. Then
$$
\lambda_1\db(x,a)^{\frac1\beta}\le \lambda_1 (C(\beta)t)^{\frac1\beta}r=\gamma r
$$
with $\gamma:=\lambda_1 (C(\beta)t)^{\frac1\beta}$.
\par Put $R:=\gamma r=\tfrac12 \gamma\diam Q$. Since $R\ge\lambda_1\db(x,a)^{\frac1\beta}$, by Proposition \reff{SP-Q},
$$
|F(x)-F(a)|\le C\,R\left(\frac{1}{|Q(a,R)|} \intl_{Q(a,R)\cap \Omega}
\|\nabla F(z)\|^q\,dz\right)^{\frac{1}{q}}
$$
so that
$$
\left(\frac{|F(x)-F(a)|}{\diam Q}\right)^q\le
C\, 2^{-q}\gamma^{q-n}
\frac{1}{|Q|}\intl_{Q(a,R)\cap \Omega}
\|\nabla F(z)\|^q\,dz.
$$
\par The lemma is proved.\bx
\par As usual, given a function $G\in L_{1,loc}(\RN)$ by
$\Mc[G]$ we denote the Hardy-Littlewood maximal function
\bel{H-L-M}
\Mc[G](x):=\sup_{t>0}\frac{1}{|Q(x,t)|}
\intl_{Q(x,t)}|G(y)|dy,~~~~x\in\RN.
\ee
\par The last auxiliary result of the section is the  following
\begin{lemma}\lbl{MF} Let $G\in L_{1,loc}(\RN)$ and let $\gamma\ge 1,\theta>0$. Then for every cube $Q\subset\RN$
we have
$$
\left(\frac{1}{|Q|}\intl_{\gamma Q}
|G(x)|\,dx\right)^\theta\le C(n,\gamma)\frac{1}{|Q|}\intl_{Q}
\Mc[G]^\theta(x)\,dx.
$$
\end{lemma}
\par {\it Proof.} Let $z\in Q$ and let $K:=Q(z,\diam(\gamma Q)).$ Then $K\supset\gamma Q$ and $|K|\sim |Q|$ so that
$$
\frac{1}{|Q|}\intl_{\gamma Q}
|G(x)|\,dx\le \frac{1}{|Q|}\intl_{K}
|G(x)|\,dx\le C\frac{1}{|K|}\intl_{K}
|G(x)|\,dx\le C\,\Mc[G](z).
$$
where $C=C(n,\gamma)$. Hence
$$
\left(\frac{1}{|Q|}\intl_{\gamma Q}
|G(x)|\,dx\right)^\theta\le
C\,\Mc[G]^\theta(z),~~~~z\in Q.
$$
Integrating this inequality over cube $Q$ we obtain the required inequality
$$
|Q|\left(\frac{1}{|Q|}\intl_{\gamma Q}
|G(x)|\,dx\right)^\theta\le C\intl_{Q}
\Mc[G]^\theta(x)\,dx.\BX
$$
\section*{{\normalsize 4. Local oscillation properties of the Whitney extension operator.}}
\setcounter{section}{4}
\setcounter{theorem}{0}
\setcounter{equation}{0}
\par In this section we study local oscillation properties of the classical Whitney operator which extends  every function defined on the Sobolev boundary of a domain to a function defined on all of the domain. We present several auxiliary results which we will use later in the proofs of the sufficiency part of the main theorems, see Sections 5 and 6.
\par Since $\Omega$ is an open subset of $\RN$, it
admits a Whitney covering $\TW$. In the next lemma we recall the main properties of this covering, see, e.g. \cite{St}, or \cite{G}.
\begin{theorem}\lbl{Wcov}
$\TW=\{Q_k\}$ is a countable family of cubes  such that
\par (i). $\Omega=\cup\{Q:Q\in \TW\}$;
\par (ii). For every cube $Q\in \TW$
\bel{DKQ} \diam Q\le \dist(Q,\DO)\le 4\diam Q; \ee
\par (iii). The covering multiplicity $M(\TW)$ of the family $\TW$ is bounded by a constant $N=N(n)$. Thus every point of $\Omega$ is covered by at most $N$ cubes from
$\TW$.
\end{theorem}
\par We are also needed certain additional properties of
Whitney's cubes which we present in the next lemma. Let $Q$ be a cube in $\Omega$ and let $Q^*:=\tfrac98 Q.$
\begin{lemma}\lbl{Wadd}
(1). If $Q,K\in \TW$ and $Q^*\cap K^*\ne\emptyset$, then
\bel{DM} \frac{1}{4}\diam Q\le \diam K\le 4\diam Q. \ee
\par (2). For every cube $K\in \TW$ there are at most
$N=N(n)$ cubes from the family
$ {\mathbb W}^*(\Omega):=\{Q^*:Q\in \TW\} $
which intersect $K^*$.
\par (3). Let $Q,K\in \TW$. Then $Q^*\cap K^*\ne\emptyset$
if and only if  $Q\cap K\ne\emptyset$.
\end{lemma}
\par We turn to the construction of the Whitney extension operator. Fix $\alpha\in(0,1]$. As usual, given a cube $Q\subset\Omega$, by $a_Q$ we denote a point of $\DO$ nearest to $Q$ in the Euclidean norm.
\par We recall that the standard Whitney's extension algorithm assigns every function $f:\DO\to\R$ a piecewise constant function $F$ which on every cube $Q\in\TW$ takes the value $f(a_Q)$. Then we smooth  $F$ using a certain smooth partition of unity subordinated to the Whitney decomposition $\TW$.
\par Let now $f:\DOA\to\R$ be a function defined on the $\alpha$-boundary of $\Omega$, see Definition \reff{B-DA}.
Observe that the same extension procedure works well whenever the $\alpha$-boundary of $\Omega$ can be identified with a subset of $\DO$. A domain satisfying the $(A_\alpha)$-condition, see Definition \reff{A-C}, provides an example of such a domain; in fact, in this case the boundary $\DO$ does not contains ``agglutinated" parts and does not split under the Cauchy completion with respect to the metric $\dt$.
\par However, in general, the point $a_Q$ may split into a finite or infinite number of elements of the $\alpha$-boundary. In this case we have to assign the pair $(a_Q,Q)$ an appropriate equivalence class $\oqa\in\DOA$ and then to proceed the Whitney algorithm using the value $f(\omega_{Q,\alpha})$ rather than $f(a_Q)$.
\par We will do this as follows. Clearly, if a cube $Q=Q(x_Q,r_Q)\subset\Omega$ then
$$Q^\circ(x_Q,\D{x_Q})\subset\Omega$$
so that the point 
$$
a_Q~~\text{is}~~Q\text{-visible in}~~\Omega,
$$
see Definition \reff{Q-V}. In particular, $[x_Q,a_Q)\subset\Omega$. We define a sequence of points $x_i\in[x_Q,a_Q)$ by letting
\bel{XNE}
x_i:=a_Q+\tfrac1i(x_Q-a_Q), ~~~i=1,2,...
\ee
Clearly, $x_i\stackrel{\EN}{\longrightarrow} a_Q$ as $i\to\infty$ so that, by Lemma \reff{DS-QV}, $(x_i)$ is {\it fundamental with respect to the metric $\dt$}.
\par Thus $(x_i)\in \CSQ$, see \rf{D-CSQ}, so that the equivalence class of $(x_i)$, the set
\bel{DYO}
\oqa=[(x_i)]_\alpha,
\ee
is an element of $\OA$, see \rf{CMP-OM}.
\par Moreover, since  $x_i\stackrel{\EN}{\longrightarrow} a_Q$, the point $a_Q$ is the common limit (in $\EN$-norm) of all sequences $(y_i)\in\oqa$, i.e.,
\bel{OQA}
\ell(\oqa)=a_Q,
\ee
see Remark \reff{R-CPM}. On the other hand, since $a_Q\in\DO$, by Definition \reff{B-DA}, the class $\oqa\in\DOA$. In other words, $\oqa$ is an element of the $\alpha$-boundary of the domain $\Omega$:
$$
\oqa\in\DOA=\OA\setminus\Omega,
$$
see \rf{B-A}.
\par We will be needed the following property of the element $\oqa$.
\begin{lemma}\lbl{LO-QA} For every cube $Q\in\TW$ and every $y\in Q$ the following inequality
\bel{FL-L}
\dtc(y,\oqa)\le C(\alpha)\dist(y,\DO)^\alpha
\ee
holds.
\end{lemma}
\par {\it Proof.}  By \rf{CDR-L}, \rf{XNE} and \rf{OQA},
$$
\dtc(x_Q,\oqa)=
\lim_{i\to\infty}\dao(x_Q,x_i)+\|x_Q-a_Q\|.
$$
\par Since $a_Q$ is a nearest point to $Q$ on $\DO$,
$$
\|x_Q-a_Q\|\le r_Q+\dist(Q,a_Q)=r_Q+\dist(Q,\DO),
$$
so that, by \rf{DKQ},
\bel{XA-Q}
\|x_Q-a_Q\|\le r_Q+4\diam Q\le 5\diam Q.
\ee
\par Let us estimate $\dao(x_Q,x_i)$. Since $a_Q$ is $Q$-visible and $x_i\in[x_Q,a_Q)$, by Lemma \reff{DS-QV}, (i),
$$
\dao(x_Q,x_i)\le\intl_{[x_i,x_Q]}
\D{z}^{\alpha-1}\,ds(z)
\le C\left(\frac{\|a_Q-x_Q\|}{\diam Q}\right)^{1-\alpha}\|x_i-x_Q\|^\alpha,
$$
so that
$$
\lim_{i\to\infty}\dao(x_Q,x_i)\le  C\left(\frac{\|a_Q-x_Q\|}{\diam Q}\right)^{1-\alpha}\|a_Q-x_Q\|^\alpha.
$$
Combining this inequality with \rf{XA-Q}, we obtain
$$
\lim_{i\to\infty}\dao(x_Q,x_i)\le  C(\diam Q)^\alpha.
$$
From this and \rf{XA-Q} it follows
$$
\dtc(x_Q,\oqa)\le  C(\diam Q)^\alpha.
$$
\par Let us estimate $\dt(x_Q,y)$. We have $\|x_Q-y\|\le r_Q$ so that
\bel{R-OXY}
\dt(x_Q,y)=\dao(x_Q,y)+\|x_Q-y\|^{\alpha}\le \dao(x_Q,y) +(\diam Q)^\alpha.
\ee
But, by \rf{DKQ},
\bel{DY-D}
\dist(y,\DO)\ge \dist(Q,\DO)\ge \diam Q
\ee
proving that
$$
\|x_Q-y\|\le \diam Q\le \max\{\D{x_Q},\D{y}\}.
$$
Hence, by \rf{DXY},
$$
\dao(x_Q,y)\le\intl_{[x_Q,y]}
\D{z}^{\alpha-1}\,ds(z)
\le\tfrac{1}{\alpha}\|x_Q-y\|^\alpha\le C(\diam Q)^\alpha,
$$
so that, by \rf{R-OXY}, $\dt(x_Q,y)\le C(\diam Q)^\alpha.$ We obtain
$$
\dtc(y,\oqa)\le \dtc(x_Q,\oqa)+\dt(x_Q,y)\le C(\diam Q)^\alpha.
$$
This inequality and \rf{DY-D} imply the required inequality \rf{FL-L}.\bx
\smallskip
\par Let $\sigma>0, \bc\in\R,$ and let $f:\DOA\to\R$ be a function defined on the $\alpha$-boundary of $\Omega$.
We put
\bel{CQ} c_Q:=\left \{
\begin{array}{ll}
f(\oqa),& \diam Q\le \sigma,\\
\bc,& \diam Q>\sigma.
\end{array}
\right. \ee
We define an extension operator
$\tf=\Ex[f;\sigma,\bc,\alpha,\Omega]$ by letting $\tf(\omega):=f(\omega),$ $\omega\in\DOA$, and
\bel{ExtOp}
\tf(x):=\sum\limits_{Q\in \TW}
c_Q\varphi_Q(x), ~~~ x\in \Omega.
\ee
\par Here $\{\varphi_Q:Q\in \TW\}$ is a smooth partition of unity subordinated to the Whitney decomposition $\TW$, see, e.g. \cite{St}. Recall the main properties of this
partition.
\begin{lemma}\lbl{PU} The partition of unity  $\{\varphi_Q:Q\in \TW\}$ has the following properties:
\par (a). $\varphi_Q\in C^\infty(\RN)$ and
$0\le\varphi_Q\le 1$ for every $Q\in \TW$;\smallskip
\par (b). $\supp \varphi_Q\subset Q^*(:=\frac{9}{8}Q),$
$Q\in \TW$;\smallskip
\par (c). $\sum\{\varphi_Q(x):~Q\in \TW\}=1$ for every
$x\in\Omega$;\smallskip
\par (d). $ \|\nabla\varphi_Q(x)\| \le C(n)/\diam Q\,\,$
for every $Q\in \TW$ and every $x\in\Omega$.
\smallskip
\end{lemma}
\begin{lemma}\lbl{C-L} For every cube $K\in\TW$ with $\diam K\le\sigma/4$ and every $x\in K$ the following inequality
$$
|\tf(x)-f(\omega_{K,\alpha})|\le C(n)\max\{ |f(\omega_{Q,\alpha})-f(\omega_{K,\alpha})|:~Q\in\TW, Q\cap K\ne\emp\}
$$
holds.
\end{lemma}
\par {\it Proof.} By \rf{ExtOp} and Lemma \reff{PU}, (c), we have
\be
|\tf(x)-f(\omega_{K,\alpha})|&=&|\sum_{Q\in \TW}
c_Q\varphi_Q(x)-f(\omega_{K,\alpha})|=|\sum_{Q\in \TW}
(c_Q-f(\omega_{K,\alpha}))\varphi_Q(x)|\nn\\&\le&\sum_{Q\in \TW}|c_Q-f(\omega_{K,\alpha})|\,\varphi_Q(x).\nn
\ee
Hence, by Lemma \reff{PU}, (b), we obtain
\be
|\tf(x)-f(\omega_{K,\alpha})|&\le&
\sum\{|c_Q-f(\omega_{K,\alpha})|\varphi_Q(x):~Q\in \TW, Q^*\ni x\}\nn\\
&\le&\sum\{|c_Q-f(\omega_{K,\alpha})|\varphi_Q(x):~Q\in \TW, Q^*\cap K^*\ne\emp\}.\nn
\ee
But $0\le\varphi_Q\le 1$, and, by Lemma \reff{Wadd}, (3), $Q^*\cap K^*\ne\emp$ iff $Q\cap K\ne\emp$, so that
$$
|\tf(x)-f(\omega_{K,\alpha})|\le
\sum\{|c_Q-f(\omega_{K,\alpha})|:~Q\in \TW, Q\cap K\ne\emp\}.
$$
By Lemma \reff{Wadd}, (2), there are at most $N(n)$ cubes $Q\in \TW$ such that $Q\cap K\ne\emp$ so that
$$
|\tf(x)-f(\omega_{K,\alpha})|\le
C\,\max\{|c_Q-f(\omega_{K,\alpha})|:~Q\in \TW, Q\cap K\ne\emp\}.
$$
Moreover, for every $Q\in \TW$, $Q\cap K\ne\emp$, by Lemma \reff{Wadd}, (1),
$$
\diam Q\le 4\diam K\le 4(\sigma/4)=\sigma,
$$
so that, by \rf{CQ}, $c_Q=f(\omega_{Q,\alpha})$. The lemma is proved.\bx
\smallskip
\par Observe that $\tf|_\Omega\in C^{\infty}(\Omega)$. The next lemma provides an estimate of the norm of gradient of $\tf$ on every Whitney cube $K\in\TW$.
\begin{lemma}\lbl{E-GR} Let $K\in\TW$ be a Whitney cube. Then
$$
\sup_K\|\nabla\tf\|\le C(n)\, (\diam K)^{-1} \sum\{|c_Q-c_K|:~Q\in\TW, ~Q\cap K\ne\emptyset\}.
$$
\end{lemma}
\par {\it Proof.} For every $x\in K$ we have
$$ \|\nabla\tf(x)\|=
\|\nabla\left(\sum_{Q\in\TW}(c_Q-c_K)\varphi_Q(x)\right)\| =
\|\sum_{Q\in\TW}
(c_Q-c_K)\nabla\varphi_Q(x)\|.
$$
Since $\supp \varphi_Q\subset Q^*,~Q\in\TW$, and $x\in K$, in the latter sum one can consider only those cubes $Q\in\TW$ for which $Q^*\cap K\ne\emptyset$. Hence
$$ \|\nabla\tf(x)\|\le
\sum\{
|c_Q-c_K|\|\nabla\varphi_Q(x)\|:~Q\in\TW, Q^*\cap K\ne\emptyset\}
$$
so that, by Lemma \reff{PU}, (d),
$$
\|\nabla\tf(x)\|\le C(n)\,
\sum\{|c_Q-c_K|(\diam Q)^{-1}:~Q\in\TW, Q^*\cap K \ne\emptyset\}.
$$
By Lemma \reff{Wadd}, (3), $Q\cap K\ne\emp$ provided
$Q^*\cap K\ne\emp$. Moreover,  by \rf{DM}, $\diam Q\sim \diam K$ for every cube $Q\in\TW$ such that $Q^*\cap K\ne\emp$. Hence
$$
\|\nabla\tf(x)\|\le C(n)\, (\diam K)^{-1} \sum\{|c_Q-c_K|:~Q\in\TW, ~Q\cap K\ne\emptyset\},~~~~x\in K,
$$
proving the lemma.\bx
\par Let us estimate the $L_p$-norm of $\nabla\tf$.
\begin{lemma}\lbl{E-GRLP}
$$
\|\nabla\tf\|^p_{L_p(\Omega)}\le C(n)\,\sum\left\{\frac {|c_Q-c_{Q'}|^p}{(\diam Q+\diam Q')^{p-n}}:~Q,Q'\in\TW, ~Q\cap Q'\ne\emptyset\right\}.
$$
\end{lemma}
\par {\it Proof.} We have
$$
\intl_{\Omega}\|\nabla\tf\|^p\,dx\le\sum_{K\in\TW}
\intl_{K}\|\nabla\tf\|^p\,dx\le\sum_{K\in\TW}
|K|\sup_K\|\nabla\tf\|^p
$$
so that, by Lemma \reff{E-GR},
\be
\intl_{\Omega}\|\nabla\tf\|^p\,dx
&\le& C\sum_{K\in\TW}|K|(\diam K)^{-p} \sum\{|c_Q-c_K|^p:Q\in\TW,Q\cap K\ne\emptyset\}\nn\\
&=&
C\sum_{K\in\TW}\sum \left\{\frac{|c_Q-c_K|^p}{(\diam K)^{p-n}}:Q\in\TW, Q\cap K\ne\emptyset\right\}.\nn
\ee
By \rf{DM}, $\diam K\sim\diam Q$ for every $K,Q\in\TW$,
$Q\cap K\ne\emptyset$, so that
\be
\intl_{\Omega}\|\nabla\tf\|^p\,dx
&\le&
C\sum_{K\in\TW}\sum \left\{\frac{|c_Q-c_K|^p}{(\diam K+\diam Q)^{p-n}}:Q\in\TW,Q\cap K\ne\emptyset\right\}\nn\\
&=&
2C\,\sum\left\{\frac {|c_Q-c_{Q'}|^p}{(\diam Q+\diam Q')^{p-n}}:~Q,Q'\in\TW, ~Q\cap Q'\ne\emptyset\right\}.\nn
\ee
\par The lemma is proved.\bx
\begin{lemma}\lbl{WP-NORM} Let $\bc:=0$, see formula \rf{CQ}. Then the following inequality
\be
&&\|\tf\|^p_{\WPO}\nn\\&\le& C\,\left(\sum\left\{\frac {|f(\omega_{Q,\alpha})-f(\omega_{Q',\alpha})|^p}{(\diam Q+\diam Q')^{p-n}}:Q,Q'\in\TW, Q\cap Q'\ne\emptyset,\,Q,Q'\subset\Oc_{5\sigma}(\DO) \right\}\right.\nn\\
&+&
\left.\sum\left\{|f(\omega_{Q,\alpha})|^p|Q|:~Q\in\TW, Q\subset\Oc_{10\sigma}(\DO)\right\}\right)\nn
\ee
holds. Here $C=C(n,\sigma)$ is a constant depending only on $n$ and $\sigma$.
\end{lemma}
\par {\it Proof.} Recall that the set $\DOEP$ is defined by \rf{OC}. By Lemma \reff{E-GRLP},
$$\|\nabla\tf\|^p_{L_p(\Omega)}\le C(I_1+I_2)$$
where
\be
I_1&:=&\sum\{|c_{Q}-c_{Q'}|^p(\diam Q+\diam Q')^{n-p}:\nn\\
&&
Q,Q'\in\TW,~ Q\cap Q'\ne\emptyset,~
\max\{\diam Q,\diam Q'\}\le \sigma\}\nn
\ee
and
\be
I_2&:=&\sum\{|c_{Q}-c_{Q'}|^p(\diam Q+\diam Q')^{n-p}:\nn\\
&&
Q,Q'\in\TW,~ Q\cap Q'\ne\emptyset,~
\max\{\diam Q,\diam Q'\}>\sigma\}.
\nn
\ee
\par Let us estimate $I_1$. Let $Q\in\TW$ and let $\diam Q\le \sigma$. Then, by \rf{DKQ}, for every $x\in Q$ we have
$$
\dist(x,\DO)\le\diam Q+\dist(Q,\DO)\le\diam Q+4\diam Q\le 5\sigma
$$
so that $Q\subset\Oc_{5\sigma}(\DO)$. By \rf{CQ},
$$
|c_{Q}-c_{Q'}|=|f(\omega_{Q,\alpha})-
f(\omega_{Q',\alpha})|
$$
provided $\diam Q\le\sigma$ and $\diam Q'\le \sigma$ so that
$$
I_1\le
\sum\left\{\frac {|f(\omega_{Q,\alpha})-f(\omega_{Q',\alpha})|^p}{(\diam Q+\diam Q')^{p-n}}:Q,Q'\in\TW, Q\cap Q'\ne\emptyset, Q,Q'\subset\Oc_{5\sigma}(\DO)\right\}.
$$
\par Let us estimate $I_2$. Let $Q\in\TW$ and let $Q\nsubseteq\Oc_{10\sigma}(\DO)$ so that there exists $x\in Q$ such that $\dist(x,\DO)>10\sigma.$ Then for every $Q'\in\TW,~Q'\cap Q\ne\emptyset,$ we have
$$
10\sigma<\dist(x,\DO)\le\diam Q+\diam Q'+\dist(Q',\DO)
$$
so that, by \rf{DKQ}  and \rf{DM},
$$
10\sigma<4\diam Q'+\diam Q'+4\diam Q'=9\diam Q'.
$$
Thus $\sigma<\diam Q'$ for every $Q'\in\TW,~Q'\cap Q\ne\emptyset,$ so that, by \rf{CQ},
\bel{QF}
c_{Q'}=0,~~~~Q'\in\TW,~Q'\cap Q\ne\emptyset,
\ee
provided $Q\nsubseteq\Oc_{10\sigma}(\DO)$. Hence
\be
I_2&:=&\sum\{|c_{Q}-c_{Q'}|^p(\diam Q+\diam Q')^{n-p}:\nn\\
&&
Q,Q'\in\TW,~ Q\cap Q'\ne\emptyset,~
\max\{\diam Q,\diam Q'\}>\sigma,Q,Q'\subset\Oc_{10\sigma}(\DO)\}.\nn
\ee
\par Observe that for each $Q\subset\Oc_{10\sigma}(\DO)$, by \rf{DKQ}, we have
$$
\diam Q\le\dist(Q,\DO)\le 10\sigma.
$$
\par Let us consider cubes $Q,Q'\in\TW$ satisfying the following conditions: $Q,Q'\subset\Oc_{10\sigma}(\DO)$, $Q\cap Q'\ne\emptyset$ and $\max\{\diam Q,\diam Q'\}>\sigma.$ Assume that $\diam Q>\sigma$. Then, by \rf{DM},
$$
\diam Q'\ge\tfrac14\diam Q\ge\sigma/4
$$
so that
$$
\sigma/4<\diam Q'\le\diam Q\le 10\sigma.
$$
Since in \rf{CQ} we put $\bc:=0$,
$$
|c_Q-c_{Q'}|\le |c_Q|+|c_{Q'}|\le |f(\omega_{Q,\alpha})|+|f(\omega_{Q',\alpha})|,
$$
so that
\be
I_2&\le&
C\sum\{(|f(\omega_{Q,\alpha})|^p+ |f(\omega_{Q',\alpha})|^p)\,(\diam Q+\diam Q')^{n-p} :~Q,Q'\in\TW,\nn\\
&&
Q\cap Q'\ne\emptyset,~Q,Q'\subset\Oc_{10\sigma}(\DO),~
\sigma/4<\diam Q',\diam Q\le 10\sigma\}.
\nn
\ee
Hence
$$
I_2\le
C(n,\sigma)\sum \{|f(\omega_{Q,\alpha})|^p\,|Q|:Q\in\TW,
Q\subset\Oc_{10\sigma}(\DO)\}.
$$
We obtain
\be
&&\|\nabla\tf\|^p_{L_p(\Omega)}\le C(I_1+I_2)\nn\\
&\le&
C\,\left(\sum\left\{\frac {|f(\omega_{Q,\alpha})-f(\omega_{Q',\alpha})|^p}{(\diam Q+\diam Q')^{p-n}}:Q,Q'\in\TW, Q\cap Q'\ne\emptyset,\,Q,Q'\subset\Oc_{5\sigma}(\DO) \right\}\right.\nn\\
&+&
\left.\sum \{|f(\omega_{Q,\alpha})|^p\,|Q|:Q\in\TW,
Q\subset\Oc_{10\sigma}(\DO)\}\right)\nn
\ee
with $C=C(n,\sigma)$.
\par Let us estimate $\|\tf\|_{L_p(\Omega)}$. By \rf{QF}, $c_Q=0$ provided $Q\nsubseteq\Oc_{10\sigma}(\DO)$. Also, since $\bc=0$, by \rf{CQ}, $|c_Q|\le|f(\omega_{Q,\alpha})|$ for every $Q\in\TW$. Hence
\be
\|\tf\|_{L_p(\Omega)}^p&\le&\sum_{Q\in\TW}|c_Q|^p\,|Q|\le
\sum\{|c_Q|^p\,|Q|:~Q\in\TW,Q\subset\Oc_{10\sigma}(\DO)\}
\nn\\
&\le&
\sum\{|f(\omega_{Q,\alpha})|^p\,|Q|: ~Q\in\TW,Q\subset\Oc_{10\sigma}(\DO)\}.\nn
\ee
\par The lemma is completely proved.\bx
\par Let us extend the notion of  $Q$-visibility (Definition \reff{Q-V}) to the case of an arbitrary domain $\Omega\subset\RN$.
\begin{definition}\lbl{Q-VA} {\em Let $\omega\in\DOA,$ $\alpha\in(0,1],$ and let $Q\subset\Omega$ be a cube.  We say that {\it $\omega$ is $(\alpha,Q)$-visible in $\Omega $} if the following conditions are satisfied:
\par (i). The point  $\ell(\omega)\in\DO$ is $Q$-visible in $\Omega$ (see Remark \rf{R-CPM} and Definition \reff{Q-V});
\par (ii).  There exists a sequence $y_i\in(\ell(\omega),x_Q],$ $i=1,2,...~,$ such that the equivalence class of $(y_i)$ with respect to $``\sa"$, see \rf{ESIM} and \rf{KE}, coincides with $\omega$:
$$
[(y_i)]_\alpha=\omega.
$$
}
\end{definition}
\par Observe that part (ii) of this definition can be replaced with one the following equivalent statement:
\par {\it (a). There exists a sequence $(y_i)\in\omega$ such that $y_i\in(\ell(\omega),x_Q]$.
\par (b). Every sequence $y_i\in(\ell(\omega),x_Q]$ such that $y_i\stackrel{\EN}{\longrightarrow} \eo$ as $i\to\infty,$ belongs to $\omega$, see Lemma \reff{DS-QV}.}
\par Also note that for each cube $Q\subset\Omega$ the element $\oqa\in\DOA$ defined by \rf{DYO} possesses the following property:
\bel{QVO}
\oqa~~\text{{\it is}}~~(\alpha,Q)\text{{\it -visible in}}~~\Omega.
\ee
\begin{lemma}\lbl{TR-SQ} Let $Q_1,Q_2$ be cubes in $\Omega$. Suppose that $Q_2\in\TW$, $Q_1\cap Q_2\ne\emptyset$, and for some $\tau\ge 1$
$$
\diam Q_1\le\diam Q_2\le \tau \diam Q_1.
$$
\par Then $\omega_{Q_2,\alpha}$ is $(\alpha,Q_1)$-visible in $\Omega$. In addition, $\ell(\omega_{Q_2,\alpha})\in (10\tau +1)Q_1$.
\end{lemma}
\par {\it Proof.} Prove that $\omega_{Q_2,\alpha}$ is $(\alpha,Q_1)$-visible in $\Omega$. Observe that
$\omega_{Q_2,\alpha}$ is $(\alpha,Q_2)$-visible in $\Omega$ and $a_{Q_2}=\ell(\omega_{Q_2,\alpha})$, see \rf{QVO}. (Recall that $a_{Q_2}$ denotes a point on $\DO$ nearest to the cube $Q_2$.) Thus $\omega_{Q_2,\alpha}=[(y_i)]_\alpha$ where  $(y_i)$ is an arbitrary sequence of points such that $y_i\in(a_{Q_2},x_{Q_2}]$, $i=1,2,...~,$ and
\bel{YTU}
\|y_i-a_{Q_2}\|\to 0~~\text{as}~~~i\to\infty.
\ee
\par Let us fix such a sequence $(y_i)$ and construct a sequence $(z_i)$ such that $z_i\in(a_{Q_2},x_{Q_1}]$ for every $ i=1,2,...~,$ and $(z_i)\sa(y_i)$. Since $Q_2\in\TW$, $Q_1\cap Q_2\ne\emptyset$ and
$$
\diam Q_1\le\diam Q_2\le \dist(Q_2,\DO),
$$
by Lemma \reff{L-QV}, $a_{Q_2}$ is $Q_1$-visible in $\Omega$ and for every $i=1,2,...~,$ there exists a point $z_i\in(a_{Q_2},x_{Q_1}]$ such that $\|y_i-z_i\|\le C\|y_i-a_{Q_2}\|,$ and
$$
\dao(y_i,z_i)\le C\|y_i-a_{Q_2}\|^\alpha
$$
where $C$ is a constant independent of $i$. Hence,
$$
\dt(y_i,z_i)=\|y_i-z_i\|^\alpha+\dao(y_i,z_i)\le C\|y_i-a_{Q_2}\|^\alpha
$$
so that, by \rf{YTU}, $\dt(y_i,z_i)\to 0$ as $i\to\infty.$
\par This proves that $(y_i)\sa(z_i)$, see \rf{ESIM}.
Since $\omega_{Q_2,\alpha}$ is the equivalence class of $(y_i)$, and $(y_i)\sa(z_i)$, we conclude that $\omega_{Q_2,\alpha}=[z_i]_\alpha$ as well.
\par Thus the point $\ell(\omega_{Q_2,\alpha})=a_{Q_2}$ is $Q_1$-visible, $(z_i)\in\omega_{Q_2,\alpha}$, and $z_i\in(\ell(\omega),x_{Q_1}]$, so that, by Definition \reff{Q-VA}, $\omega_{Q_2,\alpha}$ is $(\alpha,Q_1)$-visible.
\par Prove that $a_{Q_2}\in (10\tau +1)Q_1$. Let $Q_1=Q(x_{Q_1},r_1)$. Since $Q_1\cap Q_2\ne\emp$,
$$
\|a_{Q_2}-x_1\|\le \dist(a_{Q_2},Q_2)+\diam Q_2+r_1=
\dist(Q_2,\DO)+\diam Q_2+r_1
$$
so that, by \rf{DKQ},
$$
\|a_{Q_2}-x_1\|\le 4\diam Q_2+\diam Q_2+r_1.
$$
But  $\diam Q_2\le \tau\diam Q_1=2\tau r_1$ so that
$
\|a_{Q_2}-x_1\|\le (10\tau +1) r_1
$
proving the required property $\ell(\omega_{Q_2,\alpha})=a_{Q_2}\in (10\tau +1)Q_1$.\bx
\par This lemma and the property \rf{DM} of Whitney's cubes imply the following
\begin{lemma}\lbl{OM-VIS} Let $Q_1,Q_2\in\TW$, $Q_1\cap Q_2\ne\emptyset$, and let $\diam Q_2\ge\diam Q_1$. Then $\omega_{Q_1,\alpha}$ and $\omega_{Q_2,\alpha}$ are $(\alpha,Q_1)$-visible. In addition, $\ell(\omega_{Q_1,\alpha}),\ell(\omega_{Q_2,\alpha})\in 41Q_1$.
\end{lemma}
\par The last auxiliary result of the section is the following
\begin{lemma}\lbl{INC} Let cubes $Q,Q'\subset\Omega$ and let $Q'\subset Q$. Then $\omega_{Q,\alpha}$ is an $(\alpha,Q')$-visible in $\Omega$.
\end{lemma}
\par The proof of this result relies on Lemma \reff{L-QV} and is very similar to the proof of Lemma \reff{TR-SQ}. We leave the details to the interested reader.
\section*{{\normalsize 5. Boundary values of Sobolev functions: restrictions and extensions.}}
\setcounter{section}{5}
\setcounter{theorem}{0}
\setcounter{equation}{0}
\par In subsection 1.4 we have formulated Theorem \reff{VC-D} and Theorem \reff{VC-WD} which provide
constructive descriptions of the trace spaces   $\LOPO|_{\DO}$ and $\WPO|_{\DO}$ whenever $\Omega$ is a domain in $\RN$ satisfying the condition $A_\alpha$ with $\alpha=(p-n)/(p-1)$, see Definition \reff{A-C}. In this section we present generalizations of these theorems to the case of {\it an arbitrary} domain $\Omega\subset\RN$.
\par We will be needed several addition definitions and notation. First of them is a definition of the metric $\dco$ introduced earlier only for the domains satisfying
the condition $A_\alpha$, see \rf{D-AO}. Given $\alpha\in(0,1]$ and $\omega_1,\omega_2\in\OA$ we put
\bel {D-CON}
\dco(\omega_1,\omega_2):=\lim_{i\to\infty}\dao(x_i,y_i)
\ee
where $(x_i)\in\omega_1,(y_i)\in\omega_2$ are arbitrary sequences. Recall that $\OA$ is the family of all equivalence classes of Cauchy sequences with respect to the metric $\dt$, see \rf{CMP-OM}. Since $\dao\le\dt$, see \rf{DEF-TDA}, $\dco$ is well defined on $\OA$.
\par  Clearly, $\dco=\dao$ on $\Omega$.  Moreover, the reader can easily see that $\dco$ coincides with the metric of the Cauchy completion of the metric space $(\Omega,\dao)$. Now equality \rf{CDR-L} can be written in the following form: for every $\omega_1,\omega_2\in\OA$
\bel{F-RD}
\dtc(\omega_1,\omega_2)=\dco(\omega_1,\omega_2) +\|\ell(\omega_1)-\ell(\omega_2)\| ^{\alpha}.
\ee
\par In turn, this equality and \rf{H-LOC} yield:
Let $f\in \tra(\LOPO)$ and let $\omega_1,\omega_2\in\DOA$, see Definition \reff{B-DA} and \rf{B-A}. Then
\bel{H-DOM} |f(\omega_1)-f(\omega_2)|\le C\|f\|_{\tra(\LOPO)}
\{\dco(\omega_1,\omega_2)^{1-\frac{1}{p}} +\|\ell(\omega_1)-\ell(\omega_2)\| ^{1-\frac{n}{p}}\}
\ee
where $\alpha =(p-n)/(p-1)$ and $C=C(n,p)$.
Observe that this inequality is equivalent to the following one:
\bel{H-ROM} |f(\omega_1)-f(\omega_2)|\le C\|f\|_{\tra(\LOPO)}\,\dtc(\omega_1,\omega_2) ^{1-\frac{1}{p}},~~~\omega_1,\omega_2\in\OA.
\ee
\par By Lemma \reff{CMP}, $\dt(x,y)\sim\dao(x,y)$ provided
$x,y\in\Omega$ and
$$
\|x-y\|\ge\min\{\D{x},\D{y}\}.
$$
This implies the following inequality:
\bel{EQ-GR1}
\dtc(\omega_1,\omega_2)\sim\dco(\omega_1,\omega_2), ~~~~\omega_1\in\DOA,\,\omega_2\in\OA.
\ee
In particular,
\bel{U-0}
\dtc(\omega_1,\omega_2)\sim\dco(\omega_1,\omega_2)~~~\text{for every}~~\omega_1,\omega_2\in\DOA,
\ee
and
$$
\dtc(\omega,x)\sim \dco(\omega,x)~~~\text{for every}~~\omega\in\DOA,\,x\in\Omega.
$$
Combining equivalence \rf{U-0} with inequality \rf{H-ROM}, we obtain
$$
|f(\omega_1)-f(\omega_2)|\le C\|f\|_{\tra(\LOPO)}\,\dco(\omega_1,\omega_2) ^{1-\frac{1}{p}},~~~~\omega_1,\omega_2\in\DOA.
$$
Thus {\it every function $f\in \tra (\LOPO)$ is continuous with respect to the metric $\dco$.}
\par Observe also that \rf{EQ-GR1} and \rf{TR-LM} provide the following definition of the trace to $\DOA$: Let $F\in\WPO$ and let $f=\tra F$. Then for every $\omega\in\DOA$ we have
\bel{TR-DLM}
f(\omega)=\lim\{F(x):~\dco(x,\omega)\to 0, x\in\Omega\}.
\ee
\begin{theorem}\lbl{VS-GN} Let $\Omega$ be a domain in
$\RN$. Let $p\in(n,\infty)$ and let $\alpha=\pn$. Let $\eta$ be a constant satisfying $\eta \ge 41$.
\par A function $f:\DOA\to\R$ is the trace to $\DOA$ of a (continuous) function $F\in L_{p}^{1}(\Omega )$ if and only if $f$ is continuous with respect to $\dco$ and there exists a constant $\lambda>0$ such that for every finite family $\{Q_{i}:i=1,...,m\}$ of pairwise disjoint cubes in $\Omega$ and every choice of $(\alpha,Q_{i})$-visible elements $\omega_i^{(1)},\omega_i^{(2)}\in\DOA$ such that
\bel{O1O2-41Q}
\ell(\omega_i^{(1)}),\ell(\omega_i^{(2)})\in (\eta Q_{i}) \cap \DO,
\ee
the following inequality
\bel{VAR-IN}
\sum_{i=1}^{m}\frac{|f(\omega_i^{(1)})-f(\omega_i^{(2)})|^p}
{(\diam Q_{i})^{p-n}}\le \lambda
\ee
holds. Moreover,
$$
\|f\|_{\tra(\LOPO)}\sim \inf\lambda ^{\frac{1}{p}}
$$
with constants of equivalence depending only on $n$, $p$ and $\eta$.
\end{theorem}
\par {\it Proof. (Necessity).} Let $F\in\LOPO$. Prove that the function $f:=\tra F$ satisfies the theorem's conditions. As we have seen above, $f$ is  a continuous function on $\DOA$ with respect to the metric $\dco$. Prove that $f$ satisfies inequality \rf{VAR-IN}.
\par We will be needed an auxiliary lemma. Given a function $g$ defined on $\Omega$ we let $g^\cw$ denote its extension by zero to all of $\RN$. Thus $g^\cw(x):=g(x), x\in\Omega,$ and $g^\cw(x):=0, x\notin\Omega.$ Also, given $q>0$ we put
\bel{DGQ1}
G(x):=(\|\nabla F\|^q)^\cw(x), ~~~~~x\in\RN.
\ee
\begin{lemma}\lbl{NG} Let $q\in(n,p]$ and let $\eta>1$. Let $Q=Q(x_Q,r_Q)$ be a cube in $\Omega$ and let $\omega\in\DOA$. Suppose that $\eo\in(\eta Q)\cap \DO$ and $\omega$ is $(\alpha,Q)$-visible. Then
\bel{F-OM}
\frac{|f(\omega)-F(x_Q)|^p}{(\diam Q)^{p-n}}\le C\,
\intl_{Q} \Mc[G]^{\frac{p}{q}}(z)\,dz
\ee
where $C=C(n,q,\eta)$.
\end{lemma}
\par {\it Proof.} Since $\omega$ is $(\alpha,Q)$-visible,
$$
\Conv\{\eo,Q\}\setminus \{\eo\}\subset\Omega,
$$
see Definition \reff{Q-V} and Definition \reff{Q-VA}. Let $y\in(\eo,x_Q]$. Clearly, $y$ is also $Q$-visible point of $\Omega$ and $y\in(\eta Q)\cap \DO$. Then, by Lemma \reff{E-QVX} with $t=\eta$,
$$
\left(\frac{|F(y)-F(x_Q)|}{\diam Q}\right)^q\le C\,\frac{1}{|Q|}\intl_{\gamma Q\cap \Omega}
\|\nabla F(z)\|^q\,dz
$$
where $C=C(n,q,\eta)$ and $\gamma=\gamma(n,q,\eta)$.
\par Applying Lemma \reff{MF} with $\theta=p/q$, we obtain
$$
\left(\frac{1}{|Q|}\intl_{\gamma Q\cap \Omega}
\|\nabla F(z)\|^q\,dz\right)^{\frac{p}{q}}=\left(\frac{1}{|Q|} \intl_{\gamma Q\cap \Omega} G(z)\,dz\right)^{\frac{p}{q}}\le C\,\frac{1}{|Q|}\intl_{Q} \Mc[G]^{\frac{p}{q}}(z)\,dz.
$$
Hence
$$
\left(\frac{|F(y)-F(x_Q)|}{\diam Q}\right)^p\le  C\,\frac{1}{|Q|}\intl_{Q} \Mc[G]^{\frac{p}{q}}(z)\,dz
$$
so that
\bel{F-YX}
\frac{|F(y)-F(x_Q)|^p}{(\diam Q)^{p-n}}\le C\,
\intl_{Q} \Mc[G]^{\frac{p}{q}}(z)\,dz.
\ee
\par Now let $y_i\in(\eo,x_Q]$ and let $(y_i)$ tends to $\eo$ in the Euclidean norm. Then, by Lemma \reff{DS-QV}, (ii), the sequence $(y_i)\in\omega$. Therefore, by \rf{TR-S},
$$
f(\omega)=(\tra F)(\omega)=\lim_{i\to\infty} F(y_i).
$$
But, by \rf{F-YX},
$$
\frac{|F(y_i)-F(x_Q)|^p}{(\diam Q)^{p-n}}\le C\,
\intl_{Q} \Mc[G]^{\frac{p}{q}}(z)\,dz,
$$
so that, letting $i$  tend to $\infty$, we obtain \rf{F-OM}.\bx
\par Using this lemma, we prove the necessity as follows. Put $q:=(n+p)/2$. Now, let $\omega^{(1)},\omega^{(2)}\in\DOA$ be  $(\alpha,Q)$-visible elements such that
$$
\ell(\omega^{(1)}),\ell(\omega^{(2)})\in (\eta Q)\cap \DO.
$$
Then, by \rf{F-OM},
$$
\frac{|f(\omega^{(1)})-f(\omega^{(2)})|^p}{(\diam Q)^{p-n}}\le 2^p\left(
\frac{|f(\omega^{(1)})-F(x_Q)|^p}{(\diam Q)^{p-n}}+
\frac{|f(\omega^{(2)})-F(x_Q)|^p}{(\diam Q)^{p-n}}\right)\nn\\
\le C\,\intl_{Q} \Mc[G]^{\frac{p}{q}}(z)\,dz.
$$
\par Finally, let $\omega_i^{(1)},\omega_i^{(2)}\in\DOA$, $i=1,...,m,$ and let  $\{Q_{i}:i=1,...,m\}$ be a family of pairwise disjoint cubes in $\Omega$ such that $\omega_i^{(1)},\omega_i^{(2)}$ are $(\alpha,Q_i)$-visible and
$$
\ell(\omega_i^{(1)}),\ell(\omega_i^{(2)})\in (\eta Q_{i}) \cap \DO
$$
for every $i=1,...,m$. Then
$$
I:=\sum_{i=1}^{m}
\frac{|f(\omega_i^{(1)})-f(\omega_i^{(2)})|^p}
{(\diam Q_{i})^{p-n}}\le C\,\sum_{i=1}^{m}
\intl_{Q_i} \Mc[G]^{\frac{p}{q}}(z)\,dz=C\,\intl_{U} \Mc[G]^{\frac{p}{q}}(z)\,dz
$$
where $U:=\cup_{i=1}^m Q_i$. Hence
\bel{FEI}
I\le C\,\intl_{\RN}\Mc[G]^{\frac{p}{q}}(z)\,dz.
\ee
Since $p/q>1$, by the Hardy-Littlewood maximal theorem,
$$
\intl_{\RN}\Mc[G]^{\frac{p}{q}}(z)\,dz\le C\,\intl_{\RN}|G|^{\frac{p}{q}}(z)\,dz
=C\,\intl_{\RN} ((\|\nabla F\|^q)^\cw)^{\frac{p}{q}}(z)\,dz=C\,\intl_{\Omega}(\|\nabla F\|^q)^{\frac{p}{q}}(z)\,dz
$$
so that
\bel{E-MG}
\intl_{\RN}\Mc[G]^{\frac{p}{q}}(z)\,dz\le
C\,\|\nabla F\|^p_{L_p(\Omega)}.
\ee
Hence $I\le C\,\|\nabla F\|^p_{L_p(\Omega)}.$ Taking in this inequality the infimum  over all functions $F\in\LOPO$ such that $f=\tra F$ we obtain the required inequality
$I^{\frac1p}\le C\,\|f\|_{\tra(\LOPO)}.$
\par The proof of the necessity is finished.
\bigskip
\par {\it (Sufficiency.)} Let $f:\DOA\to\R$ be a continuous function with respect to the metric $\dco$, see \rf{D-CON}. Let $\lambda$ be a positive constant such that
for every finite family $\{Q_{i}:i=1,...,m\}$ of pairwise disjoint cubes in $\Omega$ and every $(\alpha,Q_i)$-visible elements $\omega_i^{(1)},\omega_i^{(2)}\in\DOA$ satisfying \rf{O1O2-41Q} the inequality \rf{VAR-IN} holds.
\par We put $\sigma=+\infty$ in formula \rf{CQ}; thus $c_Q=f(\oqa)$ for {\it every} $Q\in\TW$. Then we define a function $\tf:\Omega\to\R$ by formula \rf{ExtOp}. Thus
\bel{ExtOp-L}
\tf(x):=\sum\limits_{Q\in \TW}
f(\oqa)\varphi_Q(x), ~~~ x\in \Omega.
\ee
\par Prove that $\tf\in\LOPO$ and $\|\tf\|_{\LOPO}\le C\lambda^{\frac1p}$. Since $c_Q=f(\oqa)$ for every $Q\in\TW$, by Lemma \reff{E-GRLP},
$$
\|\nabla\tf\|^p_{L_p(\Omega)}\le C\, V_p(f;\Omega)
$$
where
\bel{VP-OM}
V_p(f;\Omega):=\sum\left\{\frac {|f(\oqa)-f(\oqp)|^p}{(\diam Q+\diam Q')^{p-n}}: ~Q,Q'\in\TW, ~Q\cap Q'\ne\emptyset\right\}.
\ee
Hence
$$
V_p(f;\Omega)\le \sum_{Q\in\TW}\sum_{Q'\in\Ac_Q}\frac {|f(\oqa)-f(\oqp)|^p}{(\diam Q)^{p-n}}.
$$
where
$$
\Ac_Q:=
\{Q'\in\TW:~Q'\cap Q\ne\emptyset, \diam Q'\ge\diam Q \}.
$$
\par Let $K_Q\in\Ac_Q$ be a cube such that
$$
\max_{Q'\in\Ac_Q}|f(\oqa)-f(\oqp)| =|f(\oqa)-f(\omega_{K_Q,\alpha})|.
$$
Since the family $\Ac_Q$ consists of at most $N(n)$ cubes, see Lemma \reff{Wadd}, (2), we have
$$
V_p(f;\Omega)\le C\,\sum_{Q\in\TW}\frac {|f(\oqa)-f(\omega_{K_Q,\alpha})|^p}{(\diam Q)^{p-n}}.
$$
By Theorem \reff{Wcov}, (iii), the family $\TW$ has the covering multiplicity $ M(\TW)\le N(n)$. Therefore this family can be partitioned into at most $N_1=N_1(n)$ families $\{\pi_j:~j=1,...,N_1\}$ of {\it pairwise disjoint cubes}. See \cite{BrK,Dol}.
\par Observe that for each $Q\in\TW$ the cube $K_Q\in\TW$, $K_Q\cap Q\ne\emp$ and $\diam Q\le\diam K_Q$, so that, by Lemma \reff{OM-VIS}, $\oqa$ and $\omega_{K_Q,\alpha}$ are $(\alpha,Q)$-visible and
$$
\ell(\oqa),\ell(\omega_{K_Q,\alpha})\in 41Q.
$$
Hence, by \rf{VAR-IN}, for every $j=1,...,N_1,$ we have
$$
\sum_{Q\in\pi_j}\frac {|f(\oqa)-f(\omega_{K_Q,\alpha})|^p}{(\diam Q)^{p-n}}\le\lambda,
$$
so that
\bel{N-FW}
V_p(f;\Omega)\le C\,\sum_{j=1}^{N_1}\sum_{Q\in\pi_j}\frac {|f(\oqa)-f(\omega_{K_Q,\alpha})|^p}{(\diam Q)^{p-n}}\le C\,\lambda.
\ee
Hence $\|\nabla\tf\|^p_{L_p(\Omega)}\le C\,\lambda$ proving that $\tf\in\LOPO$. This also proves that the trace $\tra\tf$ is well-defined.
\par Prove that $\tra\tf=f$. Let $\omega\in\DOA$ and let a sequence $(y_i)\in\omega$. Thus
$$
\lim_{i\to\infty}\dtc(y_i,\omega)=0.
$$
By \rf{F-RD},
$$
\dtc(y_i,\omega)=
\dco(y_i,\omega) +\|y_i-\ell(\omega)\| ^{\alpha},
$$
so that $\dco(y_i,\omega)\to 0$ and $\|y_i-\ell(\omega)\| \to 0$ as $i\to\infty.$ Recall that, by \rf{TR-S},
\bel{DT-SF}
\tra \tf(\omega):=\lim_{i\to\infty} \tf(y_i).
\ee
\par We let $K_i\in\TW$ denote a Whitney's cube such that $K_i\ni y_i, i=1,2,...~.$ By Lemma \reff{LO-QA},
$$
\dtc(y_i,\omega_{K_i,\alpha})\le C\,\dist(y_i,\DO)^\alpha.
$$
\par Recall that the element $\omega_{K_i,\alpha}\in\DOA$ is defined as an equivalence class  $\omega_{K_i,\alpha}=[(x_i)]_\alpha$ where
$$
x_i:=a_{K_i}+\tfrac1i(x_{K_i}-a_{K_i}), ~~~i=1,2,...~,
$$
see \rf{XNE} and \rf{DYO}. Since $\ell(\omega)\in\DO,$ we obtain
\bel{DY-DO}
\dist(y_i,\DO)\le \|y_i-\ell(\omega)\|\to 0~~~\text{as}~~~i\to\infty,
\ee
so that
$$
\dco(y_i,\omega_{K_i,\alpha})
\le\dtc(y_i,\omega_{K_i,\alpha})\le C\,\dist(y_i,\DO)^\alpha\to 0~~~\text{as}~~~i\to\infty.
$$
Hence
\bel{DO2}
\dco(\omega_{K_i,\alpha},\omega)\le
\dco(\omega_{K_i,\alpha},y_i)+\dco(y_i,\omega)\to 0~~~\text{as}~~~i\to\infty.
\ee
\par Let us prove that
\bel{FW2}
\lim_{i\to\infty}|f(\omega_{K_i,\alpha})-\tf(y_i)|=0.
\ee
\par Put
$$
I(K_i):=\{Q\in\TW:~Q\cap K_i\ne\emp\}.
$$
By Lemma \reff{C-L},
\bel{M-F0}
|\tf(y_i)-f(\omega_{K_i,\alpha})|\le C\,\max_{Q\in I(K_i)} |f(\omega_{Q,\alpha})-f(\omega_{K_i,\alpha})|.
\ee
On the other hand, by Lemma \reff{OM-VIS}, for every cube $Q\in I(K_i)$ with $\diam Q\ge \diam K_i$ the elements
$\omega_{Q,\alpha}$ and $\omega_{K_i,\alpha}$ are $(\alpha,K_i)$-visible. In addition, $\ell(\omega_{Q,\alpha}),
\ell(\omega_{K_i,\alpha})\in 41 K_i$.
\par Put $\omega_1^{(1)}:=\omega_{Q,\alpha}$, $\omega_1^{(2)}:=\omega_{K_i,\alpha}$ and $Q_1:=K_i$. Then the triple $\omega_1^{(1)},\omega_1^{(2)},\{Q_1\}$ satisfies the conditions of Theorem \reff{VS-GN} (with $m=1$) so that, by the assumption, the inequality \rf{VAR-IN} holds for this triple. By this inequality,
$$
\frac{|f(\omega_1^{(1)})-f(\omega_1^{(2)})|^p}
{(\diam Q_{1})^{p-n}}\le \lambda,
$$
so that
$$
|f(\omega_{Q,\alpha})-f(\omega_{K_i,\alpha})|\le
\lambda^{\frac1p}\,(\diam K_i)^{1-\frac{n}{p}}
$$
provided $Q\in I(K_i)$ and $\diam Q\ge\diam K_i$.
\par If $Q\in I(K_i)$ and $\diam Q<\diam K_i$, in the same way we prove that
$$
|f(\omega_{Q,\alpha})-f(\omega_{K_i,\alpha})|\le
\lambda^{\frac1p}\,(\diam Q)^{1-\frac{n}{p}}.
$$
\par But $\diam Q\le 4\diam K_i$ for every $Q\in I(K_i),$ see \rf{DM}, so that
$$
|f(\omega_{Q,\alpha})-f(\omega_{K_i,\alpha})|\le
C\,\lambda^{\frac1p}\,(\diam K_i)^{1-\frac{n}{p}}~~\text{ for every}~~Q\in I(K_i).
$$
Hence, by \rf{M-F0},
$$
|\tf(y_i)-f(\omega_{K_i,\alpha})|\le C\,\lambda^{\frac1p}\,(\diam K_i)^{1-\frac{n}{p}}.
$$
Since $K_i\in\TW$ and $y_i\in K_i$, by \rf{DKQ},
$$
\diam K_i\le\dist(K_i,\DO)\le\dist(y_i,\DO)
$$
so that
$$
|\tf(y_i)-f(\omega_{K_i,\alpha})|\le C\,\lambda^{\frac1p}\,\dist(y_i,\DO)^{1-\frac{n}{p}}.
$$
But, by \rf{DY-DO}, $\dist(y_i,\DO)\to 0$ as $i\to\infty$, proving \rf{FW2}.
\par It remains to note that the function $f:\DOA\to\R$ is continuous with respect to the metric $\dco$ so that, by \rf{DO2},
$$
\lim_{i\to\infty}f(\omega_{K_i,\alpha})=f(\omega).
$$
Combining this equality with \rf{FW2}, we conclude that
$$
\lim_{i\to\infty}\tf(y_i)=f(\omega)
$$
so that, by \rf{DT-SF}, $\tra \tf(\omega)=f(\omega)$.
\par Theorem \reff{VS-GN} is completely proved.\bx
\par Our next result, Theorem \reff{VCR-AO}, extends the trace criterion for the Sobolev space given in Theorem \reff{VC-WD} to the case of an {\it arbitrary domain in $\RN$.}
\par Let $\theta>1$ and let $\Qc=\{Q\}$ be a covering of $\Omega$ by non-overlapping cubes such that
\bel{TH-WTN}
\frac{1}{\theta}\diam Q\le\dist(Q,\DO)\le{\theta}\diam Q.
\ee
By $T_{\Qc}:\Omega\to\DOA$ we denote a mapping defined by the following formula:
\bel{DF-T}
T_{\Qc}|_Q:=\omega_{Q,\alpha},~~~~Q\in\Qc.
\ee
Recall that the element $\omega_{Q,\alpha}$ is defined by equalities \rf{XNE} and \rf{DYO}.
\begin{theorem}\lbl{VCR-AO} Let $\Omega$ be a domain in
$\RN$ and let $p\in(n,\infty)$. Fix constants $\ve>0$, $\theta>1$, $\eta\le 22\theta^2$, and an arbitrary covering $\Qc$ of $\Omega$ consisting of
non-overlapping cubes $Q\subset\Omega$ each satisfying inequality \rf{TH-WTN}.
\par A function $f:\DOA\to\R$ is an element of\,\,
$\tra(\WPO)$ if and only if $f$ is continuous with respect to the metric $\dco$, the function $f\circ T_{\Qc}\in L_{p}(\DOEP)$, and there exists a constant $\lambda>0$ such that for every finite family $\{Q_{i}:i=1,...,m\}$ of pairwise disjoint cubes contained in
$\DOEP$ and every choice of $(\alpha,Q_i)$-visible elements $\omega_i^{(1)},\omega_i^{(2)}\in\DOA$ such that
\bel{WP-THQ}
\ell(\omega_i^{(1)}),\ell(\omega_i^{(2)})\in (\eta Q_{i}) \cap \DO,
\ee
the following inequality
\bel{WPV-IN}
\sum_{i=1}^{m}\frac{|f(\omega_i^{(1)})-f(\omega_i^{(2)})|^p}
{(\diam Q_{i})^{p-n}}\le \lambda
\ee
holds. Moreover,
\bel{NW-LO}
\|f\|_{\tra(\WPO)}\sim \|f\circ T_{\Qc}\|_{L_p(\DOEP)}+\inf\lambda ^{\frac{1}{p}}
\ee
with constants of equivalence depending only on $n,p,
\ve,\theta$ and $\eta$.
\end{theorem}
\par {\it Proof. (Necessity).} Let $F\in\WPO$ and let $f=\tra F$. Since $F\in\LOPO$, the function  $f$ is continuous with respect to the metric $\dco$. In turn, we obtain inequality \rf{WPV-IN} by repeating the proof of the necessity part of Theorem \reff{VS-GN}.
\par Thus it remains to show that $f\circ T_{\Qc}\in L_{p}(\DOEP)$ and
$$
\|f\circ T_{\Qc}\|_{L_p(\DOEP)}\le C\|F\|_{\WPO}.
$$
\par We put
$$
\WQ:=\{Q\in\Qc:~\diam Q\le\theta\ve\}.
$$
Let $Q\in\Qc$ be a cube such that $Q\cap \DOEP\ne\emp$. Then $\dist(Q,\DO)\le\ve$ so that
$$
\diam Q\le \theta\dist(Q,\DO)\le\theta\ve
$$
proving that
$$
\DOEP\subset \cup\{Q:~Q\in\WQ\}.
$$
\par Hence
$$
\|f\circ T_{\Qc}\|_{L_p(\DOEP)}^p\le\sum_{Q\in\WQ}\,\,\intl_Q(f\circ T_{\Qc})^p(x)\,dx= \sum_{Q\in\WQ}|Q|\,|f(\omega_{Q,\alpha})|^p.
$$
\par We let $F_Q:=|Q|^{-1}\int_QF\,dx$ denote the average of $F$ over cube $Q$. Then
\be
\|f\circ T_{\Qc}\|_{L_p(\DOEP)}^p&\le& C\,\left(\sum_{Q\in\WQ}
|Q|\,|f(\omega_{Q,\alpha})-F(x_Q)|^p\right.
\nn\\
&+&\left.\sum_{Q\in\WQ}|Q|\,|F(x_Q)-F_Q|^p
+\sum_{Q\in\WQ}|Q|\,|F_Q|^p\right)= C\,(I+J+K).
\nn\ee
\par Let us consider the element  $\omega_{Q,\alpha}\in\DOA$ defined by formulas \rf{XNE} and \rf{DYO}. We recall that  $\omega_{Q,\alpha}$ is an $(\alpha,Q)$-visible element, see \rf{QVO}, and $\ell(\omega_{Q,\alpha})=a_Q$, i.e.,  $\ell(\omega_{Q,\alpha})$ is a point nearest to $Q$ on $\DO$ in the Euclidean norm. Hence,
$$
\|\ell(\omega_{Q,\alpha})-x_Q\|=\|a_Q-x_Q\|\le \dist(a_Q,Q)+r_Q =\dist(Q,\DO)+r_Q,
$$
so that, by \rf{TH-WTN},
$$
\|\ell(\omega_{Q,\alpha})-x_Q\|\le \theta\diam Q+r_Q= (2\theta+1)r_Q\le \eta\, r_Q.
$$
(Recall that $\eta\ge 22\theta^2$ and $\theta\ge 1$.) Thus $\ell(\omega_{Q,\alpha})=a_Q\in (\eta Q)\cap\DO.$
\par We put $q:=(n+p)/2$ and apply Lemma \reff{NG} to the cube $Q$ and the element $\omega_{Q,\alpha}$. We obtain
$$
\frac{|f(\omega_{Q,\alpha})-F(x_Q)|^p}{(\diam Q)^{p-n}}\le C\,\intl_{Q} \Mc[G]^{\frac{p}{q}}(z)\,dz,~~~Q\in\WQ.
$$
\par Recall that $G$ is a function defined by \rf{DGQ1}. Since $\diam Q\le \theta\ve$ for every $Q\in\WQ$, we have
\be
I&:=& \sum_{Q\in\WQ}
|Q|\,|f(\omega_{Q,\alpha})-F(x_Q)|^p\le
C\,\sum_{Q\in\WQ}(\diam Q)^p\,
\intl_{Q} \Mc[G]^{\frac{p}{q}}(z)\,dz\nn\\&\le& C\,(\theta\ve)^p\sum_{Q\in\WQ}\,
\intl_{Q} \Mc[G]^{\frac{p}{q}}(z)\,dz\le C\, \intl_{\Omega} \Mc[G]^{\frac{p}{q}}(z)\,dz= C\, \intl_{\RN} \Mc[G]^{\frac{p}{q}}(z)\,dz.
\nn\ee
\par Thus we have obtained the same estimate of $I$ as in inequality \rf{FEI}. Hence
$$
I\le C\, \|\nabla F\| ^p _{L_p(\Omega)},
$$
see \rf{E-MG}. Let us estimate the quantity
$$
J:=\sum_{Q\in\WQ}|Q|\,|F(x_Q)-F_Q|^p.
$$
By the Sobolev-Poincar\'e inequality \rf{SPE},
\be
|F(x_Q)-F_Q|^p\,|Q|&\le& |Q|\,\sup_{x,y\in Q} |F(x)-F(y)|^p\nn\\&\le& C\,(\diam Q)^p \intl_Q
\|\nabla F(z)\|^p\,dz\le C\, \intl_Q
\|\nabla F(z)\|^p\,dz,\nn
\ee
so that
$$
J\le C\, \sum_{Q\in\WQ}\,\intl_Q
\|\nabla F(z)\|^p\,dz.
$$
Since every two cubes of the family $\Qc$ are pairwise disjoint and $\WQ\subset \Qc$, we obtain
$$
J\le C\,\intl_{\Omega}\|\nabla F(z)\|^p\,dz=
C\,\|\nabla F\|^p_{L_p(\Omega)}.
$$
It remains to estimate the quantity
$K:=\sum\{|F_Q|^p\,|Q|:~Q\in\WQ\}.$ We have
$$
K\le\sum_{Q\in\WQ}
|Q|\left(\frac{1}{|Q|}\intl_{Q}|F|\,dz\right)^p\le \sum_{Q\in\WQ}\,\intl_{Q}|F|^p\,dz
$$
so that
$$
K\le
\intl_{\Omega}|F|^p\,dz=\|F\|^p_{L_p(\Omega)}.
$$
\par Summarizing the estimates for $I,J$ and $K$, we  finally obtain
$$
\|f\circ T_{\Qc}\|_{L_p(\DOEP)}^p\le C\,(I+J+K)\le C\,(\|\nabla F\|_{L_p(\Omega)}^p+
\|F\|_{L_p(\Omega)}^p)
\le C\|F\|^p_{\WPO}.
$$
\par The necessity part of the theorem is proved.
\medskip
\par {\it (Sufficiency.)} Let $f:\DOA\to\R$ be a continuous function with respect to the metric $\dco$ and let $f\circ T_{\Qc}\in L_{p}(\DOEP)$. Also assume that there exists  $\lambda>0$ such that
for every finite family $\{Q_{i}:i=1,...,m\}$ of pairwise disjoint cubes contained in $\DOEP$ and every $(\alpha,Q_i)$-visible elements $\omega_i^{(1)},\omega_i^{(2)}\in\DOA$ satisfying \rf{WP-THQ} the inequality \rf{WPV-IN} holds.
\par We put
\bel{D-DLW}
\sigma:=\sigma(\ve,\theta)=\ve/(80\theta),
\ee
and $\bc:=0$ in formula \rf{CQ}; thus $c_Q:=f(\oqa)$ if $Q\in\TW$ and  $\diam Q\le \sigma$ and $c_Q:=0$ if $\diam Q>\sigma$. Then we define a function $\tF:\Omega\to\R$ by formula \rf{ExtOp}:
$$
\tF(x):=\sum\limits_{Q\in \TW}
c_Q\varphi_Q(x)=\sum\limits_{Q\in \TW,\,\diam Q\le\sigma}
f(\oqa)\varphi_Q(x), ~~~ x\in \Omega.
$$
\par It can be easily seen that the extension $\tf$ defined by formula \rf{ExtOp-L} and the extension $\tF$ coincide in a $\sigma/2$-neighborhood of $\DO$. In fact, assume that  $\dist(x,\DO)< \sigma/2$. By Lemma \reff{PU}, (b), if $\varphi_Q(x)\ne 0$ then $Q^*=(9/8)Q\ni x$ so that
$$
\dist(Q,\DO)\le\dist(x,Q)+\dist(x,\DO)\le \diam Q/8 +\sigma/2.
$$
But, by \rf{DKQ}, $\diam Q\le \dist(Q,\DO)$ so that
$$
\diam Q\le\diam Q/8+\sigma/2,
$$
proving that $\diam Q<\sigma$.
\par Recall that in the sufficiency part of the proof of Theorem \reff{VS-GN} we have shown that $f=\tra \tf$. The proof of this equality relies only on inequality \rf{VAR-IN} which we apply to cubes $\{Q_i\}$ contained in a small neighborhood of $\DO$. For such cubes corresponding inequality of Theorem \reff{VCR-AO}, i.e., inequality \rf{WPV-IN}, holds as well. Since $\tF$ coincides with $\tf$ in a neighborhood of $\DO$ and $f=\tra \tf$, we conclude that $f=\tra \tF$.
\par Let us estimate the $\WPO$-norm of $\tF$. By Lemma
\reff{WP-NORM},
$$
\|\tF\|^p_{\WPO}\le C\,(I_1+I_2)
$$
where
$$
I_1:=\sum\left\{\frac {|f(\omega_{Q,\alpha})-f(\omega_{Q',\alpha})|^p}{(\diam Q+\diam Q')^{p-n}}:Q,Q'\in\TW, Q\cap Q'\ne\emptyset,\,Q,Q'\subset\Oc_{5\sigma}(\DO)\right\},
$$
and
\bel{I2-W}
I_2:=\sum\left\{|f(\omega_{Q,\alpha})|^p|Q|:~Q\in\TW, Q\subset\Oc_{10\sigma}(\DO)\right\}.
\ee
Since $5\sigma<\ve$, we obtain $I_1\le V_p(f;\DOEP)$ where
\be
V_p(f;\DOEP)&:=&\sum\{
|f(\omega_{Q,\alpha})-f(\omega_{Q',\alpha})|^p\,(\diam Q+\diam Q')^{n-p}:\nn\\
&&Q,Q'\in\TW, Q\cap Q'\ne\emptyset,
\,Q,Q'\in\DOEP\}.\nn
\ee
\par Observe that the definition of $V_p(f;\DOEP)$ is similar to that of the quantity $V_p(f;\Omega)$ where
the cubes $Q,Q'$ run over all cubes from $\TW$ such that $Q\cap Q'\ne\emptyset$, see \rf{VP-OM}. In turn,
the definition of $V_p(f;\DOEP)$ involves the same family of cubes with the additional requirement $Q,Q'\in\DOEP$.
\par This allows us to repeat the proof presented in the sufficiency part of Theorem \reff{VS-GN} and to show that an analog of inequality \rf{N-FW} holds for the quantity $V_p(f;\DOEP)$ as well. In other words,
$V_p(f;\DOEP)\le C\,\lambda$ proving that $I_1\le C\,\lambda$ where $C=C(n,p)$.
\par It remains to estimate the quantity $I_2$ defined by \rf{I2-W}. Let $Q\in\TW$ and let $Q\subset \Oc_{10\sigma}(\DO)$. Since $\Qc$ is a covering of $\Omega$, there exists a cube $K_Q\in\Qc$ such that $K_Q\cap Q\ne\emp$.
\par We let $S_Q$ denote a cube of diameter 
\bel{CN1}
\diam S_Q:=\min(\diam Q,K_Q)
\ee
such that
\bel{CN2}
S_Q\subset K_Q~~~\text{and}~~~S_Q\cap Q\ne\emp.
\ee
\par Let us compare diameters of $K_Q$ and $Q$. Let $y\in K_Q\cap Q$. Then
$$
\diam Q\le\dist(Q,\DO)\le\dist(y,\DO)
\le \diam K_Q+\dist(K_Q,\DO)
$$
so that, by \rf{TH-WTN},
\bel{DMQK}
\diam Q\le(1+\theta)\diam K_Q.
\ee
On the other hand,
\bel{A-E}
\dist(K_Q,\DO)\le\dist(y,\DO)
\le \diam Q+\dist(Q,\DO)
\ee
so that, by \rf{DKQ},
\bel{W-ER}
\dist(K_Q,\DO)\le 5\diam Q.
\ee
Note that \rf{A-E} and \rf{DKQ} also imply the following:
\bel{KF-E}
\dist(K_Q,\DO)\le \diam Q+\dist(Q,\DO)\le 2\dist(Q,\DO).
\ee
\par Now, by \rf{W-ER} and  \rf{TH-WTN},
$$
\diam K_Q\le\theta\dist(K_Q,\DO)\le 5\theta\diam Q
$$
proving that
\bel{KQ-C}
\frac{1}{1+\theta}\diam Q\le\diam K_Q\le 5\theta\diam Q.
\ee
\par Recall that $\diam S_Q=\diam Q$ provided $\diam Q\le \diam K_Q$ so that in this case
\bel{DKS}
\diam K_Q\le 5\theta\diam Q=5\theta\diam S_Q.
\ee
If $\diam K_Q\le \diam Q$, then $\diam S_Q=\diam K_Q$, so that \rf{DKS} holds as well.
\par Let us estimate the distance from $S_Q$ to the points
$a_Q$ and $a_{K_Q}$. By \rf{TH-WTN} and \rf{DKS},
\be
\|x_{S_Q}-a_{K_Q}\|&\le& \diam K_Q+\dist(K_Q,\DO)\le\diam K_Q+\theta\diam K_Q\nn\\
&\le&(1+\theta)(5\theta\diam S_Q)=10\theta(1+\theta)r_{S_Q}\nn
\ee
proving that $a_{K_Q}\in \gamma_1 S_Q$ where $\gamma_1:=10\theta(1+\theta).$
\par Recall that $Q\cap S_Q\ne\emp$ and $\diam S_Q=\min\{\diam Q,\diam K_Q\}$. Hence, by \rf{DMQK},
\bel{SQ-C}
\diam S_Q\le \diam Q\le(1+\theta)\diam S_Q.
\ee
\par By Lemma \reff{TR-SQ},
$\omega_{Q,\alpha}$ is $(\alpha,S_Q)$-visible in $\Omega$. Moreover,
$$
a_Q=\ell(\omega_{Q,\alpha})\in (10\tau +1)S_Q~~~\text{where}~~ \tau=(1+\theta).
$$
\par Observe that the element $\omega_{K_Q,\alpha}$ is $(\alpha,K_Q)$-visible, see \rf{QVO}. Since $S_Q\subset K_Q$, by Lemma \reff{INC}, $\omega_{K_Q,\alpha}$ is $(\alpha,S_Q)$-visible.
\par Summarizing the properties of the elements  $\omega_{Q,\alpha}$ and $\omega_{K_Q,\alpha}$, we conclude that $\omega_{Q,\alpha}$ and $\omega_{K_Q,\alpha}$ are $(\alpha,S_Q)$-visible and
$\ell(\omega_{Q,\alpha}),\ell(\omega_{K_Q,\alpha})
\in\eta S_Q$. (Recall that $\eta\ge 22\,\theta^2$.)
\par We also note that $Q\subset \Oc_{10\sigma}(\DO)$. Since $\diam S_Q \le \diam Q$, we have $S_Q\subset \Oc_{20\sigma}(\DO).$ Since $20\sigma<\ve$, see \rf{D-DLW}, we obtain $S_Q\subset \Oc_{\ve}(\DO)$, $Q\in\TW.$
\par Now can we estimate the quantity $I_2$, see \rf{I2-W}, as follows. Let $Q\in\TW$ and let  $Q\subset\Oc_{10\sigma}(\DO)$. We have
$$
|f(\omega_{Q,\alpha})|^p|Q|\le C(|f(\omega_{Q,\alpha})-f(\omega_{K_Q,\alpha})|^p|Q|+ |f(\omega_{K_Q,\alpha})|^p|Q|).
$$
\par Since $\diam Q\le \dist(Q,\DO)\le 10\sigma,$ by \rf{SQ-C},
$$
\diam S_Q\le\diam Q\le 10\sigma.
$$
Hence
\be
|f(\omega_{Q,\alpha})-f(\omega_{K_Q,\alpha})|^p|Q|&=&
\frac{|f(\omega_{Q,\alpha})-f(\omega_{K_Q,\alpha})|^p}{(\diam S_Q)^{p-n}}(\diam S_Q)^{p-n}|Q|\nn\\&\le& C(\sigma)\,\frac{|f(\omega_{Q,\alpha})
-f(\omega_{K_Q,\alpha})|^p}{(\diam S_Q)^{p-n}}.\nn
\ee
\par On the other hand, by \rf{KQ-C},
$$
|f(\omega_{K_Q,\alpha})|^p|Q|\le\,(1+\theta)^n\,
|f(\omega_{K_Q,\alpha})|^p\,|K_Q|.
$$
Hence,
$$
I_2:=\sum\left\{|f(\omega_{Q,\alpha})|^p|Q|:~Q\in\TW, Q\subset\Oc_{10\sigma}(\DO)\right\}\le C(J_1+J_2)
$$
where
$$
J_1:=\sum\left\{\frac{|f(\omega_{Q,\alpha})
-f(\omega_{K_Q,\alpha})|^p}{(\diam S_Q)^{p-n}}:~Q\in\TW, Q\subset\Oc_{10\sigma}(\DO)\right\},
$$
and
\bel{D-J2}
J_2:=\sum\left\{|f(\omega_{K_Q,\alpha})|^p\,|K_Q|:~Q\in\TW, Q\subset\Oc_{10\sigma}(\DO)\right\}.
\ee
\par Recall that  $Q\to S_Q$ is a mapping defined on the family of cubes
$$
\TW_{10\sigma}:=\{Q\in\TW:~Q\subset\Oc_{10\sigma}(\DO)\}
$$
and satisfying conditions \rf{CN1} and \rf{CN2}. Without loss of generality we may assume that this mapping $Q\to S_Q$ is one-to-one so that the converse mapping $S\to Q_S$ is well defined on the family of cubes
$$
\Sc:=\{S:~\exists\, Q\in\TW_{10\sigma}~~\text{such that}~~S=S_Q\}.
$$
\par Put $\omega^{(1)}_S:=\omega_{Q_S,\alpha}$ and $\omega^{(2)}_S:=\omega_{K_{Q_S},\alpha}$. Then
\bel{JF}
J_1=\sum_{S\in\Sc}\,\frac{|f(\omega^{(1)}_S)
-f(\omega^{(2)}_S)|^p}{(\diam S)^{p-n}}.
\ee
\par Since the family $\Qc$ is Whitney's type decomposition of $\Omega$, it has finite covering multiplicity $M(\Qc)\le N(n,\theta)$. But $S_Q\subset K_Q\in\Qc$ so that
$$
M(\Sc)\le M(\Qc)\le N(n,\theta).
$$
Consequently, the family $\Sc$ can be partitioned into at most $N_1(n,\theta)$ families of pairwise disjoint cubes, see \cite{BrK,Dol}. This allows us to assume that the family $\Sc$ itself consists of pairwise disjoint cubes.
\par Since the cubes $S\in \Sc$ and the elements $\omega^{(1)}_S,\omega^{(2)}_S$ from \rf{JF} satisfy the conditions of Theorem \reff{VCR-AO}, we can apply inequality \rf{WPV-IN} of this theorem to the quantity $J_1$. We obtain $J_1\le C\,\lambda.$
\par It remains to estimate the quantity $J_2$ defined by \rf{D-J2}. We let $\Kc$ denote a family of cubes 
$$
\Kc:=\left\{K_Q\in\Qc:~Q\in\TW, Q\subset\Oc_{10\sigma}(\DO) \right\}.
$$
Given $K\in\Kc$ we put
$$
G_K:=\left\{Q\in\TW:~Q\subset\Oc_{10\sigma}(\DO),K_Q=K \right\}.
$$
Since $\diam K_Q\sim\diam Q$, see \rf{KQ-C}, the family  $G_K$ consists of at most $N(n,\theta)$ cubes. Hence
$$
J_2=\sum_{K\in\,\Kc}\,\sum_{Q\in G_K}|f(\omega_{K,\alpha})|^p\,|K|\le C\,N(n,\theta)\sum_{K\in\,\Kc}
|f(\omega_{K,\alpha})|^p\,|K|.
$$
On the other hand, by definition \rf{DF-T} of the mapping $T_{\Qc}:\Omega\to\DOA$,
$$
|f(\omega_{K,\alpha})|^p\,|K|=\intl_{K}
|f\circ T_{\Qc}|^p(x)\,dx
$$
so that
$$
J_2\le C\,\sum_{K\in\,\Kc}\,\,
\intl_{K}|f\circ T_{\Qc}|^p(x)\,dx=C\,\intl_{U}|f\circ T_{\Qc}|^p(x)\,dx
$$
where $U:=\cup\{K:~K\in\Kc\}.$
\par Prove that $U\subset\Oc_{\ve}(\DO)$. Recall that for every $K\in\Kc$ there exists a cube $Q\in\TW$ such that $K=K_Q$. Thus  $Q\subset\Oc_{10\sigma}(\DO)$, $K\cap Q\ne\emp$,
$$
\dist(K,\DO)\le 2\dist(Q,\DO)~~\text{and}~~\diam K\le 5\theta \diam Q,
$$
see \rf{KF-E} and \rf{KQ-C}. Hence, for each $y\in K$ we have
\be
\dist(y,\DO)&\le&\dist(K,\DO)+\diam K\le 2\dist(Q,\DO)+5\theta\diam Q\nn\\&\le& 2\dist(Q,\DO)+5\theta\dist (Q,\DO) =(5\theta+2)\dist(Q,\DO).\nn
\ee
But $Q\subset\Oc_{10\sigma}(\DO)$ so that $\dist(Q,\DO)\le 10\sigma$. Hence, $\dist(y,\DO)\le(5\theta+2)10\sigma$
proving that
$$
K\subset \Oc_{\xi}(\DO)~~~\text{with}~~\xi:=10(5\theta+2)\sigma.
$$
Since $\theta\ge 1$, we have $\xi\le 70\theta\sigma$ so that $\xi<\ve$, see \rf{D-DLW}. Consequently, $K\subset\Oc_{\ve}(\DO)$ for every $K\in\Kc$ proving that $U:=\cup\{K:~K\in\Kc\}\subset\Oc_{\ve}(\DO)$.
\par Hence,
$$
J_2\le C\,\intl_{\Oc_{\ve}(\DO)}|f\circ T_{\Qc}|^p(x)\,dx=
C\,\|f\circ T_{\Qc}\|_{L_p(\DOEP)}^p.
$$
\par Finally, summarizing the estimates for the quantities $I_1$, $J_1$ and $J_2$, we obtain
$$
\|\tF\|^p_{\WPO}\le C\,(I_1+I_2)\le C\,(I_1+J_1+J_2)\le
C\,\left(\lambda+\lambda+\|f\circ T_{\Qc}\|_{L_p(\DOEP)}^p\right).
$$
\par Theorem \reff{VCR-AO} is completely proved.\bx
\section*{{\normalsize 6. Sharp maximal functions on the Sobolev boundary of a domain.}}
\setcounter{section}{6}
\setcounter{theorem}{0}
\setcounter{equation}{0}
\par We turn to the last result of the paper, Theorem \reff{MF-C}, which is a generalization of Theorem \reff{MF-OM} formulated in Section 1.
\par Fix $q\in(n,p)$ and put $\beta:=\qn$ and $\alpha:=\pn$. By $\dlbeta$ we denote a quasi-metric on $\OB=\Omega\cup\DOB,$ see \rf{B-A},
defined by the formula
$$
\dlbeta(\omega_1,\omega_2)
:=d^{\frac{1}{\beta}}_{\beta,\clO}(\omega_1,\omega_2),
~~~\omega_1,\omega_2\in\OB.
$$
\par Given $x\in\Omega$ and $r>0$ by $B(x,r:\dlbeta)$ we denote the closed ball in the quasi-metric space $(\OB,\dlbeta)$ with center $x$ and radius $r$:
\bel{B-DL}
B(x,r:\dlbeta):=\{\omega\in\OB:~
\dlbeta(x,\omega)\le r\}.
\ee
\par Let $\omega=[(x_i)]_\beta\in\OB$. Since $0<\beta<\alpha$, by Corollary \reff{AB-CH},
$$
\dt(x,y)\le C \db(x,y),~~~x,y,\in \Omega.
$$
Consequently, every Cauchy sequence $(x_i)$ with respect to $\db$ is a Cauchy sequence with respect to $\dt$. Moreover, by definition \rf{ESIM},
$$
(x_i)\sbs (y_i)~~ \Longrightarrow~~ (x_i)\sa (y_i).
$$
\par Thus for every $(x_i),(y_i)\in\omega$ we have
$[(x_i)]_\alpha=[(y_i)]_\alpha$
so that all sequences of the equivalence  class $\omega$ (with respect to $``\sbs"$) belong to the same equivalence class with respect to $``\sa"$. We denote this equivalence class by $\omega^{[\alpha]}$; thus
\bel{D-OAH}
\omega^{[\alpha]}:=[(x_i)]_\alpha\in \OA,~~~(x_i)\in\omega.
\ee
Moreover, if $\omega\in\DOB$, then $\omega^{[\alpha]}\in\DOA$ which shows that it is defined a mapping
$$
\DOB\ni\omega\to\omega^{[\alpha]}\in\DOA.
$$
(Observe that in general $\omega\varsubsetneq\omega^{[\alpha]}$ (as families of sequences) while $\ell(\omega)=\ell(\omega^{[\alpha]})$.)
\par This enables us given a function $f:\DOA\to\R$ to define its fractional sharp maximal function $f_{\infty,\beta,\Omega}^{\sharp}$ on $\Omega$ as follows: for every $x\in\Omega$ we put
\bel{SH-MF}
f_{\infty,\beta,\Omega}^{\sharp }(x):=
\sup\left\{\frac{
|f(\omega_1^{[\alpha]})-f(\omega_2^{[\alpha]})|}{r}:r>0,\,
\omega_1,\omega_2\in  B(x,r:\dlbeta)\cap\DOB\right\}.
\ee
\begin{theorem}\lbl{MF-C} Let $\Omega$ be a domain in $\RN$ and let $n<q<p$, $\beta=\qn$ and $\alpha=\pn$.
\par (i). A function $f\in\tra(\LOPO)$ if and only if
$f$ is continuous on $\DOA$ (with respect to $\dco$) and
$f_{\infty,\beta,\Omega}^{\sharp}\in L_p(\Omega)$.
Moreover,
\bel{E-LG}
\|f\|_{\tra(\LOPO)}\sim
\|f_{\infty,\beta,\Omega}^{\sharp }\|_{L_p(\Omega)}.
\ee
\par (ii). Fix $\ve>0,\theta>1$ and a covering $\Qc$ of $\Omega$ consisting of non-overlapping cubes $Q\subset\Omega$ satisfying inequality \rf{TH-WTN}.
Let $T_{\Qc}$ be the mapping defined by \rf{DF-T}.
\par A function $f:\DOA\to\R$ is an element of \,\,
$\tra(\WPO)$ if and only if  $f$ is continuous on $\DOA$ (with respect to $\dco$), and
$$
f\circ T_{\Qc}~~\text{and}~~
f_{\infty,\beta,\Omega}^{\sharp }~~\text{are both in}~~
L_{p}(\DOEP).
$$
Furthermore,
\bel{E-WG}
\|f\|_{\tra(\WPO)}\sim \|f\circ T_{\Qc}\|_{L_p(\DOEP)}+ \|f_{\infty,\beta,\Omega}^{\sharp }\|_{L_p(\DOEP)}.
\ee
\par The constants of equivalence in \rf{E-LG} depend only
on $n,p$ and $q$, and in \rf{E-WG} they only depend on
$n,p,q,\ve$ and $\theta$.
\end{theorem}
\par {\it Proof. (Necessity)}. (i). Let $F\in\LOPO$ and let $f=\tra F.$
\par Recall that, by \rf{TR-DLM}, for every $\omega\in\DOA$
\bel{DG-T}
f(\omega)=\lim\{F(x):\dco(x,\omega)\to 0, x\in\Omega\}.
\ee
\par Let
$$
G(x):=(\|\nabla F\|^q)^\cw(x), ~~~~~x\in\RN.
$$
Prove that
\bel{E-SM}
f_{\infty,\beta,\Omega}^{\sharp }(x)\le C(n,q)\,(\Mc[G])^{\frac{1}{q}}(x),~~~x\in\Omega.
\ee
(Recall that $\Mc$ stands for the Hardy-Littlewood maximal operator, see \rf{H-L-M}.)
\par Let $\omega_1,\omega_2\in\DOB,$ $\omega_1\ne\omega_2$, and let $\omega_1,\omega_2\in  B(x,r:\dlbeta).$ Hence
\bel{DBC}
\dcb(x,\omega_1),\,\dcb(x,\omega_2)\le r^\beta.
\ee
\par Since $\omega_1,\omega_2\in\DOB$, there exist  sequences $(x_i)\in\omega_1,$ $(y_i)\in\omega_2$  such that
\bel{D1}
\dtcb(x_i,\omega_1)\to 0~~\text{and}~~\dtcb(y_i,\omega_2)\to 0~~\text{as}~~i\to\infty.
\ee
Recall that $\omega_1^{[\alpha]}=[(x_i)]_\alpha$ and $\omega_2^{[\alpha]}=[(y_i)]_\alpha$ so that
$\dtc(x_i,\omega^{[\alpha]}_1)$ and  $\dtc(y_i,\omega^{[\alpha]}_2)$ tend to $0$  as $i\to\infty$.
Consequently, by \rf{DG-T},
\bel{FLM1}
f(\omega^{[\alpha]}_1)=\lim_{i\to\infty}F(x_i)~~\text{and}~~
f(\omega^{[\alpha]}_2)=\lim_{i\to\infty}F(y_i).
\ee
\par By \rf{D1},
$$
\db(x_i,y_i)=\dtcb(x_i,y_i)\to\dtcb(\omega_1,\omega_2)
~~\text{as}~~i\to\infty,
$$
so that there exists $N_0\in\N$ such that
$$
\db(x_i,y_i)\le 2\dtcb(\omega_1,\omega_2), ~~~i\ge N_0.
$$
Combining this inequality with \rf{U-0}, we obtain
$$
\db(x_i,y_i)\le C\,\dcb(\omega_1,\omega_2), ~~~i\ge N_0,
$$
so that, by \rf{DBC},
$$
\db(x_i,y_i)\le C\,(\dcb(\omega_1,x)+\dcb(x,\omega_2))\le C\,(r^\beta+r^\beta)\le C\,r^\beta, ~~~i\ge N_0.
$$
\par Let $\lambda_1=\lambda_1(n,q)$ be the constant from Proposition \reff{SP-Q}. Then
$$
\lambda_1\,\db(x_i,y_i)^{\frac1\beta}\le \lambda_1C^{\frac1\beta}\,r, ~~~i\ge N_0.
$$
By this proposition,
\bel{FXY}
|F(x_i)-F(y_i)|\le C\,R\left(\frac{1}{|Q(x_i,R)|} \intl_{Q(x_i,R)}
G(z)\,dz\right)^{\frac{1}{q}}, ~~~i\ge N_0,
\ee
provided $R:=\lambda_1 C^{\frac1\beta}\,r$. On the other hand,
$$
\|x_i-x\|^\beta\le\db(x_i,x)\le\dtcb(x_i,\omega_1)
+\dtcb(x,\omega_1)
$$
so that, by \rf{D1}, there exists $N_1\in\N$ such that
$$
\|x_i-x\|^\beta\le 2\,\dtcb(x,\omega_1),~~~i\ge N_1.
$$
Combining this inequality with \rf{EQ-GR1} and \rf{DBC}, we obtain
$$
\|x_i-x\|^\beta\le C\,\dcb(x,\omega_1)\le Cr^\beta,~~~i\ge N_1,
$$
proving that
$$
\|x_i-x\|\le Cr\le\gamma R,~~~i\ge N_1,
$$
with $\gamma=\gamma(n,q,\beta)$. By this inequality,
$$
Q(x_i,R)\subset\tQ:=Q(x,(\gamma+1)R)
$$
provided $i\ge N_1$. In addition, $|Q(x_i,R)|\sim |\tQ|$, so that, by \rf{FXY},
$$
|F(x_i)-F(y_i)|\le C\,R\left(\frac{1}{|\tQ|} \intl_{\tQ}
G(z)\,dz\right)^{\frac{1}{q}}, ~~~i\ge N_2:=\max\{N_0,N_1\}.
$$
\par Now, letting $i$ tend to $\infty$, by \rf{FLM1}, we have
$$
|f(\omega^{[\alpha]}_1)-f(\omega^{[\alpha]}_2)|\le C\,R\left(\frac{1}{|\tQ|} \intl_{\tQ}
G(z)\,dz\right)^{\frac{1}{q}}\le C\,R\left(\Mc[G](x)\right)^{\frac{1}{q}}\le C\,r\left(\Mc[G](x)\right)^{\frac{1}{q}}
$$
so that
$$
|f(\omega^{[\alpha]}_1)-f(\omega^{[\alpha]}_2)|/r\le C\,\left(\Mc[G](x)\right)^{\frac{1}{q}}.
$$
Taking the supremum in the left-hand side of this inequality over all $\omega_1,\omega_2\in\DOB$ satisfying \rf{DBC}, we obtain inequality \rf{E-SM}.
\par By this inequality,
$$
\|f_{\infty,\beta,\Omega}^{\sharp }\|_{L_p(\Omega)}\le C\,\|(\Mc[G])^{\frac{1}{q}}\|_{L_p(\Omega)}\le C\,\|(\Mc[G])^{\frac{1}{q}}\|_{L_p(\RN)}.
$$
But, by \rf{E-MG},
$$
\|(\Mc[G])^{\frac{1}{q}}\|_{L_p(\RN)}\le C\,\|\nabla F\|_{L_p(\Omega)}
$$
so that $\|f_{\infty,\beta,\Omega}^{\sharp }\| _{L_p(\Omega)} \le C\,\|\nabla F\|_{L_p(\Omega)}.$
Taking the infimum over all functions $F\in\LOPO$ such that $\tra F=f$, we obtain the required inequality
$$
\|f_{\infty,\beta,\Omega}^{\sharp }\|_{L_p(\Omega)}\le C\,\|f\|_{\tra(\LOPO)}.
$$
\par (ii). The latter inequality yields
$$
\|f_{\infty,\beta,\Omega}^{\sharp }\|_{L_p(\DOEP)}\le\|f_{\infty,\beta,\Omega}^{\sharp }\|_{L_p(\Omega)}\le C\,\|f\|_{\tra(\LOPO)}\le C\,\|f\|_{\tra(\WPO)}.
$$
\par In turn, by Theorem \reff{VCR-AO},
$$
\|f\circ T_{\Qc}\|_{L_p(\DOEP)}\le C\,\|f\|_{\tra(\WPO)},
$$
see \rf{NW-LO}. These two inequalities show that the right-hand side of equivalence \rf{E-WG} is bounded by $C\,\|f\|_{\tra(\WPO)}$.
\par The necessity part of the statements (i) and (ii) is proved.
\bigskip
\par {\it (Sufficiency).}
\par (i). Fix a constant $\tau\ge 1$. Let $Q\subset\Omega$ be a cube and let $\omega\in\DOA$ be an $(\alpha,Q)$-visible element such that $\ell(\omega)\in (\tau Q)\cap \DO.$
\par Let $(y_i)$ be a sequence of points in $\Omega$ such that $y_i\in(\ell(\omega),x_Q]$, $i=1,2,...~,$ and
$\omega=[(y_i)]_\alpha,$ see Definition \reff{Q-VA}. Recall that $\ell(\omega)$ is a $Q$-visible point, the line segment $(\ell(\omega),x_Q]\subset\Omega$, see Definition \reff{Q-V}, and $\lim_{i\to\infty}y_i=\ell(\omega).$
\par By Lemma \reff{DS-QV}, $(y_i)$ is a Cauchy sequence with respect to the metric $\db$. We put
$$
\tom=[(y_i)]_\beta.
$$
\par Observe that, by Lemma \reff{DS-QV}, $\tom$ is well defined and does not depend on the choice of the sequence $(y_i)\in\omega$ such that  $y_i\in(\ell(\omega),x_Q]$, $i=1,2,...~$. Also note that $\tom\in\DOB$. Since $\omega=[(y_i)]_\alpha$, by \rf{D-OAH},
$\omega=\tom^{[\alpha]}.$
\par Let us estimate the distance $\dcb(\tom,x_Q)$.
By Lemma \reff{DS-QV}, (i),
$$
\dbo(y_i,x_Q)\le\intl_{[y_i,x_Q]}
\D{z}^{\beta-1}\,ds(z)
\le C(\beta)\left(\frac{\|\ell(\omega)-x_Q\|}{\diam Q}\right)^{1-\beta}\|y_i-x_Q\|^\beta.
$$
Recall that $\ell(\omega)\in\tau Q$.
Since $y_i\in(\ell(\omega),x_Q]$, the point $y_i\in\tau Q$
as well so that
$$
\|\ell(\omega)-x_Q\|\le \tau r_Q=\tau\diam Q/2,
$$
and
$$
\|y_i-x_Q\|\le \tau r_Q=\tau\diam Q/2,~~~i=1,2,...~.
$$
Hence
$$
\dbo(y_i,x_Q)\le  C(\diam Q)^\beta,
$$
where $C=C(\beta,\tau).$ Since $\tom=[(y_i)]_\beta$, we have
$$
\dcb(\tom,x_Q)=\lim_{i\to\infty}\dbo(y_i,x_Q)\le C(\diam Q)^{\beta}.
$$
\par Now, let $y$ be an arbitrary point in $Q$.
Then
$$
\|x_Q-y\|\le r_Q\le\dist(x_Q,\DO)\le\max\{\dist(x_Q,\DO)\dist(y,\DO)\},
$$
so that, by Lemma \reff{Q-OM},
$$
\dbo(x_Q,y)\le\intl_{[x_Q,y]}
\D{z}^{\beta-1}\,ds(z)
\le \tfrac1\beta\,\|x_Q-y\|^\beta\le\tfrac1\beta\,(\diam Q/2)^\beta.
$$
Hence
$$
\dcb(\tom,y)\le\dcb(\tom,x_Q)+\dbo(x_Q,y)\le C\,(\diam Q)^\beta,
$$
proving that
$$
\dlbeta(\tom,y):=\dcb^{\frac1\beta}(\tom,y)\le R:=C^{\frac1\beta}\diam Q.
$$
Thus
\bel{B-DO}
\tom\in B(y,R:\dlbeta)~~\text{where}~~R=C^{\frac1\beta}\diam Q,
\ee
see \rf{B-DL}.
\par Let $\omega_1,\omega_2\in\DOA$ be two $(\alpha,Q)$-visible elements such that
$$
\ell(\omega_1),\ell(\omega_2)\in(\tau Q)\cap(\DO).
$$
\par Prove that
\bel{M-FQ}
I:=\frac{|f(\omega_1)-f(\omega_2)|^p}
{(\diam Q)^{p-n}}
\le C\,|Q|\,(f_{\infty,\beta,\Omega}^{\sharp})^p(y).
\ee
\par In fact, $\omega_1=\tom_1^{[\alpha]}$, $\omega_2=\tom_2^{[\alpha]}$, and, by \rf{B-DO},
\bel{OM12}
\tom_1,\tom_2\in B(y,R:\dlbeta)~~\text{with}~~R=C^{\frac1\beta}\diam Q.
\ee
Hence
\be
I&:=&\frac{|f(\omega_1)-f(\omega_2)|^p}
{(\diam Q)^{p-n}}=
\frac{|f(\tom_1^{[\alpha]})-f(\tom_2^{[\alpha]})|^p}
{(\diam Q)^{p-n}}\nn\\&\le& C\,|Q|\left(\frac{|f(\tom_1^{[\alpha]})-f(\tom_2^{[\alpha]})|}
{\diam Q}\right)^p\le C\,|Q|\left(\frac{|f(\tom_1^{[\alpha]})-f(\tom_2^{[\alpha]})|}
{R}\right)^p.\nn
\ee
\par By \rf{OM12} and \rf{SH-MF},
$$
|f(\tom_1^{[\alpha]})-f(\tom_2^{[\alpha]})|/R\le
f_{\infty,\beta,\Omega}^{\sharp}(y),
$$
and inequality \rf{M-FQ} follows.
\par Integrating this inequality over cube $Q$, we obtain
\bel{IN-MF}
I:=\frac{|f(\omega_1)-f(\omega_2)|^p}
{(\diam Q)^{p-n}}\le C\intl_{Q}(f_{\infty,\beta,\Omega}^{\sharp})^p(y)\,dy
\ee
where $C=C(n,p,\beta,\tau).$
\par Let $\Ac=\{Q_i:~i=1,...,m\}$ be a finite family of pairwise disjoint cubes in $\Omega$. Let $\omega_i^{(1)},\omega_i^{(2)}\in\DOA$
be $(\alpha,Q_{i})$-visible elements such that
$$
\ell(\omega_i^{(1)}),\ell(\omega_i^{(2)})\in (\tau Q)\cap \DO,~~~i=1,...,m,
$$
where $\tau=41$. Then, by \rf{IN-MF},
$$
J(f;\Ac):=
\sum_{i=1}^{m}\frac{|f(\omega_i^{(1)})-f(\omega_i^{(2)})|^p}
{(\diam Q_{i})^{p-n}}\le C\,\sum_{i=1}^{m}
\intl_{Q_i}(f_{\infty,\beta,\Omega}^{\sharp})^p(y)\,dy.
$$
Since the cubes of the family $\Ac$ are pairwise disjoint, we have
\bel{Z1}
J(f;\Ac)\le C\,\|f_{\infty,\beta,\Omega}^{\sharp}\|^p_{L_p(U)}
\ee
where $U:=\cup\{Q_i:~i=1,...,m\}$ and $C=C(n,p,\beta)$. Hence $J(f;\Ac)\le \lambda$ provided $\lambda:=C\,\|f_{\infty,\beta,\Omega}
^{\sharp}\|^p_{L_p(\Omega)}$.
\par Now, let $f:\DOA\to\R$ be a continuous  function (with respect to the metric $\dco$) and let $f_{\infty,\beta,\Omega}^{\sharp}\in L_p(\Omega)$. Then, by Theorem \reff{VS-GN}, $f\in\tra(\LOPO)$ and the following inequality
$$
\|f\|_{\tra(\LOPO)}\le C\,\lambda^{\frac1p}\le C\,
\|f_{\infty,\beta,\Omega}^{\sharp }\|_{L_p(\Omega)}
$$
holds.
\par The sufficiency part of the statement (i) of Theorem \reff{MF-C} is proved.
\bigskip
\par (ii). By inequality \rf{Z1} with $\tau=22\theta^2$, we have
$$
J(f;\Ac)\le C\,\|f_{\infty,\beta,\Omega}^{\sharp}\|^p_{L_p(\DOEP)}
$$
provided $\Ac=\{Q_i:~i=1,...,m\}$ is a finite family of pairwise disjoint cubes contained in $\DOEP$. Hence $J(f;\Ac)\le\lambda$ where $\lambda:=C\,\|f_{\infty,\beta,\Omega}
^{\sharp}\|^p_{L_p(\DOEP)}$.
\par Consequently, if $f:\DOA\to\R$ is a continuous  function (with respect to the metric $\dco$), and the functions $f\circ T_{\Qc}$ and $f_{\infty,\beta,\Omega}^{\sharp}$ are both in $L_{p}(\DOEP)$, then, by Theorem \reff{VCR-AO}, the function $f\in\tra(\WPO)$. Moreover,
$$
\|f\|_{\tra(\WPO)}\le C(\|f\circ T_{\Qc}\|_{L_p(\DOEP)}+\lambda^{\frac1p})\le C(\|f\circ T_{\Qc}\|_{L_p(\DOEP)}+ \|f_{\infty,\beta,\Omega}^{\sharp }\|_{L_p(\DOEP)}).
$$
\par Theorem \reff{MF-C} is completely proved.\bx
\par {\bf Acknowledgement.} I am very thankful to M. Cwikel and V. Maz'ya for useful suggestions and remarks. I am also very grateful to C. Fefferman, N. Zobin and all participants of "Whitney Problems Workshop", College of William and Mary, August 2009, for stimulating discussions
and valuable advice.


\begin{thebibliography}{ABCD}
\bibitem [1]{AHHL} K. Astala, K. Hag, P. Hag, V.
    Lappalainen, Lipschitz Classes and the Hardy-Littlewood
    Property, Mh. Math 115 (1993) 267--279.
\bibitem [2] {AS} N. Aronszajn, P. Szeptycki,
Theory of Bessel potentials. IV,  Ann. Inst. Fourier, 25
(1975) 27--69.
\bibitem [3] {Be} O. V. Besov,  The behavior of
differentiable functions on a non-smooth surface.
Studies in the theory of differentiable functions
of several variables and its applications, IV.
Trudy Mat. Inst. Steklov. 117 (1972) 3--10. (Russian)
\bibitem [4]{Br1} Yu. A. Brudnyi, Spaces that are
    definable by means of local approximations, Trudy
    Moscov. Math. Obshch., 24 (1971) 69--132; English
    transl. in Trans. Moscow Math. Soc. 24 (1974) 73--139.
\bibitem [5]{BrK} Yu. A. Brudnyi, B. D. Kotljar, A certain
    problem of combinatorial geometry,  Sibirsk. Mat. Z. 11        (1970) 1171–-1173; English transl. in Siberian Math. J.    11 (1970) 870–-871.
\bibitem [6]{BKos} S. Buckley, P. Koskela, Criteria for
    imbeddings of Sobolev-Poincar\'{e} type, Internat.
    Math. Res. Notices 18 (1996) 881--902.
\bibitem [7]{BSt2} S. Buckley, A. Stanoyevitch, Weak
    slice conditions and H\"older
    imbeddings, J. London Math. Soc. 66 (2001) 690--706.
\bibitem [8]{BSt3} S. Buckley, A. Stanoyevitch, Weak
    slice conditions, product domains, and quasiconformal
    mappings, Rev. Math. Iberoam. 17 (2001) 1--37.
\bibitem [9]{C1} A. P. Calder\'{o}n, Estimates for
    singular integral operators in terms of maximal
    functions, Studia Math. 44 (1972) 563--582.
\bibitem [10]{CS}  A. P. Calder\'{o}n, R. Scott, Sobolev
    type inequalities for $p>0$,  Studia Math. 62 (1978)
    75--92.
\bibitem [11]{Dol} V. L. Dolnikov, The partitioning of
    families of convex bodies, Sibirsk. Mat. Z.
    12 (1971) 664-–667 (Russian); English transl. in
    Siberian Math. J. 12 (1971) 473-–475.
\bibitem [12]{Ga} E. Gagliardo, Caratterizzazioni delle
tracce sulla frontiera relative ad alcune classi di
funzioni in $n$ variabili, Rend. Sem. Mat. Univ. Padova
27 (1957) 284--305.(Italian)
\bibitem [13]{GM} F.W. Gehring, O. Martio, Lipschitz
    classes and quasiconformal mappings, Ann. Acad. Sci.
    Fenn. Ser. AI Math. 10 (1985) 203–-219.
\bibitem [14]{Gr} P. Grisvard,  Elliptic problems in
    non-smooth domains. Monographs and Studies in
    Mathematics, 24. Pitman, Boston, 1985. xiv+410 pp.
\bibitem [15]{G} M. de Guzm\'{a}n, Differentiation of
    integrals in $\RN$, Lect. Notes in Math. 481,
    Springer-Verlag, 1975.
\bibitem [16]{Jn} P. W. Jones, Quasiconformal mappings and
    extendability of functions in Sobolev spaces, Acta
    Math. 147 (1981) 71--78.
\bibitem [17]{J1}  A. Jonsson, The trace of potentials on
    general sets, Ark. Mat. 17 (1979) 1--18.
\bibitem [18]{J2} A. Jonsson, Besov spaces on surfaces
    with singularities, Manuscripta Math. 71 (1991)
    121--152.
\bibitem [19]{JW} A. Jonsson, H. Wallin, Function Spaces on Subsets of $\RN$, Harwood Acad. Publ., London, 1984,
    Mathematical Reports, Volume 2, Part 1.
\bibitem [20]{L}  V. Lappalainen, $Lip_h$-extension
    domains, Ann. Acad. Sci. Fenn. Ser. A I Math.
    Dissertations 56 (1985) 1-52.
\bibitem [21]{LM} J-L. Lions, E. Magenes, Non-homogeneous
    boundary value problems and applications, Berlin, New
    York, Springer-Verlag, 1972.
\bibitem [22]{M}  V.G. Maz'ja,  Sobolev spaces, Springer-Verlag, Berlin, 1985, xix+486 pp.
\bibitem [23]{MP1} V. Maz'ya,  S. Poborchi, Traces of functions from S. L. Sobolev spaces on small and large components of the boundary,  Mat. Zametki 45 (1989), no. 4, 69--77, 126; Engl. transl. in Math. Notes 45 (1989), no. 3-4, 312--317.
\bibitem [24]{MP2} V. Maz'ya,  S. Poborchi, Traces of functions with a summable gradient in a domain with a cusp at the boundary, Mat. Zametki 45 (1989), no. 1, 57--65, 140; Engl. transl. in Math. Notes 45 (1989), no. 1-2, 39--44
\bibitem [25]{MP3} V. Maz'ya,  S. Poborchi, Boundary traces of functions from Sobolev spaces on a domain with a cusp, Trudy Inst. Mat. (Novosibirsk) 14 (1989) 182--208; Engl. transl. in Siberian Adv. Math. 1 (1991), no. 3, 75--107.
\bibitem [26]{MP} V. Maz'ya,  S. Poborchi, Differentiable
Functions on Bad Domains, Word Scientific, River Edge, NJ, 1997.
\bibitem [27]{MPN} V. Maz'ya,  S. Poborchi, Yu. Netrusov, Boundary values of functions from Sobolev spaces in some non-Lipschitzian domains. (Russian) Algebra i Analiz 11 (1999), no. 1, 141--170; Eng. transl. in St. Petersburg Math. J. 11 (2000), no. 1, 107--128.
\bibitem [28]{N} S. M. Nikol'skii,  Boundary properties of
functions defined on a region with angular points. I--III.
Mat. Sb. 40 (1956) 303--318, 43 (1957) 127--144, 45 (1958)
181--194. (Russian); Eng. transl.in Amer. Math. Soc.
Transl., Series 2, 83 (1958) 101--158.
\bibitem [29]{P1} S. Poborchi, Continuity of the boundary
trace operator  $W\sp 1\sb p(\Omega)\to L\sb Q(\DO)$ for a  domain with outward peak. Vestnik St. Petersburg Univ.
Math. 38 (2005), no. 3, 37--44 (2006).
\bibitem [30]{S3} P. Shvartsman, Sobolev $W^1_p$-spaces
on closed subsets of $\RN$, Advances in Math. 220 (2009) 1842–-1922.
\bibitem [31]{S4} P. Shvartsman, On Sobolev extension
domains in $\RN$, J. Funct. Anal. 258 (2010) 2205–-2245.
\bibitem [32]{St} E. M. Stein, Singular integrals and
differentiability properties of functions, Princeton
Univ. Press, Princeton, New Jersey, 1970.
\bibitem [33]{V} M. Yu. Vasil'chik,  The boundary behavior
of functions of Sobolev spaces defined on a planar domain
with a peak vertex on the boundary [Translation of Mat.
Tr. 6 (2003), no. 1, 3--27; MR1985623]. Siberian Adv.
Math. 14 (2004), no. 2, 92--115.
\bibitem [34]{Y1} G. N. Yakovlev, Boundary properties of
functions of the class $W\sb{p}{}\sp{(l)}$ in regions with
corners. (Russian) Dokl. Akad. Nauk SSSR 140 (1961)
73--76.; Eng. transl. in Soviet Math. 2 (1961) 1177--1180.
\bibitem [35]{Y2} G. N. Yakovlev, Dirichlet problem for a
region with a non-Lipschitz boundary.
Differencial'nye Uravnenija 1 (1965) 1085--1098 (Russian); Eng. transl. in Differential Equations 1 (1965), 847–-858.
\end{thebibliography}
\end{document}